\pdfoutput=1
\documentclass[11pt]{article}

\usepackage{amsmath, amssymb, amsthm, mathtools,authblk}
\usepackage[margin=1in]{geometry}
\usepackage{hyperref}
\usepackage[round]{natbib}
\usepackage{enumitem}
\usepackage{graphicx}
\usepackage{xcolor}
\usepackage{booktabs}

\mathtoolsset{showonlyrefs}

\newtheorem{theorem}{Theorem}[section]
\newtheorem{lemma}[theorem]{Lemma}
\newtheorem{proposition}[theorem]{Proposition}
\newtheorem{corollary}[theorem]{Corollary}
\theoremstyle{remark}
\newtheorem{remark}[theorem]{Remark}
\theoremstyle{definition}
\newtheorem{assumption}{Assumption}

\DeclareMathOperator{\tr}{tr}
\DeclareMathOperator{\Ran}{Ran}
\newcommand{\E}{\mathbb{E}}
\newcommand{\B}{\mathcal{B}}

\newcommand{\cH}{\mathcal{H}}

\newcommand{\ipH}[2]{\langle #1,\, #2 \rangle_{\cH}}
\newcommand{\Id}{\boldsymbol{I}}

\begin{document}
\title{Regularization vs.~Convolution in Kernel-Based Hypothesis Testing}

\author[1]{Omar Hagrass}
\author[2]{Bharath K. Sriperumbudur}
\affil[1]{Department of Statistics, University of Oxford\protect\\
\texttt{omar.hagrass@stats.ox.ac.uk}}
\affil[2]{Department of Statistics, Pennsylvania State University\protect\\
\texttt{bks18@psu.edu}}
\date{}
\maketitle

\begin{abstract}

The maximum mean discrepancy (MMD) test has become a popular nonparametric
testing approach that is widely used on general (i.e., non-Euclidean) domains.
Two distinct lines of work have investigated the optimality of the MMD test
through structurally different approaches. The first, a \emph{convolution}
approach on $\mathbb{R}^d$, uses translation-invariant kernels with a shrinking
bandwidth and is minimax optimal over Sobolev balls under $L^2$-separation
\citep{li2019optimality, schrab2023mmd}; the second, a \emph{regularization}
approach on general measurable spaces, uses spectral regularization of the
kernel covariance operator and is minimax optimal over Hellinger- and
$\chi^2$-separated alternatives defined through the range space of a fractional
power of the integral operator \citep{hagrass2023Gof, hagrass2024a_spectral_twosample}. In this paper, we study the optimality of the regularization approach over Sobolev balls under $L^2$-separation, which reveals the connection between the two approaches. We show that, with a translation-invariant kernel on $\mathbb{R}^d$, both are ways of meeting a single condition on the spectral density of the kernel. Specifically, the spectral regularized discrepancy captures the $L^2$-separation $\Delta$ of a Sobolev alternative of smoothness $s$ as soon as the regularization parameter $\lambda$ is below the kernel's spectral density on frequencies up to $\Delta^{-1/(2s)}$. A kernel with a shrinking bandwidth $h$ meets this condition with $\lambda \asymp h^d$. Under this coupling, the constant regularizer, which gives the MMD test itself, attains the Sobolev minimax rate. Therefore, this result explains the need for a shrinking bandwidth in the convolution approach from a spectral perspective, as a way of meeting the required spectral condition with a constant regularizer. Interestingly, with a non-constant spectral regularizer, our analysis requires the bandwidth to shrink only when the spectral density decays faster than any polynomial, as for the Gaussian kernel. Thus, for kernels whose spectral density decays polynomially, such as the Mat\'ern and Laplace kernels, a fixed bandwidth with $\lambda \to 0$ attains the Sobolev minimax rate. Moreover, we show that an adaptive spectral regularized test aggregating over a grid of bandwidths and regularization parameters is minimax optimal up to a $\log\log$ factor over both classes without prior knowledge of the smoothness.

\end{abstract}
\section{Introduction}

Given independent samples $\mathbb{X}_n := (X_i)_{i=1}^n
\overset{\mathrm{iid}}{\sim} P$ and $\mathbb{Y}_m :=
(Y_j)_{j=1}^m \overset{\mathrm{iid}}{\sim} Q$ on a measurable space
$\mathcal{X}$, the two-sample testing problem is to test
$H_0: P = Q$ against $H_1: P \neq Q$. A closely related setting is
the goodness-of-fit testing problem, in which $Q = P_0$ is a known
reference distribution, where the test has to be performed based on only $\mathbb{X}_n$. Both problems are classical in statistics
and have been studied extensively in the parametric and
nonparametric literature \citep{lehmann2005testing}. Over the last
decade, tests based on the maximum mean discrepancy (MMD)---a
kernel-based distance between probability measures---have become a
method of choice for nonparametric testing on general (i.e.,
non-Euclidean) domains due to their flexibility and computational
tractability \citep{NipsGretton,gretton2012kernel}.

The performance of a test against shrinking alternatives is
quantified through the framework of minimax separation rates
\citep{Burnashev,ingster1987minimax}. Given a class
$\mathcal{C}$ of pairs of probability measures and a probability
metric $\rho$, define the class of $\Delta$-separated alternatives as
\[
  \mathcal{P}_\Delta \coloneqq \bigl\{(P, Q) \in \mathcal{C}^2 :
  \rho(P, Q) \geq \Delta\bigr\}.
\]
Fix a level $\alpha$ and a tolerance $\delta \in (0, 1-\alpha)$. The minimax
separation rate $\Delta^*$ is the smallest $\Delta = \Delta_{n,m}$ (smallest in
the sense of its order of decay as $n, m \to \infty$) for which some
level-$\alpha$ test has Type-II error at most $\delta$ uniformly over
$\mathcal{P}_\Delta$; for any $\Delta$ decaying strictly faster, no such test
exists. A test achieving this over $\mathcal{P}_\Delta$ with
$\Delta \asymp \Delta^*$, i.e.\ up to constants independent of $n$ and $m$, is
called minimax optimal. The choice of class $\mathcal{C}$ and metric $\rho$
specifies the notion of regularity and separation under consideration, and
different choices give rise to different minimax rates and to different
optimal tests.

\paragraph{The MMD test.}
Let $\mathcal{H}$ be a reproducing kernel Hilbert space (RKHS) \citep{Aronszajn} with bounded
measurable kernel $k$ satisfying $\sup_{x\in\mathcal{X}}k(x,x)\le\kappa$.
The \emph{kernel mean embedding} of a probability measure $P$ is the mean
element $\mu_P:=\int_{\mathcal{X}}k(\cdot,x)\,dP(x)\in\mathcal{H}$, which
satisfies $\langle f,\mu_P\rangle_{\mathcal{H}}=\int f\,dP$ for all
$f\in\mathcal{H}$, and the maximum mean discrepancy between $P$ and $Q$ is the
RKHS distance between the corresponding mean elements,
\begin{equation}\label{eq:mmd}
D_{\mathrm{MMD}}(P,Q):=\|\mu_P-\mu_Q\|_{\mathcal{H}}
=\sup_{\|f\|_{\mathcal{H}}\le 1}\Big|\int f\,dP-\int f\,dQ\Big|,
\end{equation}
so that
$D^2_{\mathrm{MMD}}(P,Q)=\mathbb{E}k(X,X')+\mathbb{E}k(Y,Y')-2\,\mathbb{E}k(X,Y)$
for mutually independent $X,X'\stackrel{iid}{\sim}P$ and
$Y,Y'\stackrel{iid}{\sim}Q$. If $k$ is characteristic
\citep{JMLR_metrics, JMLR_universal}, then $D_{\mathrm{MMD}}(P,Q)=0$ if and
only if $P=Q$, which makes $D_{\mathrm{MMD}}$ a natural discrepancy on which
to base a test. Replacing the expectations above by their unbiased two-sample
$U$-statistic estimator yields
\begin{equation}\label{eq:mmd-hat}
\widehat{D}^2_{\mathrm{MMD}}
:=\frac{1}{n(n-1)}\sum_{i\ne j}k(X_i,X_j)
+\frac{1}{m(m-1)}\sum_{i\ne j}k(Y_i,Y_j)
-\frac{2}{nm}\sum_{i,j}k(X_i,Y_j).
\end{equation}
The asymptotic null distribution of \eqref{eq:mmd-hat} does not have a simple
closed form: it is an infinite weighted sum of centered chi-squared random
variables, with weights given by the eigenvalues of the integral operator
associated with the centered version of $k$ w.r.t.~the common distribution under $H_0$. Consequently,
several approximations of the null quantile have been proposed, each yielding
a consistent test \citep[e.g.,][]{gretton2012kernel, FastKernel}.
Alternatively, one may reject $H_0$ when $\widehat{D}^2_{\mathrm{MMD}}$
exceeds the $(1-\alpha)$ quantile of its permutation distribution
\citep{permutations}, which has the advantage of being exactly level $\alpha$
in finite samples. Taking $\mathcal{C}$ to be all pairs of probability measures on $\mathcal{X}$
and $\rho=D_{\mathrm{MMD}}$ gives the class of MMD-separated alternatives
\begin{equation}\label{eq:mmd-class}
\mathcal{P}^{\mathrm{MMD}}_{\Delta}
:=\big\{(P,Q):D_{\mathrm{MMD}}(P,Q)\ge\Delta\big\},
\end{equation}
over which \citet{KimSchrab2026} recently showed the minimax separation rate
to be $\Delta^{*}\asymp (m+n)^{-1/2}$; this rate is attained by
the MMD test. Hence, the MMD test is minimax optimal for
MMD-separated alternatives.

However, recent work, as discussed below, has shown that the MMD test with a fixed kernel, i.e., the
kernel is not dependent on the sample size, e.g., the bandwidth of the
Gaussian kernel on $\mathbb{R}^d$ is fixed and is not dependent on the sample
size, is not minimax optimal under natural classes of $L^2$ or
$\chi^2$-separated alternatives, and has identified two structurally distinct
modifications of the MMD test that restore minimax optimality w.r.t.~these
other alternatives. We refer to these as the \emph{convolution} approach and
the \emph{regularization} approach, and briefly describe each below.

\paragraph{The convolution approach.} On $\mathcal{X} =
\mathbb{R}^d$, \citet{li2019optimality} and \citet{schrab2023mmd}
observed that an MMD test with a fixed-bandwidth kernel cannot
adapt its spatial scale to the regularity of the density difference
$p-q$, and proposed using a translation-invariant kernel
$\tilde{k}_h(x, y) = h^{-d}\varphi((x-y)/h)$ with a shrinking bandwidth
$h \to 0$, where $\varphi : \mathbb{R}^d \to \mathbb{R}$ is a fixed bounded function taken as a product of one-dimensional translation-invariant kernels. The resulting test statistic behaves as a convolution-based estimator of $\|p-q\|_2^2$, i.e., a plug-in estimator of $\int (\tilde{k}_h\star (p-q))(x)\,d(P-Q)(x)$, and an appropriate choice of bandwidth tied to the smoothness of $p-q$ yields a test that is minimax optimal over the Sobolev class of $L^2$-separated alternatives,
\begin{equation}\label{eq:sobolev-alt}
  \mathcal{F}_{s,\Delta} \coloneqq \left\{
  (P, Q) : \max(\|p\|_\infty, \|q\|_\infty) \leq M,\;
  p - q \in \mathcal{S}_d^s(L),\;
  \|p - q\|_2^2 \geq \Delta
  \right\},
\end{equation}
where $\star$ denotes convolution operation and $p,q$ denote the densities of $P$ and $Q$, respectively, with respect to the Lebesgue measure on $\mathbb{R}^d$. Here, $\mathcal{S}_d^s(L)$ is the Sobolev ball with smoothness $s > 0$ and
radius $L > 0$ on $\mathbb{R}^d$, defined as
\[
  \mathcal{S}_d^s(L) \coloneqq
  \left\{f \in L^1(\mathbb{R}^d) \cap L^2(\mathbb{R}^d) :
  \int_{\mathbb{R}^d}\|\xi\|_2^{2s}|\hat{f}(\xi)|^2\,d\xi
  \leq (2\pi)^d L^2\right\},
\]
with $\hat{f}(\xi) = \int f(x)e^{-i\langle x,\xi\rangle}dx$ being the Fourier transform of $f$;
elements of $\mathcal{S}_d^s(L)$ are functions whose Fourier
coefficients decay at a rate determined by $s$.

In particular, under $L^2$-separation, the modified MMD test (i.e., with the bandwidth varying with the sample size) attains the minimax rate
$\Delta \asymp (m+n)^{-4s/(4s+d)}$.
Therefore, the convolution approach
relies essentially
on the translation invariance of the kernel and the Fourier-analytic regularity of $p-q$ in the Sobolev ball on $\mathbb{R}^d$.

\paragraph{The regularization approach.} On a general
measurable space $\mathcal{X}$,
\citet{Krishna, hagrass2023Gof, hagrass2024a_spectral_twosample} observed that the
suboptimality of MMD has a different structural source. In particular in \citet{hagrass2024a_spectral_twosample}, letting
$R = (P+Q)/2$, $u = dP/dR - 1$, and $\mathcal T : L^2(R) \to L^2(R)$ the
centered integral operator of the kernel with respect to $R$, the squared
MMD defined in \eqref{eq:mmd} admits the spectral expansion
$$D^2_{\mathrm{MMD}}(P,Q) = 4\langle\mathcal{T}u,u\rangle_{L^2(R)}=4\sum_{i\geq 1}
\lambda_i\langle u, \tilde{\phi}_i\rangle^2_{L^2(R)},$$
in which the high-frequency Fourier coefficients of $u$, i.e., $\langle u,\tilde{\phi}_i\rangle_{L^2(R)}$ for large $i$ with $\tilde{\phi}_i$ being the eigenfunctions of $\mathcal{T}$, are
downweighted by the decaying eigenvalues $\lambda_i$ of $\mathcal T$---under the boundedness of $k$, $\mathcal{T}$ can be shown to be trace class, which means the eigenvalues $\lambda_i\rightarrow 0$ as $i\rightarrow\infty$. This means, the differences between $P$ and $Q$ can be hidden in the high-frequency components of $u$, and because of the down-weighting of these coefficients, more samples would be needed to detect these differences, which intuitively explains the suboptimality of MMD.

To correct this downweighting,
\citet{hagrass2024a_spectral_twosample} introduced the spectral
regularized discrepancy
\[
  \eta_\lambda \coloneqq 4\langle g_\lambda(\mathcal{T})\mathcal{T}u,u\rangle_{L^2(R)}=
  \bigl\|g_\lambda^{1/2}(\Sigma_R)(\mu_P - \mu_Q)\bigr\|_{\mathcal{H}}^2,
\]
where \[
  \Sigma_R \coloneqq \int_\mathcal{X}
  \bigl(k(\cdot, x) - \mu_R\bigr)
  \otimes_\mathcal{H}
  \bigl(k(\cdot, x) - \mu_R\bigr)\, dR(x):\mathcal{H}\rightarrow\mathcal{H},
\] is the centered covariance operator, and $g_\lambda$ is a spectral regularizer \citep{Engl.et.al},
a function comparable, up to constants, to $x \mapsto (x+\lambda)^{-1}$, so that
$g_\lambda(\Sigma_R)$ acts as a regularized inverse of $\Sigma_R$. The
discrepancy $\eta_\lambda$ is thus a regularized Mahalanobis-type distance
between the mean embeddings, in the spirit of Hotelling's $T^2$-statistic
\citep{lehmann2005testing}; with the Tikhonov regularizer
$g_\lambda(x) = (x+\lambda)^{-1}$, it is closely related to the kernel Fisher
discriminant statistic of \citet{Harchaoui}, which uses the pooled within-sample
covariance in place of $\Sigma_R$. As $\lambda \to 0$,
the regularizer $g_\lambda(\mathcal T)$ flattens the spectral weighting and
$\eta_\lambda$ becomes comparable to $\|u\|^2_{L^2(R)} = \rho^2(P, Q)$.  \citet{hagrass2024a_spectral_twosample}
proved that the resulting test (based on $\eta_\lambda$) is minimax optimal over the
alternative class
\begin{equation}\label{eq:rkhs-alt}
  \mathcal{P}_{\theta,\Delta} \coloneqq \left\{
  (P,Q)\;:\; u \in \Ran(\mathcal T^\theta),\;\;
  \rho^2(P,Q) \geq \Delta
  \right\},
\end{equation}
under separation in the Hellinger-type distance
$\rho^2(P,Q) = \tfrac{1}{2}\int (dP-dQ)^2/(dP+dQ)$, where
$\theta > 0$ controls the smoothness of the class and, for $\theta < 1/2$,
$\Ran(\mathcal T^\theta)$ is an interpolation space between $L^2(R)$ and
the RKHS $\mathcal{H}$ \citep{Steinwart2012MercersTO}. The smoothness in $\mathcal{P}_{\theta,\Delta}$ is
determined by the interaction between the kernel and the reference
measure rather than by Fourier regularity, making the class natural for general domains. The same construction covers the goodness-of-fit problem: taking $R = P_0$ and
$u = dP/dP_0 - 1$, the separation measure $\|u\|^2_{L^2(R)}$ becomes the
$\chi^2$-divergence $\chi^2(P \,\|\, P_0)$, and
\citet{hagrass2023Gof} establish the minimax optimality of the corresponding
spectral regularized test over the analogue of $\mathcal{P}_{\theta,\Delta}$
under $\chi^2$-separation. Note that $\eta_\lambda$ is precisely the squared MMD associated with the regularized kernel $k_\lambda(x,y):=\sum_{i\ge 1}\lambda_i g_\lambda(\lambda_i)\tilde\phi_i(x)\tilde\phi_i(y)$, $x,y\in\mathcal{X}$, in place of the centered kernel $\bar k(x,y):=\langle k(\cdot,x)-\mu_R,k(\cdot,y)-\mu_R\rangle_{\mathcal{H}}=\sum_{i\ge 1}\lambda_i\tilde\phi_i(x)\tilde\phi_i(y)$, whose MMD coincides with that of $k$; here $(\lambda_i,\tilde\phi_i)_{i\ge 1}$ are the eigenpairs of $\mathcal{T}$, and the second series is the Mercer representation of $\bar k$ \citep{Steinwart2012MercersTO}.

\medskip
\noindent The two approaches thus address different sources of MMD's suboptimality under
different structural assumptions: the convolution approach is confined to
$\mathbb{R}^d$ but attains Sobolev minimax optimality through bandwidth
selection alone, with no regularization parameter; the regularization approach
applies on general measurable spaces and attains Hellinger-type optimality
through spectral regularization at a fixed bandwidth. Neither analysis speaks
to the other's alternative class, and the two tuning mechanisms have so far
been studied in isolation.

This raises a structural question: are the bandwidth and the regularization
parameter genuinely distinct devices, or is one an instance of the other? We
show that, on $\mathbb{R}^d$, both serve to satisfy a single condition. For a
translation-invariant kernel $k(x,y)=K(x-y)$ with spectral density
$\widehat K$, the spectral regularized discrepancy captures the separation of
the Sobolev class, in the sense that $\eta_\lambda\gtrsim\|p-q\|_2^2$ for every
$(P,Q)\in\mathcal{F}_{s,\Delta}$, as soon as
\[
  \lambda \;\lesssim\; \min_{\|\omega\|_2\,\leq\, c\,\Delta^{-1/(2s)}}\widehat K(\omega),
\]
with $c$ depending only on $s$ and $L$; that is, as soon as $\lambda$ lies below
the spectral density of the kernel on all frequencies up to order
$\Delta^{-1/(2s)}$ (Theorem~\ref{thm:signal}). For kernels of finite smoothness
the condition is also necessary, up to constants: if $\lambda$ exceeds the
spectral density at frequencies of order $\Delta^{-1/(2s)}$ by a large factor,
then $\eta_\lambda$ is much smaller than $\|p-q\|_2^2$ for some alternative in
$\mathcal{F}_{s,\Delta}$ (Proposition~\ref{prop:necessity}). There are two ways
of meeting this condition. The first is to shrink the bandwidth: for $K=\varphi(\cdot/h)$
the spectral density is of order $h^d$ on frequencies up to $1/h$, so the
condition holds with $h\asymp\Delta^{1/(2s)}$ and $\lambda\asymp h^d$, the
regularization parameter decaying together with the bandwidth. The second is to
let $\lambda$ decay with the kernel held fixed, which is possible when the
spectral density decays only polynomially: for kernels of finite smoothness,
such as the Mat\'ern, Laplace and Wendland kernels, the spectral regularized test
with a fixed bandwidth attains the Sobolev minimax rate with only $\lambda$
decaying, under an eigenvalue condition met, for instance, by distributions
supported in a bounded set. The convolution approach is an instance of the first
route. With the constant regularizer $g_\lambda\equiv\lambda^{-1}$, the spectral
regularized test is the MMD test, for every $\lambda$, and our analysis applies
to it along the coupling $\lambda\asymp h^d$ (Proposition~\ref{prop:mmd}), so
the optimality of the MMD test with a shrinking bandwidth over
$\mathcal{F}_{s,\Delta}$ \citep{li2019optimality, schrab2023mmd} is recovered as
a special case of ours. This separates two reasons for shrinking the bandwidth.
For the MMD test, whose statistic does not depend on $\lambda$, the bandwidth is
the only device available, whatever the kernel; for a spectral regularized test,
our analysis requires a shrinking bandwidth only when the spectral density
decays faster than any polynomial, as for the Gaussian kernel, for which the
required $\lambda$ is too small to be used with a covariance operator estimated
from the sample. On a
general domain, where there is no bandwidth to tune, a decaying $\lambda$ is
also what delivers optimality over $\mathcal{P}_{\theta,\Delta}$. The spectral regularized test, on a general measurable space, is thus
minimax optimal over both alternative classes.

\medskip
\noindent\textbf{Contributions.} The main contributions of this paper are as follows.

\begin{enumerate}[label=(\roman*)]

\item \emph{A spectral condition for Sobolev alternatives.}
Theorem~\ref{thm:signal} shows that, for any translation-invariant kernel with
positive spectral density, $\eta_\lambda \gtrsim \|p-q\|_2^2$ as soon as
$\lambda \lesssim m_K(1/t)=\min_{\|\omega\|_2\le 1/t}\widehat K(\omega)$, i.e., the smallest value of the kernel's spectral density
on frequencies up to $1/t$, and $\|p-q\|_2^2 \gtrsim t^{2s}$. The frequency
cutoff $1/t$ is set by the separation: at separation $\Delta$ the condition reads
$\lambda \lesssim m_K(c\,\Delta^{-1/(2s)})$, that is, $\lambda$ must not exceed
the kernel's spectral density on the frequencies up to order $\Delta^{-1/(2s)}$. Conversely,
Proposition~\ref{prop:necessity} shows that, for kernels of finite smoothness,
the condition is also necessary, up to constants, for
$\eta_\lambda \gtrsim \|p-q\|_2^2$ to hold uniformly over
$\mathcal{F}_{s,\Delta}$, so that it characterizes when the discrepancy captures
the separation. A shrinking bandwidth and a decaying regularization parameter
are two ways of meeting this condition.

\item \emph{The convolution approach as a spectral regularized test.} For a
bandwidth family, $m_K(1/h) \asymp h^d$, and Theorem~\ref{thm:separation} shows
that the spectral regularized test with $h \asymp N^{-2/(4s+d)}$ and
$\lambda \asymp h^d$ attains the Sobolev minimax rate $N^{-4s/(4s+d)}$ for every
regularizer obeying a two-sided bound on $g_\lambda(x)(x+\lambda)$. The constant
regularizer $g_\lambda\equiv\lambda^{-1}$ obeys this bound along the coupling
$\lambda\asymp h^d$ on the spectrum of $\Sigma_R$, and with it the test is the MMD test, for every $\lambda$,
with no covariance samples needed (Proposition~\ref{prop:mmd}). The
convolution test of \citet{li2019optimality, schrab2023mmd} is therefore a member
of the spectral regularized family, and its optimality over
$\mathcal{F}_{s,\Delta}$ follows, with a short proof, for any
translation-invariant kernel with positive spectral density.

\item \emph{Regularization without shrinking the bandwidth.} In
Theorem~\ref{thm:sobolev-fixed} we show that a fixed kernel of finite smoothness,
whose spectral density decays no faster than $(1+\|\omega\|_2^2)^{-\tau}$ for some
$\tau > d/2$, such as the Mat\'ern, Laplace and Wendland kernels, attains the
Sobolev minimax rate with only the regularization parameter decaying,
$\lambda \asymp N^{-4\tau/(4s+d)}$, under an eigenvalue decay condition that holds
for distributions supported in a fixed bounded set. The regularizer then does the
work of the bandwidth.

\item \emph{Adaptation over both classes, $\mathcal{F}_{s,\Delta}$ and $\mathcal{P}_{\theta,\Delta}$.} In Theorem~\ref{thm:adaptation} we
construct an adaptive test that aggregates over a grid of $(\lambda, h)$ values
and is minimax optimal, up to a $\log\log N$ factor, simultaneously over
$\mathcal{F}_{s,\Delta}$ and $\mathcal{P}_{\theta,\Delta}$, without knowledge of
$s$ or $\theta$. 

\item The analysis covers both the two-sample and the goodness-of-fit
problems. All results are stated for a general covariance reference measure
$R$, subject only to the dominance condition
$\Sigma_P,\Sigma_Q\preccurlyeq c\,\Sigma_R$, and apply verbatim to
goodness-of-fit testing on setting $Q=P_0$ for a known reference distribution
$P_0$. Section~\ref{sec:reference-measure} collects the natural choices of $R$
in each setting and records the separation rate that each entails. 

\end{enumerate}

\medskip
\noindent Together, these results show that the bandwidth and the
regularization parameter are two ways of meeting a single condition on the
spectral density of the kernel,
that the convolution approach is the member of the spectral regularized family
with the constant regularizer, and that for kernels of finite smoothness
regularization alone attains the Sobolev rate. A single spectral test, defined on a
general measurable space, is minimax optimal over both the operator-theoretic
class and the Sobolev class on $\mathbb{R}^d$; Table~\ref{tab:routes} summarizes
the routes to the Sobolev rate. The proofs of all results are provided in the
Appendix.
\section{Setup}\label{sec:setup}

We first introduce the notation and definitions used throughout, and then
the spectral regularized test statistic and its permutation test.

\subsection{Notation}\label{sec:notation}

For a measurable space $\mathcal{X}$, let $L^r(\mathcal{X}, \mu)$
denote the Banach space of $r$-power integrable functions under
measure $\mu$. Let $\mathcal{H}$ be a reproducing kernel Hilbert
space (RKHS) on $\mathcal{X}$ with reproducing kernel
$k : \mathcal{X} \times \mathcal{X} \to \mathbb{R}$, satisfying
$\sup_{x \in \mathcal{X}} k(x, x) \leq \kappa$. For a bounded linear
operator $S : \mathcal{H} \to \mathcal{H}$, let $\|S\|_{\mathrm{op}}$
and $\|S\|_{\mathcal{L}^2(\mathcal{H})}$ denote the operator and
Hilbert--Schmidt norms, respectively, and, for $S$ self-adjoint, let
$\mathrm{spec}(S)$ denote its spectrum; for a compact positive $S$, this is
the set of eigenvalues of $S$, together with $0$ when $\mathcal{H}$ is
infinite-dimensional. For $f, g \in \mathcal{H}$, the
tensor product $f \otimes_\mathcal{H} g$ is identified with the
operator on $\mathcal{H}$ acting as
$(f \otimes_\mathcal{H} g)z = \langle z, g\rangle_\mathcal{H} f,\,\forall\,z\in\mathcal{H}$.

The kernel mean embedding of a probability measure $P$ is
$\mu_P = \int k(\cdot, x)\, dP(x) \in \mathcal{H}$. For a probability
measure $R$ on $\mathcal{X}$, the covariance operator is
\[
  \Sigma_R \coloneqq \int_\mathcal{X}
  \bigl(k(\cdot, x) - \mu_R\bigr)
  \otimes_\mathcal{H}
  \bigl(k(\cdot, x) - \mu_R\bigr)\, dR(x),
\]
and we write $\Sigma_{R,\lambda} \coloneqq \Sigma_R + \lambda I$ for
$\lambda > 0$.

A spectral regularizer is a function
$g_\lambda : [0, \infty) \to (0, \infty)$, applied to a compact, positive,
self-adjoint operator $\mathcal{D}$ on a separable Hilbert space $H$ as
\begin{equation}\label{Eq:glambda}
  g_{\lambda}(\mathcal{D}) := \sum_{i\geq 1} g_{\lambda}(\sigma_i)
  (\phi_i \otimes_H \phi_i) + g_{\lambda}(0)\left(\Id -
  \sum_{i\geq 1} \phi_i \otimes_H \phi_i\right),
\end{equation}
where $\mathcal D$ has the spectral representation
$\mathcal D = \sum_{i}\sigma_i \phi_i \otimes_H \phi_i$ with
$(\sigma_i,\phi_i)_i$ being the eigenvalues and eigenfunctions of
$\mathcal D$. The Tikhonov regularizer $g_\lambda(x) = (x + \lambda)^{-1}$
and the Showalter regularizer $g_\lambda(x) = (1 - e^{-x/\lambda})/x$, with
$g_\lambda(0) = \lambda^{-1}$ by continuity, are canonical examples. The only condition we impose on $g_\lambda$ is a two-sided
bound on $g_\lambda(x)(x+\lambda)$, introduced in Section~\ref{sec:power}: the
power analysis uses both halves, and the lower bounds on the discrepancy in
Section~\ref{sec:minimax} use only the lower half. In particular, we need no
approximation (qualification) condition of the kind that is standard in the
inverse problem literature (see, for instance, \citealp{Engl.et.al}) and has
been employed in kernel ridge regression \citep{BAUER200752} and in earlier
analyses of kernel tests \citep{hagrass2023Gof,hagrass2024a_spectral_twosample}
(Remark~\ref{rem:qualification}). The population
discrepancy associated with $g_\lambda$ is defined as
\begin{equation}\label{eq:eta-pop}
  \eta \coloneqq
  \bigl\|g_\lambda^{1/2}(\Sigma_R)(\mu_P - \mu_Q)\bigr\|_\mathcal{H}^2.
\end{equation}

For sequences $(a_r)_r, (b_r)_r$, we write $a_r \lesssim b_r$ (resp.\ $a_r \gtrsim b_r$)
if there exists a constant $c > 0$, not depending on $r$ (but possibly on
the fixed parameters of the problem), with $a_r \leq cb_r$ (resp.\
$a_r \geq cb_r$), and $a_r \asymp b_r$ if both hold. Throughout,
$m \leq n \leq Dm$ for some constant $D > 1$, and
$N \coloneqq m + n$.

\subsection{Test statistic}\label{sec:teststat}

Let $\{X_i\}_{i=1}^n \overset{\mathrm{iid}}{\sim} P$,
$\{Y_j\}_{j=1}^m \overset{\mathrm{iid}}{\sim} Q$, and
$\mathbb{Z}_l \coloneqq \{Z_k\}_{k=1}^l \overset{\mathrm{iid}}{\sim} R$ be mutually
independent---see below for the discussion about the probability measure $R$.
In the goodness-of-fit setting, the samples $\{Y_j\}_{j=1}^m$ are drawn from
the known reference $P_0$. The covariance samples $\{Z_k\}_{k=1}^l$ are
obtained differently in the two settings. In the two-sample setting, $R$
depends on the unknown $P$ and $Q$, so the $Z_k$ must be constructed by sample
splitting: for the mixture $R=(P+Q)/2$, one sets aside a portion of each
sample and draws each $Z_k$ from the held-out $P$-samples or the held-out
$Q$-samples according to an independent fair coin flip, so that
$Z_k \overset{\mathrm{iid}}{\sim} R$ and $\{Z_k\}_{k=1}^l$ is independent of
the samples retained for the test statistic. In the goodness-of-fit setting
with $R=P_0$, no splitting is needed: the $Z_k$ are additional samples drawn
directly from the known $P_0$, at no cost to the test budget. Throughout, the reference measure $R$ is treated as a free parameter; concrete
choices, including the goodness-of-fit case, are discussed in
Section~\ref{sec:reference-measure}.

Let
$\hat{\Sigma}_R$ denote the U-statistic estimator of $\Sigma_R$ based on
$\{Z_k\}_{k=1}^l$
\citep[Section~4.1]{hagrass2024a_spectral_twosample}, and define
$A(x) \coloneqq g_\lambda^{1/2}(\hat{\Sigma}_R)\, k(\cdot, x)$.

The
test statistic is defined as
\begin{equation}\label{eq:test-stat}
  \hat{\eta}
  \coloneqq
  \frac{1}{n(n-1)}\sum_{i\neq j}\ipH{A(X_i)}{A(X_j)}
  + \frac{1}{m(m-1)}\sum_{i\neq j}\ipH{A(Y_i)}{A(Y_j)}
  - \frac{2}{nm}\sum_{i,j}\ipH{A(X_i)}{A(Y_j)},
\end{equation}
which is a two-sample U-statistic estimator of $\eta$ in
\eqref{eq:eta-pop}. Let $\Pi_{n+m}$ denote the set of all permutations of the
pooled sample $(\mathbb{X}_n, \mathbb{Y}_m)$, and for $\pi \in \Pi_{n+m}$ let
$\hat\eta^\pi$ be the statistic computed on the permuted samples. The
permutation distribution function and its $(1-\alpha)$ quantile are
\[
  F(x) \coloneqq \frac{1}{|\Pi_{n+m}|}\sum_{\pi \in \Pi_{n+m}}
  \mathbf{1}(\hat\eta^\pi \leq x),
  \qquad
  \hat q_{1-\alpha} \coloneqq \inf\{q \in \mathbb{R} : F(q) \geq 1-\alpha\},
\]
and the test rejects $H_0$ when $\hat\eta > \hat q_{1-\alpha}$. By the
exchangeability of the pooled sample under $H_0$, this test has level $\alpha$
exactly, for any $n$ and $m$.

\begin{remark}[Monte Carlo approximation of the quantile]
\label{rem:practical-quantile}
Computing $\hat q_{1-\alpha}$ exactly requires all $(n+m)!$ permutations and is
infeasible in practice. One therefore approximates it by
\[
  \hat q^{B}_{1-\alpha} \coloneqq \inf\Big\{q \in \mathbb{R} :
  \frac{1}{B}\sum_{i=1}^{B}\mathbf{1}(\hat\eta^{\pi^i} \leq q) \geq 1-\alpha
  \Big\},
\]
where $\pi^1, \dots, \pi^B$ are drawn uniformly from $\Pi_{n+m}$. By the
Dvoretzky--Kiefer--Wolfowitz inequality \citep{DKW, Massart}, $\hat
q^{B}_{1-\alpha}$ approaches $\hat q_{1-\alpha}$ up to a shift of the level by an
arbitrarily small fixed fraction of $\alpha$, with probability at least $1-\delta$, once
$B \gtrsim \alpha^{-2}\log(2\delta^{-1})$. All results stated below for
$\hat q_{1-\alpha}$ therefore extend to $\hat q^{B}_{1-\alpha}$ for $B$ of this
order.
\end{remark}

\medskip\noindent
The kernel $k$ is taken to be bounded throughout, but is otherwise unrestricted. For the Sobolev class in Section~\ref{sec:minimax}, $k$ is a translation-invariant kernel on $\mathbb{R}^d$, either taken from a bandwidth family or held fixed; the relation between the bounded bandwidth parameterization and the convolution kernel $h^{-d}\varphi((x-y)/h)$ used by \citet{schrab2023mmd} is discussed in Section~\ref{sec:mmd}.

\section{Power analysis}\label{sec:power}

In this section, we analyze the spectral regularized permutation
test in a class-agnostic manner. We bound the variance of the test
statistic, control the permutation threshold on the scale of this
bound, and combine the two to obtain a power criterion stated
as a lower bound on the population discrepancy $\eta$ relative to
the regularization parameter $\lambda$, the covariance sample size
$l$, and the test sample size $N$. The criterion does not depend
on the specific alternative class; applying it to a class reduces
to lower-bounding $\eta$ by the corresponding separation measure
under the regularity condition defining the class, which is the
content of Section~\ref{sec:minimax}.

\subsection{Variance and threshold bounds}

Define the effective dimension
\[
  \mathcal{N}_2(\lambda)
  \coloneqq
  \bigl\|\Sigma_{R,\lambda}^{-1/2}
  \Sigma_R
  \Sigma_{R,\lambda}^{-1/2}
  \bigr\|_{\mathcal{L}^2(\mathcal{H})},
\]
which controls the null-regime fluctuation of $\hat\eta$.

\begin{proposition}[Variance bound]\label{prop:var-bound}
Let $\mathcal B = g_\lambda^{1/2}(\hat{\Sigma}_R)\,\Sigma_{R,\lambda}^{1/2}$, $\zeta \coloneqq \E(\hat{\eta} \mid \mathbb{Z}_l)$, and suppose
$\Sigma_P \preccurlyeq c\,\Sigma_R$ and
$\Sigma_Q \preccurlyeq c\,\Sigma_R$ for some constant $c > 0$. Then
\[
  \mathrm{Var}(\hat{\eta} \mid \mathbb{Z}_l)
  \;\lesssim\;
  \frac{\|\mathcal B\|_{\mathrm{op}}^4\,\mathcal{N}_2^2(\lambda)}
       {N^2}
  \;+\;
  \frac{\|\mathcal B\|_{\mathrm{op}}^2\,\zeta}
       {N}.
\]
\end{proposition}

The variance decomposes into a null-regime term
$\mathcal{N}_2^2(\lambda)/N^2$ that captures the fluctuation of
$\hat\eta$ under $P = Q$, and a discrepancy-dependent term
$\zeta/N$.

Write $C_0$ for the constant implicit in Proposition~\ref{prop:var-bound}, and
\begin{equation}\label{eq:V-def}
  V \;\coloneqq\; C_0\left(
  \frac{\|\mathcal B\|_{\mathrm{op}}^4\,\mathcal{N}_2^2(\lambda)}{N^2}
  + \frac{\|\mathcal B\|_{\mathrm{op}}^2\,\zeta}{N}
  + \frac{\zeta^2}{N^2}\right)
  \;\geq\; \mathrm{Var}(\hat{\eta} \mid \mathbb{Z}_l)
\end{equation}
for the resulting upper bound on the variance; the last term is not needed
for the variance, but accounts for the spread of the permutation distribution
under $H_1$. The following result shows that the permutation quantile $\hat q_{1-\alpha}$ is bounded, with high
probability under $H_1$, on the same scale.

\begin{lemma}[Threshold bound]\label{lem:threshold}
For $0 < \alpha \leq e^{-1}$, $\delta > 0$,
there exists a constant $c_1 > 0$ such that
\[
  P_{H_1}\!\left\{\hat{q}_{1-\alpha} \leq c_1
  \sqrt{\frac{V}{\delta}}\right\}
  \geq 1-\frac{\delta}{2}.
\]
\end{lemma}

\subsection{Power criterion}\label{sec:power-criterion}

Combining the variance and threshold bounds with a Chebyshev-type argument
yields the power criterion below. Define 
\[
  \mathcal{B}_0 \;\coloneqq\; g_\lambda^{1/2}(\Sigma_R)\,\Sigma_{R,\lambda}^{1/2}.
\]  For constants $C_u \geq C_\ell > 0$, we write
$\mathcal{E} = \mathcal{E}(C_u, C_\ell)$ for the event on which
\begin{equation}\label{eq:event}
  \|\mathcal{B}\|_{\mathrm{op}} \vee \|\mathcal{B}_0\|_{\mathrm{op}}
  \;\leq\; C_u,
  \qquad
  \|\mathcal{B}^{-1}\|_{\mathrm{op}} \;\leq\; C_\ell^{-1}.
\end{equation}

The upper bound on
$\|\mathcal{B}\|_{\mathrm{op}}$ controls the variance through
Proposition~\ref{prop:var-bound}; the bounds on
$\|\mathcal{B}^{-1}\|_{\mathrm{op}}$ and $\|\mathcal{B}_0\|_{\mathrm{op}}$
together ensure that the conditional mean
$\zeta = \E(\hat\eta \mid \mathbb{Z}_l)$ does not fall below the population
discrepancy $\eta$.
\begin{theorem}[Power criterion]\label{thm:power}
Let $0 < \alpha \leq e^{-1}$ and $\delta > 0$. Suppose
$\Sigma_P \preccurlyeq c\,\Sigma_R$ and $\Sigma_Q \preccurlyeq c\,\Sigma_R$,
and that $P(\mathcal{E}) \geq 1 - \delta/2$. There exist constants $c_2 > 0$,
depending only on $c$, $D$, $\alpha$, $C_u$ and $C_\ell$, and $N_0 > 0$,
depending only on $c$, $D$, $\alpha$ and $\delta$, such that if
$N \geq N_0$ and
\begin{equation}\label{cond:eta}
  \eta \;\geq\; \frac{c_2\,\bigl(\mathcal{N}_2(\lambda) \vee 1\bigr)}{\delta\,N},
\end{equation}
then $P_{H_1}\{\hat\eta > \hat q_{1-\alpha}\} \geq 1 - 2\delta$.
\end{theorem}

The event $\mathcal{E}$ is guaranteed by a two-sided bound on the regularizer
together with enough covariance samples. Suppose that
\begin{equation}\label{eq:g-two-sided}
   C_\ell^2 \;\leq\; g_\lambda(x)\,(x+\lambda) \;\leq\; C_u^2
  \qquad \text{for all } x \in [0,2\kappa],
\end{equation}
as holds for the Tikhonov regularizer with $C_\ell = C_u = 1$ and for the
Showalter regularizer with $C_\ell = 1$ and $C_u^2 \leq 4/3$ (the interval
$[0,2\kappa]$ contains the spectra of $\Sigma_R$ and $\hat\Sigma_R$, since
$\|\hat\Sigma_R\|_{\mathrm{op}} \leq l\kappa/(l-1)$), and that
\begin{equation}\label{eq:lambda-range}
  \frac{140\kappa}{l}\log\frac{128\kappa l}{\delta}
  \;\leq\; \lambda \;\leq\; \|\Sigma_R\|_{\mathrm{op}}.
\end{equation}
Then $P\bigl(\mathcal{E}(\sqrt{2}\,C_u, \sqrt{2/3}\,C_\ell)\bigr) \geq 1-\delta/2$
by Lemma~\ref{lem:B-general} in Appendix~\ref{app:event}, so
Theorem~\ref{thm:power} applies. The lower end of~\eqref{eq:lambda-range} asks
for $l \gtrsim \lambda^{-1}\log(\lambda^{-1}\delta^{-1})$ covariance samples.
When $R$ is known, these are drawn at no cost to the test budget; when $R$
depends on $(P,Q)$, they must be split from it, and this requirement is the
source of the sample-splitting regimes in the separation rates of
Section~\ref{sec:minimax}.

\section{Minimax optimality}\label{sec:minimax}

In this section, we instantiate Theorem~\ref{thm:power} for the two
alternative classes $\mathcal{F}_{s,\Delta}$ and $\mathcal{P}_{\theta,\Delta}$.
Applying it to a class requires a lower bound on the population discrepancy
$\eta$ in terms of the separation measure of that class. For the Sobolev class
this is Theorem~\ref{thm:signal} (Section~\ref{sec:eta-lb}), which holds for any
translation-invariant kernel and reduces to a single condition on its spectral
density: $\lambda$ must not exceed the spectral density on the frequencies up to
order $\Delta^{-1/(2s)}$, a cutoff set by the separation. A shrinking bandwidth
and a decaying regularization parameter are two ways of meeting this condition,
and they lead to the two Sobolev separation theorems:
Theorem~\ref{thm:separation} (Section~\ref{sec:shrinking}), in which the
bandwidth shrinks with the sample size, and Theorem~\ref{thm:sobolev-fixed}
(Section~\ref{sec:fixed-kernel}), in which the kernel is held fixed and only the
regularization parameter decays. In between, Section~\ref{sec:mmd} identifies the
convolution approach as the member of the family obtained with a constant
regularizer, and explains in spectral terms why it must shrink the bandwidth. For completeness, we also state the
results of the general-domain class \citep{hagrass2024a_spectral_twosample,hagrass2023Gof} under the current framework in Section~\ref{sec:general-domain}.

In each separation theorem we distinguish two scenarios for the covariance
reference measure $R$: \emph{independent covariance samples}, in which $R$ is
fully known and fresh samples can be drawn at no cost to the test budget; and
\emph{sample splitting}, in which $R$ depends on $P, Q$ and estimating
$\Sigma_R$ consumes part of the test samples. Concrete choices of $R$
corresponding to each scenario, including the goodness-of-fit setting where
$Q = P_0$ is known, are discussed in Section~\ref{sec:reference-measure}.

For the Sobolev class we work with translation-invariant kernels on $\mathbb{R}^d$.

\begin{assumption}[Translation-invariant kernel]\label{ass:TI}
$k(x,y)=K(x-y)$, where $K\in C(\mathbb{R}^d)\cap L^1(\mathbb{R}^d)$ is real-valued and positive definite, and $\widehat K(\omega)>0$ for all $\omega\in\mathbb{R}^d$.
\end{assumption}

Under Assumption~\ref{ass:TI}, $k$ is bounded with $\kappa=K(0)$, and its RKHS is \citep[Theorem~10.12]{Wendland2005}
\begin{equation}\label{eq:native}
\mathcal{H}=\Big\{f\in L^2(\mathbb{R}^d)\cap C(\mathbb{R}^d):\ \widehat f/\sqrt{\widehat K}\in L^2(\mathbb{R}^d)\Big\},
\qquad
\|f\|_{\mathcal{H}}^2=\frac{1}{(2\pi)^d}\int_{\mathbb{R}^d}\frac{|\widehat f(\omega)|^2}{\widehat K(\omega)}\,d\omega .
\end{equation}
The quantity that governs the lower bound on $\eta$ in Section~\ref{sec:eta-lb} is the smallest value of the spectral density on a ball of frequencies,
\begin{equation}\label{eq:mK}
m_K(\rho):=\min_{\|\omega\|_2\le\rho}\widehat K(\omega),\qquad \rho>0,
\end{equation}
which is positive and non-increasing in $\rho$. We consider two instances of Assumption~\ref{ass:TI}.
\begin{enumerate}
\item[(a)] \emph{Bandwidth family.} $K=\varphi(\cdot/h)$ for a fixed $\varphi$ satisfying Assumption~\ref{ass:TI}. Then $\widehat K(\omega)=h^d\widehat\varphi(h\omega)$, the bound $\kappa=\varphi(0)$ does not depend on $h$, and
\begin{equation}\label{eq:mK-bandwidth}
m_K(1/h)=c_\varphi\, h^d,\qquad c_\varphi:=\min_{\|u\|_2\le 1}\widehat\varphi(u)>0 .
\end{equation}
Boundedness uniformly in $h$ is needed for the operator-theoretic analysis underlying Theorem~\ref{thm:power}; the relation between this bounded parameterization and the unbounded convolution kernel $h^{-d}\varphi((x-y)/h)$ used by \citet{schrab2023mmd} is discussed in Section~\ref{sec:mmd}.
\item[(b)] \emph{Fixed kernel of finite smoothness.} $K$ does not depend on the sample size, and for some $\tau>d/2$ and $c_K>0$,
\begin{equation}\label{eq:Khat-lower}
\widehat K(\omega)\;\ge\;c_K\,(1+\|\omega\|_2^2)^{-\tau},\qquad \omega\in\mathbb{R}^d,
\end{equation}
so that $m_K(\rho)\ge c_K(1+\rho^2)^{-\tau}$. Examples are the Mat\'ern kernel with smoothness $\nu=\tau-d/2$ and any length-scale, the Laplace kernel $e^{-\|x-y\|_2}$ ($\tau=(d+1)/2$), and Wendland's compactly supported kernels $\varphi_{d,k}$, $k\ge1$ ($\tau=(d+1)/2+k$) \citep[Theorem~10.35]{Wendland2005}. 
\end{enumerate}

The two instances differ in how the effective dimension $\mathcal{N}_2(\lambda)$
of Section~\ref{sec:power} is controlled. Without any condition on the
eigenvalues of $\Sigma_R$, bounding the squared Hilbert--Schmidt norm by the
operator norm, which is at most one, times the trace gives
\[
  \mathcal{N}_2^2(\lambda)
  \;\leq\; \tr\bigl(\Sigma_{R,\lambda}^{-1/2}\Sigma_R\Sigma_{R,\lambda}^{-1/2}\bigr)
  \;\leq\; \frac{\tr(\Sigma_R)}{\lambda}
  \;\leq\; \frac{\kappa}{\lambda}.
\]
For a bandwidth family this generic bound suffices, since the relevant
$\lambda$ is of order $h^d$. For a fixed kernel with $\lambda\to0$ it is too
crude, and, as in \citet{hagrass2024a_spectral_twosample}, we use instead a
polynomial decay of the eigenvalues $\lambda_1(\Sigma_R)\ge\lambda_2(\Sigma_R)\ge\cdots$ of $\Sigma_R$: for some $A>0$ and the same $\tau$ as in \eqref{eq:Khat-lower},
\begin{equation}\label{eq:eigdecay}
\lambda_i(\Sigma_R)\le A\,i^{-2\tau/d},\qquad i\ge 1,
\end{equation}
uniformly over the alternatives under consideration. Under \eqref{eq:eigdecay}, $\mathcal{N}_2(\lambda)\lesssim\lambda^{-d/(4\tau)}$ \citep[Lemma B.9]{sriperumbudur2022approx}. Conditions \eqref{eq:Khat-lower} and \eqref{eq:eigdecay} are complementary: the first says that $\widehat K$ decays no faster than $(1+\|\omega\|_2^2)^{-\tau}$, the second that the eigenvalues of $\Sigma_R$ decay at least as fast as $i^{-2\tau/d}$. The second is not implied by the first, since it also depends on $R$. For the kernels listed in (b), $\widehat K$ is also bounded above by a multiple of $(1+\|\omega\|_2^2)^{-\tau}$, and \eqref{eq:eigdecay} then holds, with $A$ depending only on $K$, $M$ and $\Omega$, whenever $R$ has a density bounded by $M$ and supported in a bounded set $\Omega$ \citep[e.g.][]{EdmundsTriebel1996}.

\subsection{Discrepancy lower bound for Sobolev alternatives}\label{sec:eta-lb}

The following result gives a lower bound on $\eta$ for Sobolev alternatives.

\begin{theorem}[Lower bound on $\eta$, Sobolev alternatives]\label{thm:signal}
Let Assumption~\ref{ass:TI} hold, let $\eta$ be as defined
in~\eqref{eq:eta-pop}, let $\psi = p - q$, and let $r$ denote the density of $R$
w.r.t.~the Lebesgue measure on $\mathbb{R}^d$. Suppose $g_\lambda$ satisfies
\begin{equation}\label{eq:g-lower}
  g_\lambda(x)\,(x+\lambda) \;\geq\; C_\ell^2
  \qquad\text{for all } x \in \mathrm{spec}(\Sigma_R),
\end{equation}
and that $\psi \in \mathcal{S}_d^s(L)$ and $\|r\|_\infty \leq M$. If, for some
$t > 0$,
\begin{equation}\label{eq:resolution}
  \lambda \;\leq\; \frac{M}{4}\, m_K(1/t)
  \qquad\text{and}\qquad
  \|\psi\|_2^2 \;\geq\; 4L^2\, t^{2s},
\end{equation}
then
\[
  \eta \;\geq\; \frac{C_\ell^2}{2M}\,\|\psi\|_2^2 .
\]
\end{theorem}

Theorem~\ref{thm:signal} transfers the Fourier smoothness of $\psi=p-q$ to a
lower bound on $\eta$ through a frequency cutoff $1/t$: the frequencies of $\psi$
up to $1/t$ are retained provided $\lambda$ does not exceed the spectral density
of the kernel on that band. The cutoff is set by the separation. Since
$m_K$ is non-increasing, the first condition is weakest for the largest
admissible $t$, so that for $(P,Q) \in \mathcal{F}_{s,\Delta}$ the two conditions
reduce to the single requirement
\begin{equation}\label{eq:spectral-condition}
  \lambda \;\leq\; \frac{M}{4}\, m_K\bigl((4L^2/\Delta)^{1/(2s)}\bigr):
\end{equation}
$\lambda$ must not exceed the kernel's spectral density on the frequencies up to
$(4L^2/\Delta)^{1/(2s)}$. A bandwidth family meets this requirement by shrinking the kernel, which
keeps $m_K(1/t)$ of order $h^d$ at $t=h$; a fixed kernel meets it by shrinking
$\lambda$, at a polynomial rate when the spectral density decays polynomially
(Remark~\ref{rem:resolution}). 

\begin{remark}[Choice of witness function]\label{rem:witness}
The proof of Theorem~\ref{thm:signal} lower-bounds $\eta$ by testing against
$f_t$ with $\widehat{f_t}=\widehat\psi\,\mathbf{1}\{\|\omega\|_2\le1/t\}$, the
convolution of $\psi$ with the ideal low-pass filter at scale $t$. The cutoff
$1/t$ is set by the separation alone, so this choice separates the frequency
cutoff from the kernel, and with it the role of the regularization parameter
from that of the bandwidth: the bandwidth enters only through the spectral
density in~\eqref{eq:resolution}. The analyses of the convolution approach
\citep{li2019optimality, schrab2023mmd} rely instead on the convolution of $\psi$
with the kernel itself, which for a bandwidth family is
$h^{-d}\int k(\cdot,y)\psi(y)\,dy=h^{-d}(\mu_P-\mu_Q)$, the MMD witness function.
Its frequency response $\widehat\varphi(h\omega)$ reaches frequencies only up to
order $1/h$, and its RKHS norm forces $\lambda\lesssim h^d$, whatever the kernel,
so an analysis based on this witness ties the frequency cutoff to the bandwidth.
The band-limited $f_t$ is what reveals condition~\eqref{eq:spectral-condition}
on the spectral density, which a decaying $\lambda$ can meet with the kernel held
fixed (Section~\ref{sec:fixed-kernel}).
\end{remark}

Condition~\eqref{eq:spectral-condition} is not only sufficient. For kernels
whose spectral density decays polynomially, as for those listed in~(b), it is
also necessary, up to constants, for $\eta$ to capture the separation uniformly
over $\mathcal{F}_{s,\Delta}$.

\begin{proposition}[Necessity of the spectral condition]\label{prop:necessity}
Let Assumption~\ref{ass:TI} hold with, for some $\tau > d/2$ and
$0 < c_K \leq C_K$,
\[
  c_K(1+\|\omega\|_2^2)^{-\tau} \;\leq\; \widehat K(\omega) \;\leq\;
  C_K(1+\|\omega\|_2^2)^{-\tau}, \qquad \omega\in\mathbb{R}^d,
\]
and let $g_\lambda$ satisfy the upper bound in~\eqref{eq:g-two-sided}. There
are constants $\Delta_0, c_3 > 0$, depending only on $d$, $s$, $L$, $M$, $\tau$,
$c_K$ and $C_K$, such that for every $\Delta \in (0, \Delta_0]$ some
$(P,Q) \in \mathcal{F}_{s,\Delta}$ with $\|p-q\|_2^2 = \Delta$ satisfies,
whatever the reference measure $R$,
\[
  \eta \;\leq\; c_3\,\frac{C_u^2}{\lambda}\,
  m_K\bigl((4L^2/\Delta)^{1/(2s)}\bigr)\,\|p-q\|_2^2 .
\]
Consequently, if $\eta \geq c\,\|p-q\|_2^2$ for some $c > 0$ and all
$(P,Q) \in \mathcal{F}_{s,\Delta}$, then
\[
  \lambda \;\leq\; \frac{c_3C_u^2}{c}\, m_K\bigl((4L^2/\Delta)^{1/(2s)}\bigr).
\]
\end{proposition}

\begin{remark}[The two instances]\label{rem:resolution}
For the bandwidth family (a), taking $t=h$ and using~\eqref{eq:mK-bandwidth}, condition~\eqref{eq:resolution} reads $\lambda\le\frac{Mc_\varphi}{4}h^d$ and $\|\psi\|_2^2\ge4L^2h^{2s}$: the regime $\lambda\lesssim h^d$, $\|\psi\|_2^2\gtrsim h^{2s}$ of the convolution approach, with explicit constants. For a fixed kernel (b), \eqref{eq:Khat-lower} shows that the first condition in~\eqref{eq:resolution} holds as soon as $\lambda\le\frac{Mc_K}{4}(1+t^{-2})^{-\tau}$, so the cutoff is reached by $\lambda$ alone, with no bandwidth involved: since $(1+t^{-2})^{-\tau}\ge2^{-\tau}t^{2\tau}$ for $t\le1$, condition~\eqref{eq:spectral-condition} holds once $\lambda\lesssim\Delta^{\tau/s}$.
\end{remark}

\subsection{Shrinking bandwidth}\label{sec:shrinking}

For a bandwidth family~(a), Theorem~\ref{thm:signal} with $t=h$ and the trace
bound on $\mathcal{N}_2(\lambda)$ combine with Theorem~\ref{thm:power} as follows.

\begin{theorem}[Separation rate, Sobolev alternatives]
\label{thm:separation}
Let $0 < \alpha \leq e^{-1}$, $\delta > 0$, $M, L > 0$, $s > 0$, and let
Assumption~\ref{ass:TI} hold with $K=\varphi(\cdot/h)$, the bandwidth
family~(a). Suppose
$\Sigma_P \preccurlyeq c\,\Sigma_R$, $\Sigma_Q \preccurlyeq c\,\Sigma_R$,
$\|r\|_\infty \leq M$, that $g_\lambda$ satisfies~\eqref{eq:g-two-sided}, and
that $\inf_{(P,Q)\,\in\,\mathcal{F}_{s,\Delta}}\|\Sigma_R\|_{\mathrm{op}}
\gtrsim h^d$.
\begin{enumerate}[label=\textnormal{(\roman*)}]
  \item \textnormal{ (No cost to the test budget.)} Choose
  \[
    h \asymp N^{-2/(4s+d)},
    \qquad
    \lambda \asymp N^{-2d/(4s+d)}.
  \]
  Then for any $\Delta \gtrsim N^{-4s/(4s+d)}$,
  \[
    \inf_{(P,Q)\,\in\,\mathcal{F}_{s,\Delta}}
    P_{H_1}\bigl\{\hat{\eta} > \hat{q}_{1-\alpha}\bigr\}
    \geq 1-2\delta.
  \]
  \item \textnormal{(Sample splitting, $l \asymp N$.)}
  \begin{enumerate}[label=\textnormal{(\alph*)}]
    \item If $s > d/4$, the choice
    $h \asymp N^{-2/(4s+d)}$, $\lambda \asymp N^{-2d/(4s+d)}$
    suffices for $\Delta \gtrsim N^{-4s/(4s+d)}$.
    \item If $s \leq d/4$, the choice
    $h \asymp (\log N / N)^{1/d}$,
    $\lambda \asymp \log N / N$
    suffices for
    $\Delta \gtrsim (\log N / N)^{2s/d}$.
  \end{enumerate}
  In both sub-cases,
  $\inf_{(P,Q)\,\in\,\mathcal{F}_{s,\Delta}}
  P_{H_1}\{\hat{\eta} > \hat{q}_{1-\alpha}\} \geq 1-2\delta.$
\end{enumerate}
\end{theorem}
In part~(i), $R$ is known and the covariance samples required
by~\eqref{eq:lambda-range}, of order $h^{-d}\log(h^{-d}\delta^{-1})$, are drawn at
no cost to the test budget; in part~(ii) they are split from it, which is what
produces the threshold $s>d/4$. The MMD test needs no covariance samples at all
and attains the rate of part~(i) for every $s>0$
(Proposition~\ref{prop:mmd} in Section~\ref{sec:mmd}).
\begin{remark}\label{rem:op-lower}
The lower bound on $\|\Sigma_R\|_{\mathrm{op}}$ is used only to ensure that
$\lambda \asymp h^d$ satisfies the upper end of~\eqref{eq:lambda-range}, and is
mild. Taking $f = k(\cdot,x_0)$ gives
$\|\Sigma_R\|_{\mathrm{op}} \geq \mathrm{Var}_R(f)/\varphi(0)$, which is of
order $h^d\,r(x_0)$ at any Lebesgue point $x_0$ of $r$; the condition
therefore holds for any fixed $R$ with a density, and for $R$ depending on
$(P,Q)$, such as the mixture $(P+Q)/2$, it asks only that the densities in
the class be bounded below near some point, uniformly.
\end{remark}

\subsection{The convolution approach as a spectral regularized test}\label{sec:mmd}

We now make precise the sense in which the convolution approach of
\citet{li2019optimality} and \citet{schrab2023mmd} belongs to the spectral
regularized family. Take the constant regularizer $g_\lambda \equiv
\lambda^{-1}$. By~\eqref{Eq:glambda}, $g_\lambda^{1/2}(\hat\Sigma_R) =
\lambda^{-1/2}I$, so $A(x) = \lambda^{-1/2}k(\cdot,x)$ and the
statistic~\eqref{eq:test-stat} is
$\hat\eta = \lambda^{-1}\widehat{D}^2_{\mathrm{MMD}}$, with
$\widehat{D}^2_{\mathrm{MMD}}$ as in~\eqref{eq:mmd-hat}. The permutation quantile
scales in the same way, so the test rejects exactly when
$\widehat{D}^2_{\mathrm{MMD}}$ exceeds its own permutation quantile: it is the MMD
permutation test, whatever $\lambda$ is, and it does not use the covariance
samples. Since rescaling the kernel by a positive constant rescales
$\widehat{D}^2_{\mathrm{MMD}}$ and its permutation quantile alike, the same holds
for the normalized kernel $h^{-d}\varphi((x-y)/h)$ of \citet{schrab2023mmd}.

\begin{proposition}[MMD with a shrinking bandwidth]\label{prop:mmd}
Let $0 < \alpha \leq e^{-1}$, $\delta > 0$, $M, L > 0$, $s > 0$, and let
$\varphi$ satisfy Assumption~\ref{ass:TI}. The MMD permutation test with kernel
$\varphi((x-y)/h)$, or with any positive multiple of it, coincides with the
spectral regularized test with $g_\lambda \equiv \lambda^{-1}$, for every
$\lambda > 0$ and every reference measure $R$. With $h \asymp N^{-2/(4s+d)}$,
it has power at least $1-2\delta$ uniformly over $\mathcal{F}_{s,\Delta}$ for
any $\Delta \gtrsim N^{-4s/(4s+d)}$.
\end{proposition}

The rate is that of \citet{li2019optimality} and \citet{schrab2023mmd}, here for
any $\varphi$ satisfying Assumption~\ref{ass:TI}, not necessarily of product
form. The proof, in Appendix~\ref{app:mmd}, runs the argument of
Theorem~\ref{thm:separation} with $R=(P+Q)/2$, which enters only the analysis:
for the constant regularizer the event $\mathcal{E}$ holds surely once
$\lambda \gtrsim \|\Sigma_R\|_{\mathrm{op}}$, so neither covariance samples nor
a lower bound on $\|\Sigma_R\|_{\mathrm{op}}$ are needed, and the sharp rate holds
for every $s > 0$. What the proposition adds to the known rate is the mechanism,
which the spectral picture makes explicit.

\begin{remark}[Why the convolution approach shrinks the bandwidth]\label{rem:why-shrink}
The constant regularizer has to meet two conditions at once. On
$\mathrm{spec}(\Sigma_R)$ it must satisfy the two-sided
bound~\eqref{eq:g-two-sided}, and since $g_\lambda(x)(x+\lambda) = 1+x/\lambda$,
this asks for $\lambda \gtrsim \|\Sigma_R\|_{\mathrm{op}}$
(Lemma~\ref{lem:B-flat}). It must also satisfy the first condition
in~\eqref{eq:resolution}, $\lambda \leq \frac{M}{4}m_K(1/t)$, at the cutoff set
by the separation. Together, these require
\[
  \|\Sigma_R\|_{\mathrm{op}} \;\lesssim\; \lambda \;\lesssim\; m_K(1/t).
\]
For a fixed kernel, $\|\Sigma_R\|_{\mathrm{op}}$ does not shrink while
$m_K(1/t) \to 0$ as $t \to 0$, so the two conditions are incompatible for small
$\Delta$. For a bandwidth family, both sides are of order $h^d$ at $t = h$, since
$\|\Sigma_R\|_{\mathrm{op}} \lesssim h^d$ (see the proof of
Theorem~\ref{thm:separation}) and $m_K(1/h) = c_\varphi h^d$
by~\eqref{eq:mK-bandwidth}, so the coupling $\lambda \asymp h^d$ meets both.
Shrinking the bandwidth is thus what lowers $\|\Sigma_R\|_{\mathrm{op}}$ to the
level of the spectral density at the cutoff. Since the constant regularizer gives
the MMD test for every $\lambda$, the MMD test has the bandwidth as its only
device, for every kernel, the Laplace kernel as much as the Gaussian, in line with
the known suboptimality of the MMD test with a fixed kernel
\citep{Krishna, hagrass2024a_spectral_twosample}.
\end{remark}

\paragraph{Kernel normalization.}
The convolution test of \citet{schrab2023mmd} works with the normalized kernel
$h^{-d}\varphi((x-y)/h)$, which is unbounded as $h \to 0$, whereas our analysis
uses the bounded kernel $k=\varphi(\cdot/h)$. For the MMD test the normalization
is immaterial, as noted above; for the spectral regularized test it amounts to a
rescaling of the regularization parameter. Let
\begin{equation}\label{eq:normalized-kernel}
  \tilde k(x, y) \coloneqq h^{-\gamma}\, k(x, y)
\end{equation}
for some $\gamma > 0$, the choice $\gamma = d$ giving the kernel of
\citet{schrab2023mmd}, and let $\tilde{\mathcal H}$ denote the RKHS of $\tilde k$,
with mean embeddings $\tilde\mu_P, \tilde\mu_Q$ and covariance operators
$\tilde\Sigma_R$, $\hat{\tilde\Sigma}_R$.

\begin{lemma}[Scaling relations]\label{lem:scaling}
The RKHS of $\tilde k$ coincides with $\mathcal{H}$ as a set, with
\[
  \langle f, g\rangle_{\tilde{\mathcal H}}
  = h^\gamma\,\langle f, g\rangle_{\mathcal H},
  \qquad
  f \otimes_{\tilde{\mathcal H}} g
  = h^\gamma\, f \otimes_{\mathcal H} g,
  \qquad f, g \in \mathcal{H}.
\]
Consequently, $\tilde\Sigma_R = h^{-\gamma}\Sigma_R$ and
$\hat{\tilde\Sigma}_R = h^{-\gamma}\hat\Sigma_R$, and for any
$\tilde\lambda > 0$,
\[
  \tilde\Sigma_{R,\tilde\lambda}
  = h^{-\gamma}\Sigma_{R, h^\gamma\tilde\lambda}.
\]
\end{lemma}

\begin{proposition}[Normalization as a rescaling of $\lambda$]\label{prop:normalization}
Suppose $g_\lambda(x) = \lambda^{-1}G(x/\lambda)$ for a function
$G:[0,\infty)\to(0,\infty)$, as for the Tikhonov ($G(u)=(1+u)^{-1}$), Showalter
($G(u)=(1-e^{-u})/u$) and constant ($G\equiv1$) regularizers. Then, for every
$\tilde\lambda > 0$, the spectral regularized test based on $\tilde k$ with
regularization $\tilde\lambda$ coincides with the test based on $k$ with
regularization $\lambda \coloneqq h^\gamma\tilde\lambda$: the two statistics are
equal, for the observed sample and for each of its permutations.
\end{proposition}

For a general $g_\lambda$ the two statistics need not coincide, but the analysis
still transfers with $\lambda = h^\gamma\tilde\lambda$
(Proposition~\ref{prop:equivalence} in Appendix~\ref{app:normalization}), so no
generality is lost by working with the bounded kernel. Under the normalization
$\gamma = d$ of \citet{schrab2023mmd}, the regularization
$\lambda \asymp h^d \asymp N^{-2d/(4s+d)}$ of Theorem~\ref{thm:separation}
corresponds to $\tilde\lambda \asymp 1$, a constant in $N$: the dependence on $N$
that appears in the bounded-kernel parameterization is absorbed into the kernel
by the $h^{-d}$ normalization.

\subsection{Fixed kernel: regularization without shrinking the bandwidth}\label{sec:fixed-kernel}

For a fixed kernel of finite smoothness~(b),
condition~\eqref{eq:spectral-condition} is met by $\lambda$ alone
(Remark~\ref{rem:resolution}), and the effective dimension is controlled
through~\eqref{eq:eigdecay}.

\begin{theorem}[Separation rate, Sobolev alternatives, fixed kernel]
\label{thm:sobolev-fixed}
Let $0 < \alpha \leq e^{-1}$, $\delta > 0$, $M, L > 0$, $s > 0$, and let
Assumption~\ref{ass:TI} hold with a fixed kernel $K$
satisfying~\eqref{eq:Khat-lower} for some $\tau > d/2$, as in~(b). Suppose
$\Sigma_P \preccurlyeq c\,\Sigma_R$, $\Sigma_Q \preccurlyeq c\,\Sigma_R$,
$\|r\|_\infty \leq M$, that $g_\lambda$ satisfies~\eqref{eq:g-two-sided},
that~\eqref{eq:eigdecay} holds with the same $A$ for all
$(P,Q)\in\mathcal{F}_{s,\Delta}$, and that
$\inf_{(P,Q)\,\in\,\mathcal{F}_{s,\Delta}}\|\Sigma_R\|_{\mathrm{op}} > 0$.
\begin{enumerate}[label=\textnormal{(\roman*)}]
  \item \textnormal{(No cost to the test budget.)} Choose
  $\lambda \asymp N^{-4\tau/(4s+d)}$. Then for any
  $\Delta \gtrsim N^{-4s/(4s+d)}$,
  \[
    \inf_{(P,Q)\,\in\,\mathcal{F}_{s,\Delta}}
    P_{H_1}\bigl\{\hat{\eta} > \hat{q}_{1-\alpha}\bigr\}
    \geq 1-2\delta.
  \]
  \item \textnormal{(Sample splitting, $l \asymp N$.)}
  \begin{enumerate}[label=\textnormal{(\alph*)}]
    \item If $s > \tau - d/4$, the choice
    $\lambda \asymp N^{-4\tau/(4s+d)}$ suffices for
    $\Delta \gtrsim N^{-4s/(4s+d)}$.
    \item If $s \leq \tau - d/4$, the choice
    $\lambda \asymp \log N / N$ suffices for
    $\Delta \gtrsim (\log N / N)^{s/\tau}$.
  \end{enumerate}
  In both sub-cases,
  $\inf_{(P,Q)\,\in\,\mathcal{F}_{s,\Delta}}
  P_{H_1}\{\hat{\eta} > \hat{q}_{1-\alpha}\} \geq 1-2\delta.$
\end{enumerate}
\end{theorem}

Here the kernel is held fixed and only $\lambda$ decays, at the rate
$\lambda\asymp N^{-4\tau/(4s+d)}$; the guarantee involves no bandwidth.
The constant regularizer is of no help here: $\lambda\to0$ lies below
$\|\Sigma_R\|_{\mathrm{op}}$, which does not shrink, and with a fixed kernel the
analysis of the MMD test cannot meet~\eqref{eq:spectral-condition} for small
$\Delta$ (Remark~\ref{rem:why-shrink}). $\Sigma_R$ must therefore be estimated, and the
sample-splitting threshold $s>\tau-d/4$ takes the place of the threshold $s>d/4$
of Theorem~\ref{thm:separation}.

For the kernels listed in~(b), the conditions
on $\Sigma_R$ hold when $P$ and $Q$ are supported in a fixed bounded set and $R$
is $(P+Q)/2$ or the uniform distribution on that set. This subclass of
$\mathcal{F}_{s,\Delta}$ contains the standard lower-bound constructions
\citep{ingster1987minimax, li2019optimality}, so the rate in parts~(i)
and~(ii)(a) is minimax optimal.

\begin{remark}[Regularization in place of the bandwidth]\label{rem:fixed-resolution}
The proof applies Theorem~\ref{thm:signal} with the cutoff $1/t$ at
$t\asymp\lambda^{1/(2\tau)}\asymp N^{-2/(4s+d)}$, which is the optimal bandwidth
of Theorem~\ref{thm:separation}, with $\lambda\asymp t^{2\tau}$; for the kernels
listed in~(b), this is the order of $\widehat K$ at frequency $1/t$, so the
regularization parameter sits at the spectral density of the kernel at the
cutoff frequency. For a bandwidth family, $m_K(1/h)\asymp h^d$
by~\eqref{eq:mK-bandwidth}, which is the coupling $\lambda\asymp h^d$ of
Theorem~\ref{thm:separation}; for a fixed kernel, the same cutoff is reached by
$\lambda$ alone. As $\tau\downarrow d/2$, the threshold $s>\tau-d/4$ and the
rate $(\log N/N)^{s/\tau}$ of part~(ii)(b) tend to the threshold $s>d/4$ and the
rate $(\log N/N)^{2s/d}$ of Theorem~\ref{thm:separation}(ii).
\end{remark}

\begin{remark}[The Gaussian kernel]\label{rem:gaussian}
A decaying $\lambda$ replaces bandwidth shrinkage only when the spectral density
at the frequencies required by~\eqref{eq:spectral-condition} is not too small.
Two requirements are involved, with $t \asymp \Delta^{1/(2s)}$ the scale set by
the separation. First, $\lambda \leq \frac{M}{4}m_K(1/t)$ must be usable with
the available covariance samples, which, when $\Sigma_R$ is estimated from
$l\lesssim N$ samples, means $\lambda\gtrsim\log N/N$ by~\eqref{eq:lambda-range}.
Second, at that $\lambda$ the effective dimension must be of order $t^{-d/2}$,
so that condition~\eqref{cond:eta} is met at the minimax separation. Since
$\mathcal{N}_2^2(\lambda)=\sum_i\lambda_i^2/(\lambda_i+\lambda)^2$, with
$\lambda_i=\lambda_i(\Sigma_R)$, is of the order of the number of eigenvalues
above $\lambda$ for the decay rates considered here, this asks that about
$t^{-d}$ eigenvalues exceed $\lambda$, the order of the volume of the ball of
frequencies up to $1/t$. For the kernels of finite smoothness both requirements
hold: $\lambda\asymp t^{2\tau}$ is polynomial in $N$, and by~\eqref{eq:eigdecay}
at most $(A/\lambda)^{d/(2\tau)}\asymp t^{-d}$ eigenvalues exceed $\lambda$, so
that $\mathcal{N}_2(\lambda)\lesssim t^{-d/2}$. Here the eigenvalues follow the
spectral density: by~\eqref{eq:eigdecay} and~\eqref{eq:Khat-lower}, the $i$-th
is at most of the order of $\widehat K$ at frequency $i^{1/d}$. For a kernel whose spectral density
decays faster than any polynomial, the first requirement fails. For the Gaussian
kernel $K(z)=e^{-\|z\|_2^2/(2h_0^2)}$, the most common example,
$m_K(1/t)=(2\pi)^{d/2}h_0^d\,e^{-h_0^2/(2t^2)}$ is exponentially small in
$(h_0/t)^2$, so $\lambda\gtrsim\log N/N$ limits the cutoff to
$t\gtrsim h_0(2\log N)^{-1/2}$ and yields only a logarithmic rate. Within this
analysis, the Gaussian bandwidth must therefore shrink both with the constant
regularizer (Remark~\ref{rem:why-shrink}) and when $\Sigma_R$ is estimated from
the sample. With unlimited covariance samples, so that $\lambda$ may be
exponentially small, the second requirement fails instead, because the
eigenvalues no longer follow the spectral density. For $R$ supported in a
bounded set, the eigenvalues of the integral operator of the Gaussian kernel,
and hence those of $\Sigma_R$ (centering does not increase them), are at most
of order $e^{-c\,i^{1/d}}$ \citep[Theorem~5]{belkin2018approximation}, a bound
that decays much more slowly than the spectral density at frequency $i^{1/d}$,
which decays like $e^{-c'i^{2/d}}$. At
$\lambda\asymp e^{-h_0^2/(2t^2)}$, up to $(\log(1/\lambda))^d\asymp t^{-2d}$
eigenvalues may therefore exceed $\lambda$, instead of $t^{-d}$, and the
resulting bound is $\mathcal{N}_2(\lambda)\lesssim t^{-d}$ rather than
$t^{-d/2}$. The additional directions add to the variance but not to the lower
bound on $\eta$, which uses only the frequencies of $\psi$ up to $1/t$:
condition~\eqref{cond:eta} now asks for $\Delta\gtrsim t^{-d}/N$ instead of
$\Delta\gtrsim t^{-d/2}/N$, and balancing it against $\Delta\gtrsim t^{2s}$ from
Theorem~\ref{thm:signal} gives $t\asymp N^{-1/(2s+d)}$ and only
$\Delta\asymp N^{-2s/(2s+d)}$, instead of $t\asymp N^{-2/(4s+d)}$ and
$\Delta\asymp N^{-4s/(4s+d)}$. We do not pursue sharper approximation arguments
for this case.
\end{remark}

\begin{remark}[The role of the regularizer for a fixed kernel]\label{rem:regimes}
In Theorem~\ref{thm:sobolev-fixed}, as in Theorem~\ref{thm:separation-rkhs}
below, the kernel is fixed, $\|\Sigma_R\|_{\mathrm{op}}$ does not shrink, and
$\lambda\to0$ is small compared with it, so the regularizers are genuinely
distinct and the constant one is unavailable. For the Sobolev class
the choice of $g_\lambda$ still does not affect the rate and no qualification
condition is needed, because Theorem~\ref{thm:signal} uses $g_\lambda$ only
through the lower bound~\eqref{eq:g-lower}: the error beyond the frequency
cutoff $1/t$ is controlled by the Sobolev condition on $\psi$, not by how closely
$xg_\lambda(x)$ approximates one. The same holds for the general-domain class:
Proposition~\ref{prop:signal-rkhs} below also uses $g_\lambda$ only
through~\eqref{eq:g-lower}, the error below the eigenvalue cutoff $\lambda$
being controlled by the range condition on $w$ (Remark~\ref{rem:qualification}).
\end{remark}

Table~\ref{tab:routes} summarizes the ways of
meeting~\eqref{eq:spectral-condition} at the minimax separation.

\begin{table}[t]
\centering
\small
\renewcommand{\arraystretch}{1.15}
\begin{tabular}{@{}lcccc@{}}
\toprule
 & \multicolumn{2}{c}{Shrinking bandwidth} & \multicolumn{2}{c}{Fixed bandwidth ($h_0$)} \\
\cmidrule(lr){2-3}\cmidrule(l){4-5}
 & constant $g_\lambda$ (MMD) & spectral $g_\lambda$ & $\widehat K$ polynomial decay & $\widehat K$ Gaussian decay \\
\midrule
$\lambda$ & $h^{d}$ & $h^{d}$ & $t^{2\tau}$ & $e^{-h_0^2/(2t^2)}$ \\
Rate, $\Sigma_R$ no cost & sharp, every $s$ & sharp, every $s$ & sharp, every $s$ & $N^{-2s/(2s+d)}$ \\
Rate, $\Sigma_R$ split & sharp, every $s$ & sharp if $s > d/4$ & sharp if $s > \tau - d/4$ & logarithmic \\
Adaptation over & $h$ (MMDAgg) & $(\lambda, h)$ & $\lambda$ & -- \\
\bottomrule
\end{tabular}
\caption{Ways of reaching the minimax
separation $\Delta \asymp N^{-4s/(4s+d)}$, that is, with frequency cutoff $1/t$,
$t \asymp N^{-2/(4s+d)}$.}
\label{tab:routes}
\end{table}

\subsection{General-domain alternatives}\label{sec:general-domain}

For the general-domain class, the lower bound on $\eta$ follows from the same
cutoff argument as Theorem~\ref{thm:signal}, applied to the eigenvalues of the
integral operator $\mathcal{T}$ instead of to Fourier frequencies. We state it
for a general reference measure $R$; here $\mathcal{T}:L^2(R)\to L^2(R)$ is the
centered integral operator of $k$ with respect to $R$, whose nonzero eigenvalues
are those of $\Sigma_R$ (Appendix~\ref{app:signal-rkhs}).

\begin{proposition}[Lower bound on $\eta$, general-domain alternatives]
\label{prop:signal-rkhs}
Let $\eta$ be as defined in~\eqref{eq:eta-pop}. Suppose $P, Q \ll R$ and let
$w \coloneqq d(P-Q)/dR \in L^2(R)$. Suppose $g_\lambda$
satisfies~\eqref{eq:g-lower}, and that $w \in \Ran(\mathcal{T}^\theta)$ for
some $\theta > 0$ with
\[
  \|w\|_{L^2(R)}^2 \;\geq\; 2\lambda^{2\theta}\,
  \bigl\|\mathcal{T}^{-\theta}w\bigr\|_{L^2(R)}^2 .
\]
Then
\[
  \eta \;\geq\; \frac{C_\ell^2}{4}\,\|w\|_{L^2(R)}^2 .
\]
In particular, for $R = (P+Q)/2$, $\eta \gtrsim \rho^2(P,Q)$, and for
$R = P_0$ in the goodness-of-fit setting, $\eta \gtrsim \chi^2(P \,\|\, P_0)$.
\end{proposition}

As in Theorem~\ref{thm:signal}, the smoothness of the alternative is
transferred to $\eta$ through a cutoff, here at the level $\lambda$ in the
spectrum of $\mathcal{T}$: the components of $w$ along eigenvalues above
$\lambda$ are retained, since~\eqref{eq:g-lower} gives
$xg_\lambda(x)\geq C_\ell^2/2$ for $x\geq\lambda$, and those along eigenvalues
below $\lambda$ carry at most
$\lambda^{2\theta}\|\mathcal{T}^{-\theta}w\|^2_{L^2(R)}$ by the range
condition. There is no requirement on the bandwidth, and only the lower half
of~\eqref{eq:g-two-sided} is used.

\begin{remark}[No qualification condition]\label{rem:qualification}
In \citet[Lemma~A.7]{hagrass2024a_spectral_twosample}, and in the corresponding
goodness-of-fit result of \citet{hagrass2023Gof}, the lower bound on $\eta$ is
obtained through the regularization bias, under a qualification condition of the
form $|B - xg_\lambda(x)|\,x^{2\nu} \leq C\lambda^{2\nu}$ for all $x$ in the
spectrum of $\mathcal{T}$ with $xg_\lambda(x) < B$ and all $\nu$ up to the
qualification $\xi$ of $g_\lambda$; the rates are stated in terms of
$\tilde\theta = \min(\theta,\xi)$, and the standard regularizers, Tikhonov and
Showalter included, satisfy the condition with $\xi = \infty$.
Proposition~\ref{prop:signal-rkhs} shows that the condition is not needed: the
lower half of~\eqref{eq:g-two-sided} suffices, for every $\theta > 0$. It in
fact implies the condition with $\xi = \infty$, since for $B \leq C_\ell^2/2$ it
confines $\{x : xg_\lambda(x) < B\}$ to $[0,\lambda)$, where
$|B - xg_\lambda(x)|\,x^{2\nu} \leq B\lambda^{2\nu}$ for every $\nu > 0$. The
rates of \citet{hagrass2024a_spectral_twosample} therefore hold with $\theta$ in
place of $\tilde\theta$ for every regularizer satisfying~\eqref{eq:g-two-sided},
and a single condition on $g_\lambda$ serves both classes.
\end{remark}

Beyond the two classes, that of MMD-separated
alternatives requires no smoothness assumption at all: there $\eta$ is bounded
below by a multiple of $D^2_{\mathrm{MMD}}(P,Q)$ at any fixed $\lambda$, and the resulting
separation rate is derived in Remark~\ref{rem:mmd-sep-optimal}.

\begin{theorem}[Separation rate, general-domain alternatives]
\label{thm:separation-rkhs}
Let $0 < \alpha \leq e^{-1}$, $\delta > 0$ and $\theta > 0$. Suppose
$\Sigma_P \preccurlyeq c\,\Sigma_R$, $\Sigma_Q \preccurlyeq c\,\Sigma_R$,
that $g_\lambda$ satisfies~\eqref{eq:g-two-sided}, that
$\inf_{(P,Q)\in\mathcal{P}_{\theta,\Delta}}\|\Sigma_R\|_{\mathrm{op}} > 0$, and
that
$\sup_{(P,Q)\in\mathcal{P}_{\theta,\Delta}}\|\mathcal{T}^{-\theta}w\|_{L^2(R)} < \infty$,
where $w$ is as defined in Proposition~\ref{prop:signal-rkhs}.
\begin{enumerate}[label=\textnormal{(\roman*)}]
  \item \textnormal{(No cost to the test budget.)} Choose
  $\lambda \asymp N^{-2/(4\theta+1)}$. Then for any
  $\Delta \gtrsim N^{-4\theta/(4\theta+1)}$,
  \[
    \inf_{(P,Q)\,\in\,\mathcal{P}_{\theta,\Delta}}
    P_{H_1}\bigl\{\hat{\eta} > \hat{q}_{1-\alpha}\bigr\}
    \geq 1-2\delta.
  \]
  \item \textnormal{(Sample splitting, $l \asymp N$.)}
  \begin{enumerate}[label=\textnormal{(\alph*)}]
    \item If $\theta > 1/4$, the choice
    $\lambda \asymp N^{-2/(4\theta+1)}$ suffices for
    $\Delta \gtrsim N^{-4\theta/(4\theta+1)}$.
    \item If $\theta \leq 1/4$, the choice
    $\lambda \asymp \log N / N$ suffices for
    $\Delta \gtrsim (\log N / N)^{2\theta}$.
  \end{enumerate}
  In both sub-cases,
  $\inf_{(P,Q)\,\in\,\mathcal{P}_{\theta,\Delta}}
  P_{H_1}\{\hat{\eta} > \hat{q}_{1-\alpha}\} \geq 1-2\delta$.
\end{enumerate}
\end{theorem}

\begin{remark}[Relation to \citealp{hagrass2024a_spectral_twosample}]
\label{rem:decay}
Theorem~\ref{thm:separation-rkhs} extends
\citet[Theorem~4.7 and Corollary~4.8]{hagrass2024a_spectral_twosample} and is
included so that the separation rates can be read side by side; its proof is
in Appendix~\ref{app:separation-rkhs}. Three differences are worth noting.
First, the rates there are stated in terms of
$\tilde\theta = \min(\theta,\xi)$, where $\xi$ is the qualification of
$g_\lambda$, under a qualification condition; the rates above hold with
$\theta$ for every regularizer satisfying~\eqref{eq:g-two-sided}, with no such
condition (Remark~\ref{rem:qualification}). Second, the reference measure there is fixed
at $R = (P+Q)/2$, whereas the statement above holds for the general $R$ of
Section~\ref{sec:teststat}, with $\|w\|_{L^2(R)}^2$ in place of
$\rho^2(P,Q)$ as the separation measure. Third, no decay condition is imposed
on the eigenvalues of $\mathcal{T}$: the statement uses only the generic bound
$\mathcal{N}_2(\lambda) \lesssim \lambda^{-1/2}$, which holds for any bounded
kernel and is also the bound used in Theorem~\ref{thm:separation}, and the rates
above therefore correspond to the endpoint $\beta = 1$ of Corollary~4.8 there.
Sharper rates follow by imposing a decay condition and using the corresponding
bound on $\mathcal{N}_2$ in condition~\eqref{cond:eta} of
Theorem~\ref{thm:power}: if $\lambda_i \asymp i^{-\beta}$ for $\beta > 1$, then
$\mathcal{N}_2(\lambda) \asymp \lambda^{-1/(2\beta)}$ and the choice
$\lambda \asymp \Delta^{1/(2\theta)}$ yields
$\Delta \asymp N^{-4\theta\beta/(4\theta\beta+1)}$ for
$\theta > \frac12 - \frac1{4\beta}$; if $\lambda_i \asymp e^{-c_0 i}$
for some $c_0 > 0$, then $\mathcal{N}_2(\lambda) \asymp \sqrt{\log(1/\lambda)}$
and $\Delta \asymp \sqrt{\log N}/N$ for $\theta > 1/2$, the rate of
Corollary~4.9 there. In each case the rate matches the minimax separation
boundary over $\mathcal{P}_{\theta,\Delta}$, so the test remains minimax
optimal under either decay condition.
\end{remark}

Theorems~\ref{thm:separation}, \ref{thm:sobolev-fixed}
and~\ref{thm:separation-rkhs} show that the spectral regularized test attains
the minimax separation rate over both alternative classes. For the Sobolev
class, the rate $\Delta\asymp N^{-4s/(4s+d)}$, established by
\citet{li2019optimality, schrab2023mmd} for the MMD test with a shrinking bandwidth,
is attained in two ways: with a bandwidth family, by letting $h$ and
$\lambda\asymp h^d$ decay jointly (Theorem~\ref{thm:separation}), the MMD test
being the member with a constant regularizer (Proposition~\ref{prop:mmd}); and
with a fixed kernel of finite smoothness, by letting only $\lambda$ decay
(Theorem~\ref{thm:sobolev-fixed}). For the general-domain class, the bandwidth
plays no role and only $\lambda$ decays, recovering the rate of
\citet{hagrass2024a_spectral_twosample} with no qualification condition on the regularizer.

\begin{remark}[MMD-separated alternatives]\label{rem:mmd-sep-optimal}
The spectral regularized test is also minimax optimal over the class
$\mathcal{P}^{\mathrm{MMD}}_\Delta$ of~\eqref{eq:mmd-class}, with no
smoothness assumption on the alternative and no decay of the regularization
parameter. Indeed, the lower bound~\eqref{eq:g-lower} gives
$g_\lambda(x) \geq C_\ell^2/(x+\lambda)$ on the spectrum of $\Sigma_R$, so that
\[
  \eta \;=\; \big\|g_\lambda^{1/2}(\Sigma_R)(\mu_P-\mu_Q)\big\|_{\mathcal{H}}^2
  \;\ge\; C_\ell^2\,
  \big\|\Sigma_{R,\lambda}^{-1/2}(\mu_P-\mu_Q)\big\|_{\mathcal{H}}^2
  \;\ge\; \frac{C_\ell^2}{\|\Sigma_R\|_{\mathrm{op}}+\lambda}\,
  D^2_{\mathrm{MMD}}(P,Q),
\]
and $\|\Sigma_R\|_{\mathrm{op}} \le \kappa$ under the boundedness of $k$.
Take $\lambda = \kappa$ and $R = (P+Q)/2$. Then $\mathcal{N}_2(\lambda) \le 1$,
and the event $\mathcal{E}$ holds surely, whatever $l$: since $[0,2\kappa]$
contains the spectra of $\Sigma_R$ and $\hat\Sigma_R$,
\eqref{eq:g-two-sided} gives
$\|\mathcal{B}\|_{\mathrm{op}}^2 \vee \|\mathcal{B}_0\|_{\mathrm{op}}^2
\le C_u^2(\kappa+\lambda)/\lambda = 2C_u^2$ and
$\|\mathcal{B}^{-1}\|_{\mathrm{op}}^2 \le (2\kappa+\lambda)/(\lambda C_\ell^2)
= 3C_\ell^{-2}$. Condition~\eqref{cond:eta} therefore reduces to
$D^2_{\mathrm{MMD}}(P,Q) \gtrsim 1/(\delta N)$, and the test attains the
separation rate $\Delta \asymp N^{-1/2}$, which is minimax optimal over
$\mathcal{P}^{\mathrm{MMD}}_\Delta$. The spectral regularized test thus retains
the MMD separation rate for every regularizer satisfying~\eqref{eq:g-two-sided},
not merely for the constant one, under which it is the MMD test itself
(Section~\ref{sec:mmd}).
\end{remark}

\section{Choice of reference measure}\label{sec:reference-measure}

The reference measure $R$ is a free parameter of the test, constrained only by
the dominance condition $\Sigma_P, \Sigma_Q \preccurlyeq c\,\Sigma_R$ of
Theorem~\ref{thm:power} and, for the Sobolev class $\mathcal{F}_{s,\Delta}$,
by $\|r\|_\infty \leq M$, together with the eigenvalue
decay~\eqref{eq:eigdecay} when the kernel is fixed
(Theorem~\ref{thm:sobolev-fixed}). The choice matters along three axes: the
separation measure it induces through $\|w\|_{L^2(R)}$, the regularity it
requires, and whether the covariance samples $\{Z_k\}_{k=1}^l$ must be split
from the test budget. The first of these concerns only the general-domain class
$\mathcal{P}_{\theta,\Delta}$, whose separation is measured by
$\|w\|^2_{L^2(R)}$ and therefore varies with $R$; the Sobolev class is
separated in $\|p-q\|_2^2$ irrespective of $R$, which enters
Theorem~\ref{thm:signal} only through $\|r\|_\infty \leq M$.

The third axis is not a property of $R$ alone. For the MMD test of
Proposition~\ref{prop:mmd}, $\Sigma_R$ is never estimated and no covariance
samples are drawn at all; $R$ then enters only the analysis, the cost of an $R$
depending on $P$ and $Q$ vanishes, and the choice is governed by the first two
axes only. This requires $\lambda \gtrsim \|\Sigma_R\|_{\mathrm{op}}$ in the
analysis, which holds along the coupling $\lambda \asymp h^d$ of a shrinking
bandwidth but is incompatible with the
decaying $\lambda$ of Theorems~\ref{thm:sobolev-fixed}
and~\ref{thm:separation-rkhs}, where the kernel is fixed and
$\|\Sigma_R\|_{\mathrm{op}}$ does not shrink. The splitting cost described below
is therefore unavoidable for the general-domain class and for a fixed kernel, and
avoidable for the Sobolev class only by shrinking the bandwidth and using the
MMD test.

\paragraph{Mixture.}
For $R = (P+Q)/2$, the law of total variance gives
$2\Sigma_R = \Sigma_P + \Sigma_Q + \tfrac12(\mu_P-\mu_Q)\otimes_{\mathcal{H}}
(\mu_P-\mu_Q) \succcurlyeq \Sigma_P, \Sigma_Q$, so the dominance condition
holds with $c = 2$ unconditionally, and
$\|r\|_\infty \leq \tfrac12(\|p\|_\infty + \|q\|_\infty) \leq M$ for
$(P,Q) \in \mathcal{F}_{s,\Delta}$. For the general-domain class, $w = 2u$ and
$\|w\|^2_{L^2(R)} = 4\rho^2(P,Q)$, so the separation is in the Hellinger-type
distance and Theorem~\ref{thm:separation-rkhs} recovers the rate of
\citet{hagrass2024a_spectral_twosample}. This choice imposes no structure on
$\mathcal{X}$ and no regularity beyond that already carried by the alternative
class, and is available for both classes.

Because $R$ depends on the unknown $P$ and $Q$, it cannot be sampled directly.
With a spectrum-dependent regularizer this places the mixture in the
sample-splitting regime of all three separation theorems; with a fixed kernel,
the sharp Sobolev rate then requires $s>\tau-d/4$
(Theorem~\ref{thm:sobolev-fixed}(ii)). For the MMD test, however, $\Sigma_R$ is
never estimated, so no samples are needed and the mixture costs nothing; with a
shrinking bandwidth this removes the $(\log N/N)^{2s/d}$ term and leaves the
sharp rate $N^{-4s/(4s+d)}$ for every $s > 0$, with no restriction on
$\mathcal{X}$ (Proposition~\ref{prop:mmd}).

\paragraph{Uniform.}
If $\mathcal{X}$ is compact, take $R = \mathrm{Uniform}(\mathcal{X})$. Then
$p(x) \leq M \leq M|\mathcal{X}|\,r(x)$ pointwise gives the dominance
condition with $c = M|\mathcal{X}|$, and $\|r\|_\infty = |\mathcal{X}|^{-1}$.
Because $R$ is known, covariance samples are drawn independently at no cost to
the test budget, so all three separation theorems apply in their part~(i), with
the sharp rates holding for all $s, \theta > 0$. For the general-domain class
this is the only route to the sharp rate, since there is no bandwidth to shrink
and the constant regularizer cannot reach the decaying $\lambda$ required there;
for the Sobolev class with a shrinking bandwidth it is an alternative to the MMD
test, buying the same rate at the price of compactness. It is also the natural
choice
for a fixed kernel: since $R$ does not depend on $(P,Q)$,
condition~\eqref{eq:eigdecay} and the lower bound on
$\|\Sigma_R\|_{\mathrm{op}}$ in Theorem~\ref{thm:sobolev-fixed} are properties of
$K$ and $\mathcal{X}$ alone, and hold for the kernels listed in~(b).

For the general-domain class this choice also changes the separation measure:
$w = |\mathcal{X}|(p-q)$, so $\|w\|^2_{L^2(R)} = |\mathcal{X}|\,\|p-q\|_2^2$
and the separation is in $L^2$ rather than in the Hellinger-type distance. The
minimax lower bound of
\citet[Theorem~3.2]{hagrass2024a_spectral_twosample} is proved for a fixed
reference measure and stated in terms of $\|\cdot\|_{L^2(R)}$, so it applies
here unchanged and the test remains minimax optimal over the resulting class.

\paragraph{Goodness-of-fit: null distribution.}
When $Q = P_0$ is known, taking $R = P_0$ gives $\Sigma_Q = \Sigma_R$ and
$w = u = dP/dP_0 - 1$, so that for the general-domain class the separation
measure is $\|w\|^2_{L^2(R)} = \chi^2(P \,\|\, P_0)$, and fresh covariance
samples are free. The dominance condition
$\Sigma_P \preccurlyeq c\,\Sigma_R$, however, is typically ensured by a uniform
bound on the likelihood ratio $dP/dP_0$ over the alternative class, which is a
genuine restriction on both classes. As for the uniform choice, $R=P_0$ does not depend
on $P$, so condition~\eqref{eq:eigdecay} is a property of $K$ and $P_0$ alone,
and Theorem~\ref{thm:sobolev-fixed}(i) applies to a fixed kernel whenever it
holds.

\paragraph{Goodness-of-fit: asymmetric mixture.}
That restriction is avoided by taking $R = \varepsilon P + (1-\varepsilon)P_0$
for $\varepsilon \in (0,1/2]$, for which $\Sigma_R \succcurlyeq \varepsilon\Sigma_P$
and $\Sigma_R \succcurlyeq (1-\varepsilon)\Sigma_{P_0}$ give
$c = \varepsilon^{-1}$ unconditionally. Since
$\varepsilon(dP + dP_0) \leq dR \leq (1-\varepsilon)(dP + dP_0)$, the
separation measure of the general-domain class satisfies
\[
  2\rho^2(P, P_0) \;\leq\; \|w\|_{L^2(R)}^2 \;\leq\; 2\chi^2(P \,\|\, P_0)
\]
for every $\varepsilon \in (0,1/2]$, so that separation in the Hellinger-type
distance implies separation in $\|w\|_{L^2(R)}$, uniformly in $\varepsilon$. With a spectrum-dependent
regularizer, drawing $l$ covariance samples requires roughly $\varepsilon l$
fresh $P$-samples, split from the test budget, together with
$(1-\varepsilon)l$ free samples from $P_0$; the constraint
$\varepsilon l \lesssim n$ is what ties $\varepsilon$ to the smoothness. The
trade-off therefore matters only where $\Sigma_R$ must be estimated, that is,
for the general-domain class and for a fixed kernel; with a shrinking
bandwidth, the MMD test avoids it.

\begin{proposition}[Asymmetric mixture]\label{prop:asymmetric}
Let $0 < \alpha \leq e^{-1}$, $\delta > 0$, $Q = P_0$ known, and
$R = \varepsilon P + (1-\varepsilon)P_0$ with $\varepsilon \in (0,1/2]$ and
$\varepsilon l \lesssim n$. Then
$\Sigma_P, \Sigma_Q \preccurlyeq \varepsilon^{-1}\Sigma_R$. Suppose further that
$g_\lambda$ satisfies~\eqref{eq:g-two-sided} and, according to the case below,
that the remaining hypotheses of Theorem~\ref{thm:separation}, of
Theorem~\ref{thm:sobolev-fixed}, or of Theorem~\ref{thm:separation-rkhs} hold,
where for
Theorem~\ref{thm:sobolev-fixed} condition~\eqref{eq:eigdecay} is required of
$\Sigma_P$, uniformly over the class, and of $\Sigma_{P_0}$, instead of
$\Sigma_R$. Then, in each of the following cases, the test has power at least
$1-2\delta$ uniformly over the corresponding class.
\begin{enumerate}[label=\textnormal{(\alph*)}]
  \item \textnormal{(Sobolev class, bandwidth family.)} If $s > d/4$, the choice
  $\varepsilon = 1/2$, $\lambda \asymp h^d \asymp N^{-2d/(4s+d)}$ suffices for
  $\Delta \gtrsim N^{-4s/(4s+d)}$. If $s \leq d/4$, the choice
  $\varepsilon \asymp N^{(4s-d)/(4s+3d)}$, up to a logarithmic factor,
  $\lambda \asymp h^d \asymp \varepsilon\log N/N$, suffices for
  $\Delta \gtrsim (\sqrt{\log N}/N)^{8s/(4s+3d)}$.
  \item \textnormal{(Sobolev class, fixed kernel.)} If $s > \tau - d/4$, the
  choice $\varepsilon = 1/2$, $\lambda \asymp N^{-4\tau/(4s+d)}$ suffices for
  $\Delta \gtrsim N^{-4s/(4s+d)}$. If $s \leq \tau - d/4$, the choice
  $\varepsilon \asymp N^{(4s+d-4\tau)/(4s+d+4\tau)}$, up to a logarithmic
  factor, $\lambda \asymp \varepsilon\log N/N$, suffices for
  $\Delta \gtrsim (\sqrt{\log N}/N)^{8s/(4s+d+4\tau)}$.
  \item \textnormal{(General-domain class.)} If $\theta > 1/4$, the
  choice $\varepsilon = 1/2$, $\lambda \asymp N^{-2/(4\theta+1)}$
  suffices for $\Delta \gtrsim N^{-4\theta/(4\theta+1)}$. If
  $\theta \leq 1/4$, the choice
  $\varepsilon \asymp N^{(4\theta-1)/(4\theta+3)}$, up to a
  logarithmic factor, $\lambda \asymp \varepsilon\log N/N$, suffices for
  $\Delta \gtrsim (\sqrt{\log N}/N)^{8\theta/(4\theta+3)}$.
\end{enumerate}
\end{proposition}

Outside the optimal regime, these rates are strictly faster than those that
sample splitting gives under the symmetric mixture, namely
$(\log N/N)^{2s/d}$, $(\log N/N)^{s/\tau}$ and $(\log N/N)^{2\theta}$, but
slower than the sharp rate available when $R = P_0$ is admissible. The asymmetric
sample budget thus partially compensates for the absence of an a priori bound
on $dP/dP_0$.

\section{Adaptive test}\label{sec:adaptation}

The optimal choices of $(\lambda, h)$ in Theorems~\ref{thm:separation},
\ref{thm:sobolev-fixed} and~\ref{thm:separation-rkhs} depend on the unknown
smoothness parameter ($s$ for the Sobolev
alternative, $\theta$ for the general-domain alternative). To
obtain a fully data-driven test, we follow the aggregation approach
of \citet{hagrass2024a_spectral_twosample} and construct an adaptive
test by taking the union of individual tests over a finite grid of
$(\lambda, h)$ values. The resulting test is minimax optimal up to
a $\log\log N$ factor, simultaneously over both alternative
classes, without prior knowledge of $s$ or $\theta$. With a fixed kernel of
finite smoothness, the bandwidth grid can be dropped without losing minimax
optimality (Corollary~\ref{cor:adaptation-fixed}), although we recommend the
joint grid in practice.

Define the geometric grids
\[
  \Lambda \coloneqq \bigl\{2^j \lambda_L : 0 \leq j \leq b_\lambda
  \bigr\},
  \qquad
  \mathcal{W} \coloneqq \bigl\{2^j h_L : 0 \leq j \leq b_h\bigr\},
\]
with $\lambda_U = 2^{b_\lambda}\lambda_L$, $h_U = 2^{b_h} h_L$,
and $b_\lambda, b_h \in \mathbb{N}$. The adaptive test rejects
$H_0$ if there exists $(\lambda, h) \in \Lambda \times \mathcal{W}$
such that
\[
  \hat{\eta}_{\lambda, h} >
  \hat{q}^{\,\lambda, h}_{1 - \alpha/|\Lambda \times \mathcal{W}|},
\]
where each threshold is computed at the Bonferroni-corrected level
$\alpha/|\Lambda \times \mathcal{W}|$. The family $\mathcal{W}$ provides the
bandwidth adaptation needed for the Sobolev case along the bandwidth route of
Theorem~\ref{thm:separation}, and any fixed $h \in \mathcal{W}$ combined with
the $\Lambda$ adaptation covers the general-domain case.

\begin{theorem}[Adaptation over $\lambda$ and $h$]
\label{thm:adaptation}
Suppose the hypotheses of Theorem~\ref{thm:separation} hold for every $s > s_l$,
and those of Theorem~\ref{thm:separation-rkhs} for every $\theta > \theta_l$. Choose
\[
  \lambda_L = r_1 \tfrac{\log N}{N}, \quad \lambda_U = r_2, \quad
  h_L = r_3 \big(\tfrac{\log N}{N}\big)^{1/d}, \quad h_U = r_4,
\]
for suitable constants $r_1, r_2, r_3, r_4 > 0$ that do not depend on $N$, $s$
or $\theta$. Under these choices,
$|\Lambda| \asymp |\mathcal{W}| \asymp \log N$, and so
$|\Lambda \times \mathcal{W}| \asymp \log^2 N$. Then there exists $N_0$ such
that for all $N \geq N_0$, the adaptive test has level $\alpha$ under $H_0$
and satisfies
\[
  P_{H_1}\!\left\{
    \bigcup_{(\lambda, h) \in \Lambda \times \mathcal{W}}
    \hat{\eta}_{\lambda, h} >
    \hat{q}^{\,\lambda, h}_{1 - \alpha/|\Lambda \times \mathcal{W}|}
  \right\} \geq 1 - 2\delta,
\]
uniformly over $(P, Q)$ in the union $\bigl(\bigcup_{s > s_l}
\mathcal{F}_{s,\Delta}\bigr) \cup \bigl(\bigcup_{\theta > \theta_l}
\mathcal{P}_{\theta,\Delta}\bigr)$, provided $\Delta$ satisfies the
following rate corresponding to the alternative class containing
$(P, Q)$:
\begin{itemize}
  \item \textnormal{(Sobolev.)} For
  $(P, Q) \in \mathcal{F}_{s,\Delta}$ with $s > s_l$,
  \[
    \Delta \;\gtrsim\; \max\!\left\{
      \left(\tfrac{\log\log N}{N}\right)^{4s/(4s+d)},\;\;
      \left(\tfrac{\log N}{N}\right)^{2s/d}
    \right\}.
  \]
  \item \textnormal{(General-domain.)} For
  $(P, Q) \in \mathcal{P}_{\theta,\Delta}$ with $\theta > \theta_l$,
  \[
    \Delta \;\gtrsim\; \max\!\left\{
      \left(\tfrac{\log\log N}{N}\right)^{4\theta/(4\theta+1)},\;\;
      \left(\tfrac{\log N}{N}\right)^{2\theta}
    \right\}.
  \]
\end{itemize}
\end{theorem}
Condition~\eqref{eq:lambda-range} is not required at every point of
$\Lambda \times \mathcal{W}$. The level statement holds by exchangeability for any
pair, and the power statement is established at the single pair
$(v^*_\lambda, v^*_h)$ of Step~3 of the proof, which
satisfies~\eqref{eq:lambda-range} by construction.
In each rate, the first term corresponds to the optimal regime ($s > d/4$ or
$\theta > 1/4$), where the adaptive test matches the minimax rate up to
a $\log\log N$ factor, the price of the Bonferroni correction across the grid.
The second term comes from the sample-splitting regime of the separation
theorems and is unaffected by aggregation. For the guarantees, aggregation over
$h$ is needed only for the Sobolev class along the bandwidth route, where the
optimal bandwidth
depends on $s$; for the
general-domain class the bandwidth plays no role, and $\mathcal{W}$ contributes
only to the size of the joint grid. With a fixed kernel of finite smoothness, the Sobolev class can instead be
handled by aggregating over $\lambda$ alone, following
Theorem~\ref{thm:sobolev-fixed}. Write $\hat\eta_\lambda$ and
$\hat q^{\,\lambda}_{1-\alpha'}$ for the statistic and permutation quantile
computed with the fixed kernel at regularization $\lambda$.

\begin{corollary}[Adaptation over $\lambda$ with a fixed kernel]
\label{cor:adaptation-fixed}
Let Assumption~\ref{ass:TI} hold with a fixed kernel $K$
satisfying~\eqref{eq:Khat-lower}, and suppose the hypotheses of
Theorem~\ref{thm:sobolev-fixed} hold for every $s > s_l$, and those of
Theorem~\ref{thm:separation-rkhs} for every $\theta > \theta_l$. Consider the
test that rejects $H_0$ if
$\hat\eta_\lambda > \hat q^{\,\lambda}_{1-\alpha/|\Lambda|}$ for some
$\lambda \in \Lambda$, with $\lambda_L = r_1\log N/N$ and $\lambda_U = r_2$
for suitable constants $r_1, r_2 > 0$ that do not depend on $N$, $s$ or
$\theta$, so that
$|\Lambda| \asymp \log N$. Then there exists $N_0$ such that for all
$N \geq N_0$ the test has level $\alpha$ under $H_0$ and power at least
$1-2\delta$ uniformly over
$\bigl(\bigcup_{s > s_l}\mathcal{F}_{s,\Delta}\bigr) \cup
\bigl(\bigcup_{\theta > \theta_l}\mathcal{P}_{\theta,\Delta}\bigr)$, provided
$\Delta$ satisfies the general-domain rate of Theorem~\ref{thm:adaptation} for
$(P,Q) \in \mathcal{P}_{\theta,\Delta}$, and
\[
  \Delta \;\gtrsim\; \max\!\left\{
    \left(\tfrac{\log\log N}{N}\right)^{4s/(4s+d)},\;\;
    \left(\tfrac{\log N}{N}\right)^{s/\tau}
  \right\}
\]
for $(P,Q) \in \mathcal{F}_{s,\Delta}$ with $s > s_l$.
\end{corollary}

With a kernel of finite smoothness, a single bandwidth and $O(\log N)$
regularization parameters thus suffice for minimax optimality over both classes.
In the optimal regime the Sobolev rate is that of Theorem~\ref{thm:adaptation};
the optimal regime now requires $s > \tau - d/4$ instead of $s > d/4$,
reflecting Theorem~\ref{thm:sobolev-fixed}(ii).

In practice we nevertheless recommend the joint grid of
Theorem~\ref{thm:adaptation}. Minimax optimality concerns the least favourable
alternatives in each class, whereas for a given alternative the most powerful
bandwidth depends on its scale and shape, which aggregation over $\mathcal{W}$
lets the data select; the experiments of
\citet[Section~5]{hagrass2024a_spectral_twosample} aggregate over both
$\lambda$ and the bandwidth. The joint grid also costs
little in theory. Since $|\Lambda\times\mathcal{W}|\asymp\log^2N$ against
$|\Lambda|\asymp\log N$, the bandwidth grid changes only the constant in the
$\log\log N$ factor, and where the hypotheses of both results hold, the Sobolev
rate of Theorem~\ref{thm:adaptation} is at least as good as that of
Corollary~\ref{cor:adaptation-fixed}, since $d/4<\tau-d/4$ and $2s/d>s/\tau$.

\section{Discussion}\label{sec:discussion}

The results of this paper place the bandwidth and the regularization parameter
within a single framework, and Theorem~\ref{thm:signal} shows that both are ways of meeting
the single condition~\eqref{eq:spectral-condition} on the spectral density of
the kernel: the discrepancy $\eta$ retains a fixed fraction of $\|p-q\|_2^2$ for
every alternative in $\mathcal{F}_{s,\Delta}$ once $\lambda \lesssim m_K(1/t)$
with $t \asymp \Delta^{1/(2s)}$, that is, once $\lambda$ is below the spectral
density on the frequencies up to $1/t$; for kernels of finite smoothness the
condition is also necessary, up to constants (Proposition~\ref{prop:necessity}).
A bandwidth family has $m_K(1/h) \asymp h^d$, so a shrinking bandwidth meets the
condition at $t = h$ with $\lambda \asymp h^d$, the two parameters decaying together. A fixed kernel of finite smoothness has
$m_K(1/t) \gtrsim t^{2\tau}$, so it meets the same condition with
$\lambda \asymp t^{2\tau}$ and no bandwidth at all
(Theorem~\ref{thm:sobolev-fixed}). The two routes pay the same price in
effective dimension, $\mathcal{N}_2(\lambda) \lesssim t^{-d/2}$, the first
through the trace bound and the second through the decay
condition~\eqref{eq:eigdecay}, and in both the Sobolev rate is the balance of
$t^{2s}$ against $t^{-d/2}/N$, attained at $t \asymp N^{-2/(4s+d)}$. They differ in what they ask of the regularizer. Along the coupling $\lambda\asymp h^d$,
every regularizer obeying the two-sided bound~\eqref{eq:g-two-sided} attains the
rate, the constant one included, under which the test is MMD with the kernel at
bandwidth $h$; the convolution test of \citet{schrab2023mmd} is the member of
the family that takes this route (Proposition~\ref{prop:mmd}). With the kernel
fixed, it is the regularizer, evaluated at $\hat\Sigma_R$, that does the work
the bandwidth did, and the constant one is unavailable
(Remark~\ref{rem:regimes}). On a general domain $\lambda$ also decays with the
kernel fixed, and there too the choice of regularizer does not enter the rate:
both lower bounds on $\eta$, Theorem~\ref{thm:signal} and
Proposition~\ref{prop:signal-rkhs}, are cutoff arguments that use only the lower
half of~\eqref{eq:g-two-sided}, in Fourier frequency for the Sobolev class and in
the spectrum of $\mathcal{T}$ for the general-domain class.

Seen this way, the two diagnoses of MMD's suboptimality in the literature are
the same diagnosis. The convolution account attributes the loss to a kernel
whose spatial scale cannot adapt to the regularity of $p-q$; the regularization
account attributes it to the eigenvalues of $\mathcal{T}$ downweighting the
high-frequency coefficients of $u$. For a translation-invariant kernel on
$\mathbb{R}^d$ both describe the decay of $\widehat K$, and the two cures act on
that decay differently: shrinking the bandwidth pushes it out to frequencies of
order $h^{-1}$, while the regularizer leaves $\widehat K$ unchanged and inverts
the spectrum down to the level $\lambda$. Theorem~\ref{thm:sobolev-fixed} shows
that, for a kernel of finite smoothness, the second cure can replace the first.
MMD cannot use it: with the constant regularizer the statistic is a multiple of the
MMD statistic, so the permutation test is the MMD test whatever $\lambda$ is,
and in its analysis condition~\eqref{eq:spectral-condition} can be met only by
changing the kernel. Within the spectral regularized framework, this is why the
convolution approach of \citet{li2019optimality} and \citet{schrab2023mmd}
shrinks the bandwidth, whatever the kernel, the Laplace kernel as much as the
Gaussian, in line with the known suboptimality of the MMD test with a fixed
kernel \citep{Krishna, hagrass2024a_spectral_twosample}; an analysis through the
kernel-smoothed witness function, as in the convolution approach, ties the
frequency cutoff to the bandwidth in the same way (Remark~\ref{rem:witness}).

For a spectrum-dependent regularizer the replacement has two limits. First,
$\widehat K$ must decay slowly enough that the required $\lambda$ is usable and
that the effective dimension at that $\lambda$ is of order $t^{-d/2}$. Both hold
for a kernel of finite smoothness, where $\lambda$ is polynomially small and the
eigenvalues of $\Sigma_R$ follow the spectral density. Within the present
analysis, neither holds for the Gaussian kernel (Remark~\ref{rem:gaussian}).
There $\lambda$ is super-polynomially small, which rules out estimating
$\Sigma_R$ from polynomially many samples, since~\eqref{eq:lambda-range} needs
$l$ of order $\lambda^{-1}$ up to a logarithm; and even with unlimited
covariance samples the eigenvalues of $\Sigma_R$ decay much more slowly than the
spectral density, so that the effective dimension is bounded only by $t^{-d}$
and the rate is $N^{-2s/(2s+d)}$. In this analysis, the Gaussian bandwidth must
therefore shrink even when the regularizer is spectrum-dependent. Second,
the trace bound used along the first route holds for every $R$, whereas the
decay condition~\eqref{eq:eigdecay} used along the second must hold uniformly
over the class, which limits how widely $R$ may spread. On a general domain there is no frequency
variable through which the kernel's spectrum can be moved, and the regularizer
is the only device available.

The framework also gives a concrete reading of which device to use
(Table~\ref{tab:routes}). On $\mathbb{R}^d$ against Sobolev alternatives, the
shrinking-bandwidth route with the MMD test of Proposition~\ref{prop:mmd} needs
no covariance estimation
and no control of the spectrum of $\Sigma_R$ beyond its trace, and it makes the
mixture $R = (P+Q)/2$ available at no cost to the test budget, which removes the
$(\log N/N)^{2s/d}$ term for small $s$; its price is a kernel that changes with
$N$, and hence a grid of bandwidths for adaptation. The fixed-kernel route needs
a kernel of finite smoothness, the decay condition~\eqref{eq:eigdecay} and an
estimate of $\Sigma_R$; in exchange the kernel is chosen once, and minimax
adaptation needs only a grid of regularization parameters
(Corollary~\ref{cor:adaptation-fixed}). Under sample splitting this route is
sharp for $s > \tau - d/4$, a threshold that decreases to $d/4$ as
$\tau \downarrow d/2$, so that a rougher kernel covers a wider range of $s$.
Whenever $\Sigma_R$ is estimated, as it must be along the fixed-kernel route and
on a general domain, the choice of $R$ matters: a known $R$, such as the uniform
measure on a bounded domain or the null distribution in the goodness-of-fit
problem, avoids the splitting cost and gives the sharp rate for every
smoothness, the mixture does not, and in the goodness-of-fit problem the
asymmetric mixture of Proposition~\ref{prop:asymmetric} recovers part of the
loss. When the regime is unknown---which is the usual situation, since it is
governed by the smoothness of the alternative---the aggregated tests of
Section~\ref{sec:adaptation} span all of these at the price of a $\log\log N$
factor. Of these, we recommend in practice the test aggregating over both
$\lambda$ and $h$, even where a grid over $\lambda$ alone suffices in theory
(Section~\ref{sec:adaptation}).

Numerically, a test of this form, aggregating the spectral regularized
statistic over a grid of $(\lambda, h)$ values, was compared against the MMDAgg
test of \citet{schrab2023mmd} in
\citet[Section~5]{hagrass2024a_spectral_twosample}, across the alternatives
considered there, with the corresponding comparison for goodness-of-fit testing
in \citet{hagrass2023Gof}. The theory of \citet{hagrass2024a_spectral_twosample}
covers aggregation over a finite family of kernels for the operator-theoretic
class; the present analysis extends the guarantees to the Sobolev class, with a
bandwidth that shrinks with $N$. The power gains reported there over a range of
alternatives are consistent with Theorem~\ref{thm:signal}, in which $h$ and
$\lambda$ enter jointly through condition~\eqref{eq:resolution}, and they
support aggregating over both parameters, as in Theorem~\ref{thm:adaptation}.

Several questions remain. The splitting cost incurred by a spectrum-dependent
regularizer with $R$ depending on $P$ and $Q$ is characterized here but not
shown to be necessary; a lower bound for tests that estimate $\Sigma_R$ from the
test budget would settle whether the degradation below $s = d/4$ for a bandwidth
family, below $s = \tau - d/4$ for a fixed kernel, and below $\theta = 1/4$ on a
general domain is intrinsic. The question matters most for a fixed kernel, where
the MMD test, which removes the cost for a bandwidth family, cannot
meet~\eqref{eq:spectral-condition} for small $\Delta$
(Remark~\ref{rem:why-shrink}). The need for a shrinking bandwidth with the
Gaussian kernel is established here only within the present analysis:
a lower bound showing that no spectral regularized test with a fixed-bandwidth
Gaussian kernel and $\Sigma_R$ estimated from the sample attains a polynomial
rate would make it a theorem, and when covariance samples are unlimited it is
open whether such a test can improve on the rate $N^{-2s/(2s+d)}$ of
Remark~\ref{rem:gaussian}. It is also not known whether a fixed kernel can
attain the Sobolev rate over the whole of $\mathcal{F}_{s,\Delta}$, on which the
trace bound alone gives only $N^{-2s/(2s+\tau)}$, or whether a restriction on
the spread of $R$, such as~\eqref{eq:eigdecay}, is the price of holding the
bandwidth fixed.

Between the two extremes of the family---a constant regularizer needing no spectral
information, and a spectrum-dependent one needing $\hat\Sigma_R$ in full---it is
not known whether the covariance requirement interpolates, so that part of the
operator-theoretic optimality could be bought at a weaker estimation cost.
Finally, forming $g_\lambda^{1/2}(\hat\Sigma_R)$ over the grid
$\Lambda \times \mathcal{W}$ is expensive. Random Fourier feature
approximations of the spectral regularized test are known to preserve minimax
optimality when the number of features is large enough, for a kernel held fixed
with the sample size \citep{mukherjee2025minimaxoptimalkerneltwosample}, where a
data-adaptive choice of the regularization parameter and the kernel is also
proposed for practical use. Under the coupling of Theorem~\ref{thm:separation} the kernel
itself varies with $N$, so a separate feature map is required at each bandwidth,
and whether the joint grid can be approximated at a cost that does not scale
with $|\mathcal{W}|$ remains open. Along the fixed-kernel route of
Theorem~\ref{thm:sobolev-fixed} the difficulty does not arise: a single feature
map, with frequencies drawn from the normalized spectral density of $K$, serves
the whole grid of Corollary~\ref{cor:adaptation-fixed}, which makes that test
the natural choice when the joint grid is too expensive. Whether the guarantees
of \citet{mukherjee2025minimaxoptimalkerneltwosample} extend to the Sobolev
class along this route is a natural next question.

\appendix

\section{Proofs}\label{app:proofs}

\subsection{Abstract power lemma}\label{app:abstract}

\begin{lemma}[Abstract power criterion]\label{lem:abstract}
Let $V \geq \mathrm{Var}(\hat{T})$, and suppose that for some constant $c > 0$,
\begin{align*}
  &P_{H_1}\!\left\{\E(\hat{T}) \geq (1+c)
  \sqrt{\frac{V}{\delta}}\right\}
  \geq 1-\frac{\delta}{2}, \\
  &P_{H_1}\!\left\{r \leq c\sqrt{\frac{V}{\delta}}\right\}
  \geq 1-\frac{\delta}{2},
\end{align*}
where $\E(\cdot)$ and $\mathrm{Var}(\cdot)$ may be conditional on a
subset of the samples, in which case $V$ may depend on that subset. Then
$P_{H_1}\{\hat{T} > r\} \geq 1-2\delta$.
\end{lemma}

\begin{proof}
\begin{align*}
P_{H_1}\{\hat{T} > r\}
&\geq P_{H_1}\!\left(\left\{\hat{T} > \mathbb{E}(\hat{T})-\sqrt{\frac{V}{\delta}}\right\} \cap \left\{r \leq \mathbb{E}(\hat{T}) - \sqrt{\frac{V}{\delta}}\right\}\right) \\
&\geq 1 - P_{H_1}\!\left\{\hat{T} \leq \mathbb{E}(\hat{T})-\sqrt{\frac{V}{\delta}}\right\} - P_{H_1}\!\left\{r > \mathbb{E}(\hat{T}) - \sqrt{\frac{V}{\delta}}\right\} \\
&\overset{(*)}{\geq} 1 - P_{H_1}\!\left\{\left|\hat{T}-\mathbb{E}(\hat{T})\right| \geq \sqrt{\frac{V}{\delta}}\right\} - \delta \\
&\geq 1-2\delta,
\end{align*}
where the last inequality follows by Chebyshev's inequality, since
$P_{H_1}\{|\hat T - \E(\hat T)| \geq \sqrt{V/\delta}\} \leq
\delta\,\mathrm{Var}(\hat T)/V \leq \delta$. Step $(*)$ follows by bounding
the second term using
\begin{align*}
P_{H_1}\!\left\{r \leq \mathbb{E}(\hat{T}) - \sqrt{\frac{V}{\delta}}\right\}
&\geq P_{H_1}\!\left\{c\sqrt{\frac{V}{\delta}} \leq \mathbb{E}(\hat{T})
    - \sqrt{\frac{V}{\delta}}\right\} \cap \left\{r \leq c\sqrt{\frac{V}{\delta}}\right\} \\
&\geq 1 - P_{H_1}\!\left\{\mathbb{E}(\hat{T}) \leq (1+c)\sqrt{\frac{V}{\delta}}\right\}
    - P_{H_1}\!\left\{r \geq c\sqrt{\frac{V}{\delta}}\right\} \\
&\geq 1-\delta,
\end{align*}
where the last step applies both hypotheses of the lemma, each contributing $\delta/2$.
\end{proof}

\subsection{Proof of Proposition~\ref{prop:var-bound}}

Defining $A(x) \coloneqq \B\Sigma_{R,\lambda}^{-1/2}k(\cdot,x)$,
the estimator $\hat{\eta}$ can be written as
\[
  \hat{\eta}
  = \frac{1}{n(n-1)} \sum_{i \neq j} \langle A(X_i), A(X_j)\rangle_{\mathcal{H}}
  + \frac{1}{m(m-1)} \sum_{i \neq j} \langle A(Y_i), A(Y_j)\rangle_{\mathcal{H}}
  - \frac{2}{nm} \sum_{i,j} \langle A(X_i), A(Y_j)\rangle_{\mathcal{H}}.
\]
Note that
\[
  \zeta \coloneqq \E(\hat{\eta} \mid \mathbb{Z}_l)
  = \bigl\|\B\,\Sigma_{R,\lambda}^{-1/2}(\mu_P - \mu_Q)\bigr\|_{\mathcal{H}}^2,
\]
and define
\[
  a(x) \coloneqq \B\,\Sigma_{R,\lambda}^{-1/2}\bigl(k(\cdot,x) - \mu_P\bigr),
  \qquad
  b(x) \coloneqq \B\,\Sigma_{R,\lambda}^{-1/2}\bigl(k(\cdot,x) - \mu_Q\bigr).
\]
Then $\hat{\eta}$ can be decomposed as
\begin{align*}
  \hat{\eta}
  &= \frac{1}{n(n-1)}\sum_{i\neq j}\langle a(X_i), a(X_j)\rangle_{\mathcal{H}}
   + \frac{1}{m(m-1)}\sum_{i\neq j}\langle b(Y_i), b(Y_j)\rangle_{\mathcal{H}} \\
  &\quad + \frac{2}{m}\sum_{i=1}^{m}\Bigl\langle b(Y_i),\,
       \B\,\Sigma_{R,\lambda}^{-1/2}(\mu_Q - \mu_P)\Bigr\rangle_{\mathcal{H}}
   - \frac{2}{nm}\sum_{i,j}\langle a(X_i), b(Y_j)\rangle_{\mathcal{H}} \\
  &\quad - \frac{2}{n}\sum_{i=1}^{n}\Bigl\langle a(X_i),\,
       \B\,\Sigma_{R,\lambda}^{-1/2}(\mu_Q - \mu_P)\Bigr\rangle_{\mathcal{H}}
   + \zeta.
\end{align*}
Writing $\hat{\eta} - \zeta = \mathsf{T}_1 + \mathsf{T}_2 + \mathsf{T}_3
+ \mathsf{T}_4 + \mathsf{T}_5$ for the five terms above (in order),
and conditioning on $\mathbb{Z}_l$, we have
\[
  \mathrm{Var}(\hat{\eta} \mid \mathbb{Z}_l)
  = \mathbb{E}\!\bigl[(\hat{\eta} - \zeta)^2 \mid \mathbb{Z}_l\bigr]
  \leq
  5\sum_{k=1}^{5}\mathbb{E}\bigl[\mathsf{T}_k^2 \mid \mathbb{Z}_l\bigr],
\]
where the last step follows by the Cauchy--Schwarz inequality.
We bound each term using
\citet[Lemmas~F.4--F.6]{hagrass2024a_spectral_twosample}, as follows.
\begin{align*}
  \mathbb{E}\bigl[\mathsf{T}_1^2 \mid \mathbb{Z}_l\bigr]
  &\leq \frac{4}{n^2}\|\B\|_{\mathrm{op}}^4
        \|\Sigma_{R,\lambda}^{-1/2}\Sigma_P
        \Sigma_{R,\lambda}^{-1/2}\|_{\mathcal{L}^2(\mathcal{H})}^2 \\
  &\leq \frac{16c^2}{n^2}\|\B\|_{\mathrm{op}}^4
        \|\Sigma_{R,\lambda}^{-1/2}\Sigma_R
        \Sigma_{R,\lambda}^{-1/2}\|_{\mathcal{L}^2(\mathcal{H})}^2,
\end{align*}
where the second inequality uses $\Sigma_P \preccurlyeq c\,\Sigma_R$.
Similarly,
\[
  \mathbb{E}\bigl[\mathsf{T}_2^2 \mid \mathbb{Z}_l\bigr]
  \leq \frac{16c^2}{m^2}\|\B\|_{\mathrm{op}}^4
       \|\Sigma_{R,\lambda}^{-1/2}\Sigma_R
       \Sigma_{R,\lambda}^{-1/2}\|_{\mathcal{L}^2(\mathcal{H})}^2,
\]
using $\Sigma_Q \preccurlyeq c\,\Sigma_R$. For the linear terms,
\begin{align*}
  \mathbb{E}\bigl[\mathsf{T}_3^2 \mid \mathbb{Z}_l\bigr]
  &\leq \frac{4}{m}\|\B\|_{\mathrm{op}}^2
        \|\Sigma_{R,\lambda}^{-1/2}\Sigma_Q
        \Sigma_{R,\lambda}^{-1/2}\|_{\mathrm{op}}
        \bigl\|\B\,\Sigma_{R,\lambda}^{-1/2}
        (\mu_P - \mu_Q)\bigr\|_{\mathcal{H}}^2 \\
  &\leq \frac{16c\,\zeta}{m}\|\B\|_{\mathrm{op}}^2,
\end{align*}
where the second inequality uses $\Sigma_Q \preccurlyeq c\,\Sigma_R$
and the definition of $\zeta$. By an analogous argument,
\[
  \mathbb{E}\bigl[\mathsf{T}_5^2 \mid \mathbb{Z}_l\bigr]
  \leq \frac{16c\,\zeta}{n}\|\B\|_{\mathrm{op}}^2,
\]
using $\Sigma_P \preccurlyeq c\,\Sigma_R$. Finally,
\[
  \mathbb{E}\bigl[\mathsf{T}_4^2 \mid \mathbb{Z}_l\bigr]
  \leq \frac{4}{nm}\|\B\|_{\mathrm{op}}^4
       \|\Sigma_{R,\lambda}^{-1/2}\Sigma_R
       \Sigma_{R,\lambda}^{-1/2}\|_{\mathcal{L}^2(\mathcal{H})}^2.
\]
Denote $\|\Sigma_{R,\lambda}^{-1/2}\Sigma_R\Sigma_{R,\lambda}^{-1/2}
\|_{\mathcal{L}^2(\mathcal{H})} \eqqcolon \mathcal{N}_2(\lambda)$.
Combining all five bounds and applying $\frac{1}{nm} \leq
\frac{1}{2}(\frac{1}{n^2} + \frac{1}{m^2})$ to absorb the cross term,
\begin{align*}
  \mathrm{Var}(\hat{\eta} \mid \mathbb{Z}_l)
  &\leq 5\left(\frac{16c^2}{n^2} + \frac{16c^2}{m^2}
    + \frac{4}{nm}\right)\|\B\|_{\mathrm{op}}^4\,\mathcal{N}_2^2(\lambda)
    + 5\left(\frac{16c}{m} + \frac{16c}{n}\right)\|\B\|_{\mathrm{op}}^2\,\zeta \\
  &\leq 5(16c^2 + 2)
    \left(\frac{1}{n^2}+\frac{1}{m^2}\right)
    \|\B\|_{\mathrm{op}}^4\,\mathcal{N}_2^2(\lambda)
    + 80c\left(\frac{1}{m} + \frac{1}{n}\right)\|\B\|_{\mathrm{op}}^2\,\zeta.
\end{align*}
Using $m \leq n \leq Dm$, which gives $\frac{1}{n^2}+\frac{1}{m^2}
\leq \frac{C_D}{(n+m)^2}$ and $\frac{1}{m}+\frac{1}{n} \leq
\frac{C_D}{n+m}$ for some constant $C_D$ depending only on $D$,
we conclude
\[
  \mathrm{Var}(\hat{\eta} \mid \mathbb{Z}_l)
  \leq
  \frac{C_{c,D}\,\|\B\|_{\mathrm{op}}^4\,\mathcal{N}_2^2(\lambda)}
  {(n+m)^2}
  +
  \frac{C_{c,D}\,\|\B\|_{\mathrm{op}}^2\,\zeta}
  {n+m},
\]
where $C_{c,D} = 5(16c^2+2)C_D$ for the first term and $80c\,C_D$
for the second. \hfill$\square$

\subsection{Proof of Lemma~\ref{lem:threshold}}

Define $a(x)$ and $b(x)$ as in the proof of Proposition~\ref{prop:var-bound}.
By \cite[Equation~59]{permutations}, conditionally on $\mathbb{Z}_l$, there
exists an absolute constant $C_6 > 0$ such that
\begin{equation}\label{eq:qs-bound}
  \hat{q}_{1-\alpha} \leq C_6\, I \log\frac{1}{\alpha}
  \quad \text{almost surely},
\end{equation}
where
\begin{align*}
  I^2 &\coloneqq
    \frac{1}{n^2(n-1)^2}\sum_{i\neq j}
      \langle a(X_i), a(X_j)\rangle_{\mathcal{H}}^2
    + \frac{1}{m^2(m-1)^2}\sum_{i\neq j}
      \langle a(Y_i), a(Y_j)\rangle_{\mathcal{H}}^2 \\
  &\quad
    + \frac{2}{n^2 m^2}\sum_{i,j}
      \langle a(X_i), a(Y_j)\rangle_{\mathcal{H}}^2.
\end{align*}
Writing $a(Y_j) = b(Y_j) + \B\Sigma_{R,\lambda}^{-1/2}(\mu_Q - \mu_P)$
and applying $\bigl(\sum_{k=1}^K a_k\bigr)^2 \leq K\sum_{k=1}^K a_k^2$,
we expand $I^2$ into terms of the same type as those in the variance
decomposition of Proposition~\ref{prop:var-bound}. Applying
\citep[Lemmas~F.4--F.6]{hagrass2024a_spectral_twosample} to bound the
conditional expectation of each term given $\mathbb{Z}_l$, and using
$m \leq n \leq Dm$, we obtain
\begin{equation}\label{eq:I2-bound}
  \mathbb{E}\bigl[I^2 \mid \mathbb{Z}_l\bigr]
  \lesssim
  \frac{\|\B\|_{\mathrm{op}}^4\,\mathcal{N}_2^2(\lambda)}{(n+m)^2}
  + \frac{\|\B\|_{\mathrm{op}}^2\,\zeta}{(n+m)^2}
  + \frac{\zeta^2}{(n+m)^2}.
\end{equation}
Since $(n+m)^{-2} \leq (n+m)^{-1}$, \eqref{eq:I2-bound} gives
$\mathbb{E}\bigl[I^2 \mid \mathbb{Z}_l\bigr] \lesssim V$ by the
definition~\eqref{eq:V-def} of $V$, with $N = n+m$.
By Markov's inequality, for some constant $C > 0$,
\[
  P_{H_1}\!\left\{I^2 \geq
  \frac{C\,V}{\delta}
  \;\middle|\; \mathbb{Z}_l\right\}
  \leq \frac{\delta}{2}.
\]
Combined with~\eqref{eq:qs-bound}, setting
$c_1 \coloneqq C_6\log(1/\alpha)\sqrt{C}$ gives
\[
  P_{H_1}\!\left\{\hat{q}_{1-\alpha} \leq
  c_1\sqrt{\frac{V}{\delta}}\right\}
  \geq 1-\frac{\delta}{2}. \qedhere
\]

\subsection{Proof of Theorem~\ref{thm:power}}

Throughout, we work on the event $\mathcal{E}$ of~\eqref{eq:event}, on which
$\|\B\|_{\mathrm{op}} \vee \|\B_0\|_{\mathrm{op}} \leq C_u$ and
$\|\B^{-1}\|_{\mathrm{op}} \leq C_\ell^{-1}$.

By Lemma~\ref{lem:abstract}, applied with $\hat T = \hat\eta$,
$r = \hat q_{1-\alpha}$ and $V$ as in~\eqref{eq:V-def}, which satisfies
$V \geq \mathrm{Var}(\hat\eta\mid\mathbb{Z}_l)$ by
Proposition~\ref{prop:var-bound}, and by Lemma~\ref{lem:threshold}, it
suffices to verify
\begin{equation}\label{eq:mean-cond}
  P_{H_1}\!\left\{\zeta \geq
  \frac{(1+c_1)\sqrt{C_0}}{\sqrt{\delta}}
  \left(
    \frac{\|\B\|_{\mathrm{op}}^2\,\mathcal{N}_2(\lambda)}{N}
        +\frac{\|\B\|_{\mathrm{op}}\sqrt{\zeta}}{\sqrt{N}}
    +\frac{\zeta}{N}
  \right)\right\} \geq 1-\frac{\delta}{2},
\end{equation}
where we have used $\sqrt{a_1 + a_2 + a_3} \leq \sqrt{a_1} + \sqrt{a_2} + \sqrt{a_3}$
to bound $\sqrt{V/\delta}$.
We establish~\eqref{eq:mean-cond} on $\mathcal{E}$. For
$N \geq N_0 \coloneqq 3(1+c_1)\sqrt{C_0}/\sqrt{\delta}$, the last term is at
most $\zeta/3$, so it is sufficient that $\zeta/3$ dominates each of the other
two terms, i.e.,
\begin{equation}\label{eq:dom1}
  \frac{\zeta}{3} \geq
  \frac{(1+c_1)\sqrt{C_0}\,\|\B\|_{\mathrm{op}}^2\,\mathcal{N}_2(\lambda)}{\sqrt{\delta}\,N},
\end{equation}
\begin{equation}\label{eq:dom2}
  \frac{\zeta}{3} \geq
  \frac{(1+c_1)\sqrt{C_0}\,\|\B\|_{\mathrm{op}}\sqrt{\zeta}}{\sqrt{\delta N}}.
\end{equation}
Condition~\eqref{eq:dom2} is equivalent to
$\zeta \geq \frac{9(1+c_1)^2C_0\|\B\|_{\mathrm{op}}^2}{\delta N}$.
On $\mathcal{E}$ we have $\|\B\|_{\mathrm{op}}^2 \leq C_u^2$, and since
$\delta \leq 1$, both~\eqref{eq:dom1} and~\eqref{eq:dom2} are implied by
\begin{equation}\label{eq:zeta-suff}
  \zeta \geq \frac{C\,\bigl(\mathcal{N}_2(\lambda) \vee 1\bigr)}{\delta N},
\end{equation}
for a constant $C = C(c_1, C_0, C_u) > 0$.
It remains to relate $\zeta$ to $\eta$. Writing
$g_\lambda^{1/2}(\Sigma_R) = \B_0\,\Sigma_{R,\lambda}^{-1/2}$ and
$\Sigma_{R,\lambda}^{-1/2} = \B^{-1}\B\,\Sigma_{R,\lambda}^{-1/2}$, we obtain
\begin{equation}\label{eq:zeta-lb}
  \eta = \bigl\|\B_0\,\Sigma_{R,\lambda}^{-1/2}(\mu_P-\mu_Q)\bigr\|_{\mathcal{H}}^2
  \leq \|\B_0\|_{\mathrm{op}}^2\,\|\B^{-1}\|_{\mathrm{op}}^2\,
  \bigl\|\B\,\Sigma_{R,\lambda}^{-1/2}(\mu_P-\mu_Q)\bigr\|_{\mathcal{H}}^2
  \leq \frac{C_u^2}{C_\ell^2}\,\zeta
\end{equation}
on $\mathcal{E}$, by the definition of $\zeta$ in
Proposition~\ref{prop:var-bound}. Substituting~\eqref{eq:zeta-lb}
into~\eqref{eq:zeta-suff}, condition~\eqref{eq:zeta-suff} holds on
$\mathcal{E}$ provided
\[
  \eta \geq \frac{C_u^2\,C\,\bigl(\mathcal{N}_2(\lambda) \vee 1\bigr)}
  {C_\ell^2\,\delta N},
\]
which is exactly condition~\eqref{cond:eta} for
$c_2 \geq C_u^2 C/C_\ell^2$. Therefore~\eqref{eq:mean-cond} holds on
$\mathcal{E}$, and since $P(\mathcal{E}) \geq 1-\delta/2$, it holds with
probability at least $1-\delta/2$. Combined with
Lemma~\ref{lem:threshold} and Lemma~\ref{lem:abstract}, this yields
$P_{H_1}\{\hat{\eta} > \hat{q}_{1-\alpha}\}\geq 1-2\delta$. \qed

\subsection{Sufficient conditions for the event \texorpdfstring{$\mathcal{E}$}{E}}\label{app:event}

\begin{lemma}[General regularizer]\label{lem:B-general}
Suppose $g_\lambda$ satisfies~\eqref{eq:g-two-sided} and $\lambda$
satisfies~\eqref{eq:lambda-range}. Then $\mathcal{E}(C'_u, C'_\ell)$ has
probability at least $1-\delta/2$, where $C'_u = \sqrt{2}\,C_u$ and
$C'_\ell = \sqrt{2/3}\,C_\ell$.
\end{lemma}

\begin{proof}
Let $\mathcal{M} \coloneqq \hat{\Sigma}_{R,\lambda}^{-1/2}\Sigma_{R,\lambda}^{1/2}$
and write $\B = \bigl[g_\lambda^{1/2}(\hat\Sigma_R)
\hat\Sigma_{R,\lambda}^{1/2}\bigr]\mathcal{M}$. Since
\[
  \bigl\|g_\lambda^{1/2}(\hat\Sigma_R)\hat\Sigma_{R,\lambda}^{1/2}
  \bigr\|_{\mathrm{op}}^2
  = \sup_{x \,\in\, \mathrm{spec}(\hat\Sigma_R)} g_\lambda(x)(x+\lambda)
  \leq C_u^2
\]
by the upper bound in~\eqref{eq:g-two-sided}, since
$\mathrm{spec}(\hat\Sigma_R) \subset [0,2\kappa]$, we obtain
$\|\B\|_{\mathrm{op}} \leq C_u\|\mathcal{M}\|_{\mathrm{op}}$. The same
computation with $\Sigma_R$ in place of $\hat\Sigma_R$ gives
$\|\B_0\|_{\mathrm{op}} \leq C_u$ deterministically. For the inverse, the
lower bound in~\eqref{eq:g-two-sided} ensures $g_\lambda(0) \neq 0$, so that
$g_\lambda^{1/2}(\hat\Sigma_R)$ is invertible, and writing
$\B^{-1} = \mathcal{M}^{-1}\bigl[\hat\Sigma_{R,\lambda}^{-1/2}
g_\lambda^{-1/2}(\hat\Sigma_R)\bigr]$ together with
\[
  \bigl\|\hat\Sigma_{R,\lambda}^{-1/2}g_\lambda^{-1/2}(\hat\Sigma_R)
  \bigr\|_{\mathrm{op}}^2
  = \sup_{x \,\in\, \mathrm{spec}(\hat\Sigma_R)}
  \frac{1}{g_\lambda(x)(x+\lambda)} \leq C_\ell^{-2}
\]
gives $\|\B^{-1}\|_{\mathrm{op}} \leq C_\ell^{-1}
\|\mathcal{M}^{-1}\|_{\mathrm{op}}$. It therefore suffices to control
$\mathcal{M}$. Noting that $\sup_{x,y}k(x,y) \leq \kappa$, the assumed range on $\lambda$ together with
\citet[Lemma~B.2\emph{(ii)}]{sriperumbudur2022approx} gives
\[
  P\!\left\{\|\mathcal{M}\|_{\mathrm{op}}^2 \geq 2\right\} \leq \frac{\delta}{4},
  \qquad
  P\!\left\{\|\mathcal{M}^{-1}\|_{\mathrm{op}}^2 \geq \frac{3}{2}\right\}
  \leq \frac{\delta}{4},
\]
and the claim follows by a union bound.
\end{proof}

\begin{lemma}[Constant regularizer]\label{lem:B-flat}
Suppose $g_\lambda \equiv \lambda^{-1}$ on $[0,2\kappa]$, and that
$\lambda \geq b\,\|\Sigma_R\|_{\mathrm{op}}$ for some constant $b > 0$. Then
$\mathcal{B} = \mathcal{B}_0 = \lambda^{-1/2}\Sigma_{R,\lambda}^{1/2}$ is
deterministic, and $\mathcal{E}(\sqrt{1+b^{-1}},\,1)$ holds surely. In
particular, no condition on the number of covariance samples $l$ is required.
\end{lemma}

\begin{proof}
If $g_\lambda \equiv \lambda^{-1}$ on an interval containing the spectra of
$\Sigma_R$ and $\hat\Sigma_R$, then
$g_\lambda^{1/2}(\hat\Sigma_R) = g_\lambda^{1/2}(\Sigma_R) = \lambda^{-1/2}I$,
whence $\B = \B_0 = \lambda^{-1/2}\Sigma_{R,\lambda}^{1/2}$, which involves
$\hat\Sigma_R$ nowhere and is therefore deterministic. Consequently
\[
  \|\B\|_{\mathrm{op}}^2
  = \frac{\|\Sigma_R\|_{\mathrm{op}} + \lambda}{\lambda}
  \leq 1 + b^{-1},
  \qquad
  \|\B^{-1}\|_{\mathrm{op}}^2
  = \sup_{x \,\in\, \mathrm{spec}(\Sigma_R)}\frac{\lambda}{x+\lambda}
  \leq 1,
\]
where the first bound uses $\lambda \geq b\,\|\Sigma_R\|_{\mathrm{op}}$.
\end{proof}

\subsection{Proof of Theorem~\ref{thm:signal}}

Write $v \coloneqq \mu_P - \mu_Q$. The proof proceeds as follows.

\medskip
We first show that, for every $f \in \mathcal{H}$,
\begin{equation}\label{eq:var-lb}
  \bigl\|\Sigma_{R,\lambda}^{-1/2}v\bigr\|_{\mathcal{H}}^2
  \;\geq\; \frac{1}{M}\Bigl(\|\psi\|_2^2 - \|\psi - f\|_2^2
  - \frac{\lambda}{M}\,\|f\|_{\mathcal{H}}^2\Bigr).
\end{equation}
For any $f \in \mathcal{H}$, completing the square gives
\[
  2\langle f, v\rangle_{\mathcal{H}}
  - \langle \Sigma_{R,\lambda}f, f\rangle_{\mathcal{H}}
  = \bigl\|\Sigma_{R,\lambda}^{-1/2}v\bigr\|_{\mathcal{H}}^2
  - \bigl\|\Sigma_{R,\lambda}^{1/2}f - \Sigma_{R,\lambda}^{-1/2}v\bigr\|_{\mathcal{H}}^2
  \;\leq\; \bigl\|\Sigma_{R,\lambda}^{-1/2}v\bigr\|_{\mathcal{H}}^2 .
\]
By the reproducing property,
$\langle f, v\rangle_{\mathcal{H}} = \int f\,d(P-Q) = \langle f, \psi\rangle_{L^2}$,
and, since $\Sigma_R$ is the centered covariance operator and
$\|r\|_\infty \leq M$,
\[
  \langle \Sigma_{R,\lambda}f, f\rangle_{\mathcal{H}}
  = \mathrm{Var}_R(f) + \lambda\|f\|_{\mathcal{H}}^2
  \;\leq\; M\|f\|_2^2 + \lambda\|f\|_{\mathcal{H}}^2 .
\]
Applying the first display to $f/M$ and using
$2\langle f,\psi\rangle_{L^2} - \|f\|_2^2 = \|\psi\|_2^2 - \|\psi - f\|_2^2$
gives~\eqref{eq:var-lb}.

\medskip
Define $f_t$ through its Fourier transform,
\[
  \widehat{f_t}(\omega) \coloneqq \widehat\psi(\omega)\,
  \mathbf{1}\{\|\omega\|_2 \leq 1/t\}, \qquad \omega \in \mathbb{R}^d .
\]
Since $\psi$ is real, $\widehat{f_t}(-\omega) = \overline{\widehat{f_t}(\omega)}$,
so $f_t$ is real-valued, and since $\widehat{f_t}$ is square-integrable with
compact support, it is integrable, so that
$f_t \in C(\mathbb{R}^d) \cap L^2(\mathbb{R}^d)$. By~\eqref{eq:native} and the
definition~\eqref{eq:mK} of $m_K$,
\[
  \|f_t\|_{\mathcal{H}}^2
  = \frac{1}{(2\pi)^d}\int_{\|\omega\|_2 \leq 1/t}
    \frac{|\widehat\psi(\omega)|^2}{\widehat K(\omega)}\,d\omega
  \;\leq\; \frac{\|\psi\|_2^2}{m_K(1/t)} < \infty ,
\]
so $f_t \in \mathcal{H}$, and by Plancherel's identity and
$\psi \in \mathcal{S}^s_d(L)$,
\[
  \|\psi - f_t\|_2^2
  = \frac{1}{(2\pi)^d}\int_{\|\omega\|_2 > 1/t}|\widehat\psi(\omega)|^2\,d\omega
  \;\leq\; \frac{t^{2s}}{(2\pi)^d}\int_{\mathbb{R}^d}
    \|\omega\|_2^{2s}|\widehat\psi(\omega)|^2\,d\omega
  \;\leq\; L^2t^{2s}.
\]
The function $f_t$ is the convolution of $\psi$ with the ideal low-pass filter
at scale $t$; its RKHS norm is controlled by the smallest value of the spectral
density on the retained band, which is where $m_K(1/t)$ enters.

Taking $f = f_t$ in~\eqref{eq:var-lb} and using both conditions
in~\eqref{eq:resolution}, yields that 
\[
  \bigl\|\Sigma_{R,\lambda}^{-1/2}v\bigr\|_{\mathcal{H}}^2
  \;\geq\; \frac{1}{M}\Bigl(\|\psi\|_2^2 - L^2t^{2s}
  - \frac{\lambda}{M\,m_K(1/t)}\,\|\psi\|_2^2\Bigr)
  \;\geq\; \frac{1}{M}\Bigl(1 - \frac14 - \frac14\Bigr)\|\psi\|_2^2
  \;=\; \frac{\|\psi\|_2^2}{2M}.
\]
Finally, $g_\lambda^{1/2}(\Sigma_R) = \B_0\,\Sigma_{R,\lambda}^{-1/2}$ with
$\B_0$ as in Section~\ref{sec:power-criterion}, and by~\eqref{eq:g-lower},
\[
  \|\B_0^{-1}\|_{\mathrm{op}}^2
  = \sup_{x \,\in\, \mathrm{spec}(\Sigma_R)}\frac{1}{g_\lambda(x)(x+\lambda)}
  \leq C_\ell^{-2},
\]
so that
$\|\Sigma_{R,\lambda}^{-1/2}v\|_{\mathcal{H}}
= \|\B_0^{-1}g_\lambda^{1/2}(\Sigma_R)v\|_{\mathcal{H}} \leq C_\ell^{-1}\eta^{1/2}$.
Therefore
$\eta \geq C_\ell^2\,\|\Sigma_{R,\lambda}^{-1/2}v\|_{\mathcal{H}}^2
\geq \frac{C_\ell^2}{2M}\|\psi\|_2^2$. \qed

\subsection{Proof of Proposition~\ref{prop:necessity}}\label{app:necessity}

By the upper bound in~\eqref{eq:g-two-sided}, $g_\lambda(x) \leq C_u^2/(x+\lambda)
\leq C_u^2/\lambda$ for $x \in [0,2\kappa]$, which contains the spectrum of
$\Sigma_R$; hence $g_\lambda(\Sigma_R) \preccurlyeq (C_u^2/\lambda)\,I$ and,
writing $\psi = p - q$,
\begin{equation}\label{eq:eta-upper}
  \eta \;\leq\; \frac{C_u^2}{\lambda}\,\|\mu_P-\mu_Q\|_{\mathcal{H}}^2
  \;=\; \frac{C_u^2}{\lambda}\,\frac{1}{(2\pi)^d}\int_{\mathbb{R}^d}
  \widehat K(\omega)\,|\widehat\psi(\omega)|^2\,d\omega ,
\end{equation}
where the identity is the Fourier form of
$\|\mu_P-\mu_Q\|_{\mathcal{H}}^2 = \iint K(x-y)\psi(x)\psi(y)\,dx\,dy$. The
bound holds whatever the reference measure $R$, and it remains to construct an
alternative in $\mathcal{F}_{s,\Delta}$ whose Fourier transform lives on
frequencies of order $\Delta^{-1/(2s)}$.

Let $a \coloneqq (2d^{3/2})^{-1}$ and
\[
  \chi(x) \coloneqq c_\chi \prod_{j=1}^{d}
  \Bigl(\frac{\sin(a x_j)}{a x_j}\Bigr)^{2d},
  \qquad x \in \mathbb{R}^d,
\]
with each factor equal to one at $x_j = 0$ and $c_\chi > 0$ chosen so that
$\|\chi\|_2 = 1$. Since the Fourier transform of $t \mapsto \sin(at)/(at)$ is
$(\pi/a)\mathbf{1}\{|\xi| \leq a\}$, that of its $2d$-th power is a $2d$-fold
convolution supported in $[-2da, 2da] = [-d^{-1/2}, d^{-1/2}]$, so
$\widehat\chi$ is supported in the cube $[-d^{-1/2}, d^{-1/2}]^d$, which lies in
the unit ball. Moreover, $|\sin(t)/t| \leq 2/(1+|t|)$ gives
$|\chi(x)| \leq C_\chi(1+\|x\|_2)^{-(d+1)}$ for a constant $C_\chi$ depending
only on $d$. For $\rho \geq 2$ and $b > 0$, let
\[
  \psi(x) \coloneqq b\,\chi(x)\cos(\rho x_1),
  \qquad\text{so that}\qquad
  \widehat\psi(\omega) = \tfrac{b}{2}\bigl(\widehat\chi(\omega-\rho e_1)
  + \widehat\chi(\omega+\rho e_1)\bigr).
\]
The two terms have disjoint supports, contained in
$\{\rho-1 \leq \|\omega\|_2 \leq \rho+1\}$, so that
$\int\psi = \widehat\psi(0) = 0$ and $\|\psi\|_2^2 = b^2/2$. Take
$b = (2\Delta)^{1/2}$, so that $\|\psi\|_2^2 = \Delta$, and
$\rho = \frac12(L^2/\Delta)^{1/(2s)}$. Since $\rho + 1 \leq 2\rho$,
\[
  \frac{1}{(2\pi)^d}\int_{\mathbb{R}^d}\|\omega\|_2^{2s}
  |\widehat\psi(\omega)|^2\,d\omega
  \;\leq\; (2\rho)^{2s}\,\|\psi\|_2^2 \;=\; L^2,
\]
that is, $\psi \in \mathcal{S}^s_d(L)$. Let $p_0$ be a density with
$\|p_0\|_\infty \leq M/2$ and $p_0(x) \geq c_0(1+\|x\|_2)^{-(d+1)}$ for some
$c_0 > 0$, such as a suitably scaled multivariate Cauchy density, and set
$p = p_0 + \psi/2$ and $q = p_0 - \psi/2$. Both integrate to one, and since
$|\psi(x)| \leq b\,C_\chi(1+\|x\|_2)^{-(d+1)}$, both are nonnegative and
bounded by $M$ once $b$ is small enough. Hence there is $\Delta_0 > 0$,
depending only on $d$, $s$, $L$ and $M$, such that for $\Delta \leq \Delta_0$
we have $\rho \geq 2$ and $(P,Q) \in \mathcal{F}_{s,\Delta}$ with
$\|p-q\|_2^2 = \Delta$.

For this alternative, \eqref{eq:eta-upper}, the support of $\widehat\psi$, the
upper bound on $\widehat K$, and $1+(\rho-1)^2 \geq (1+\rho^2)/4$ for
$\rho \geq 2$ give
\[
  \eta \;\leq\; \frac{C_u^2}{\lambda}
  \max_{\rho-1\leq\|\omega\|_2\leq\rho+1}\widehat K(\omega)\;\|\psi\|_2^2
  \;\leq\; \frac{C_u^2}{\lambda}\,4^\tau C_K\,(1+\rho^2)^{-\tau}\,\|\psi\|_2^2 .
\]
Since $\rho^* \coloneqq (4L^2/\Delta)^{1/(2s)} = 2\cdot4^{1/(2s)}\rho$, the lower
bound on $\widehat K$ gives
$m_K(\rho^*) \geq c_K(1+4^{1+1/s}\rho^2)^{-\tau}
\geq c_K4^{-(1+1/s)\tau}(1+\rho^2)^{-\tau}$, so the claim holds with
$c_3 \coloneqq 4^{(2+1/s)\tau}C_K/c_K$. Finally, if $\eta \geq c\|p-q\|_2^2$ for
this alternative, dividing both bounds by $\|p-q\|_2^2 = \Delta$ gives the bound
on $\lambda$. \qed

\subsection{Proof of Proposition~\ref{prop:signal-rkhs}}\label{app:signal-rkhs}

Let $\mathcal{I} : \mathcal{H} \to L^2(R)$, $f \mapsto [f - \mathbb{E}_R f]$, be
the centered inclusion operator, with adjoint
$\mathcal{I}^* : L^2(R) \to \mathcal{H}$,
$f \mapsto \int k(\cdot,x)f(x)\,dR(x) - \mu_R\,\mathbb{E}_R f$, so that
$\Sigma_R = \mathcal{I}^*\mathcal{I}$ and
$\mathcal{T} = \mathcal{I}\mathcal{I}^*$; in particular, the nonzero
eigenvalues of $\mathcal{T}$ and $\Sigma_R$ coincide. Since
$\mathbb{E}_R w = 0$, we have $\mu_P - \mu_Q = \mathcal{I}^*w$, and therefore
$\eta = \langle \mathcal{T}g_\lambda(\mathcal{T})w, w\rangle_{L^2(R)}$
\citep[Lemma A.8]{hagrass2024a_spectral_twosample}. Let
$(\lambda_i, \tilde\phi_i)_i$ be the eigenpairs of $\mathcal{T}$ with
$\lambda_i > 0$, and write $w = \mathcal{T}^\theta f$ with
$f = \mathcal{T}^{-\theta}w$. Since $w \in \Ran(\mathcal{T}^\theta)$ has no
component in $\ker\mathcal{T}$,
\[
  \eta = \sum_i \lambda_i g_\lambda(\lambda_i)\,
  \langle w, \tilde\phi_i\rangle_{L^2(R)}^2,
  \qquad
  \|w\|_{L^2(R)}^2 = \sum_i \langle w, \tilde\phi_i\rangle_{L^2(R)}^2 .
\]
Every term of the first sum is nonnegative, and for $\lambda_i \geq \lambda$,
\eqref{eq:g-lower} gives
$\lambda_i g_\lambda(\lambda_i) \geq C_\ell^2\lambda_i/(\lambda_i+\lambda)
\geq C_\ell^2/2$. Since
$\langle w, \tilde\phi_i\rangle_{L^2(R)} = \lambda_i^{\theta}
\langle f, \tilde\phi_i\rangle_{L^2(R)}$,
\[
  \sum_{\lambda_i < \lambda} \langle w, \tilde\phi_i\rangle_{L^2(R)}^2
  = \sum_{\lambda_i < \lambda} \lambda_i^{2\theta}
  \langle f, \tilde\phi_i\rangle_{L^2(R)}^2
  \;\leq\; \lambda^{2\theta}\,\|f\|_{L^2(R)}^2 .
\]
Therefore
\[
  \eta \;\geq\; \frac{C_\ell^2}{2}\sum_{\lambda_i \geq \lambda}
  \langle w, \tilde\phi_i\rangle_{L^2(R)}^2
  \;\geq\; \frac{C_\ell^2}{2}\Bigl(\|w\|_{L^2(R)}^2
  - \lambda^{2\theta}\bigl\|\mathcal{T}^{-\theta}w\bigr\|_{L^2(R)}^2\Bigr)
  \;\geq\; \frac{C_\ell^2}{4}\,\|w\|_{L^2(R)}^2 ,
\]
where the last step uses the assumed lower bound on $\|w\|^2_{L^2(R)}$. For
$R=(P+Q)/2$, $w = 2u$ and $\|w\|^2_{L^2(R)} = 4\rho^2(P,Q)$; for $R=P_0$,
$\|w\|^2_{L^2(R)} = \chi^2(P\,\|\,P_0)$. \qed

\subsection{Proof of Theorem~\ref{thm:separation-rkhs}}\label{app:separation-rkhs}

The argument is that of Theorem~\ref{thm:separation}, with
Proposition~\ref{prop:signal-rkhs} in place of Theorem~\ref{thm:signal}. Let
$S \coloneqq \sup\|\mathcal{T}^{-\theta}w\|_{L^2(R)}$ and
$\sigma \coloneqq \inf\|\Sigma_R\|_{\mathrm{op}}$, the supremum and the infimum
being over $\mathcal{P}_{\theta,\Delta}$, so that $S < \infty$ and
$\sigma > 0$ by assumption. If $\lambda \leq \sigma$ and the lower end
of~\eqref{eq:lambda-range} holds, then $P(\mathcal{E}) \geq 1-\delta/2$ by
Lemma~\ref{lem:B-general}, and Theorem~\ref{thm:power} applies. If moreover
$\Delta \geq 2S^2\lambda^{2\theta}$, Proposition~\ref{prop:signal-rkhs} gives
$\eta \geq C_\ell^2\Delta/4$, and since
$\mathcal{N}_2(\lambda) \vee 1 \leq (\kappa/\lambda)^{1/2}$ for
$\lambda \leq \kappa$, condition~\eqref{cond:eta} holds as soon as
$\Delta \gtrsim \lambda^{-1/2}/(\delta N)$. The two requirements on $\Delta$
balance at $\lambda \asymp N^{-2/(4\theta+1)}$, where both read
$\Delta \gtrsim N^{-4\theta/(4\theta+1)}$, and this $\lambda$ is at most
$\sigma$ for $N$ large. In part~(i), the covariance samples are drawn at no
cost to the test budget, so the lower end of~\eqref{eq:lambda-range} is met by
taking $l$ large enough. In part~(ii), $l \asymp N$ and the lower end
of~\eqref{eq:lambda-range} requires $\lambda \gtrsim \log N/N$. For
$\theta > 1/4$ this holds at $\lambda \asymp N^{-2/(4\theta+1)}$ for $N$ large,
which gives~(a). For $\theta \leq 1/4$, take $\lambda \asymp \log N/N$: the two
requirements become $\Delta \gtrsim (\log N/N)^{2\theta}$ and
$\Delta \gtrsim (N\log N)^{-1/2}$, and the first implies the second, which
gives~(b). \qed

\subsection{Proof of Theorem~\ref{thm:separation}}

The result follows from Theorem~\ref{thm:power}, whose condition we verify
below, together with the regularity event $\mathcal{E}$ supplied by
Lemma~\ref{lem:B-general}.

Since $(P,Q) \in \mathcal{F}_{s,\Delta}$, we have $\|\psi\|_2^2 \geq \Delta$.
We apply Theorem~\ref{thm:signal} with $t=h$. By~\eqref{eq:mK-bandwidth},
the two conditions in~\eqref{eq:resolution} read
$\lambda\le\frac{Mc_\varphi}{4}h^d$ and $\|\psi\|_2^2\ge4L^2h^{2s}$; taking
$\lambda \asymp h^d$ and $h \asymp \Delta^{1/(2s)}$ with sufficiently small
implicit constants, both hold. The lower bound~\eqref{eq:g-lower} required
there is implied by~\eqref{eq:g-two-sided}. Hence
$\eta \geq \frac{C_\ell^2}{2M}\|\psi\|_2^2 \geq \frac{C_\ell^2}{2M}\Delta$.
Moreover, $\mathcal{N}_2^2(\lambda)\le\tr(\Sigma_R)/\lambda\le\kappa/\lambda$
with $\kappa=\varphi(0)$ independent of $h$, so
$\mathcal{N}_2(\lambda) \vee 1 \lesssim \lambda^{-1/2}$ for $\lambda\le\kappa$,
and condition~\eqref{cond:eta} of Theorem~\ref{thm:power} is implied by
\[
  \Delta \;\gtrsim\; \frac{1}{\delta N \sqrt{\lambda}}
  \;=\; \frac{1}{\delta N\, h^{d/2}} ,
\]
which, after substituting $h \asymp \Delta^{1/(2s)}$ and absorbing $\delta$
into the implicit constant, reads $\Delta^{1+d/(4s)} \gtrsim N^{-1}$, i.e.
\begin{equation}\label{eq:sob-main}
  \Delta \;\gtrsim\; N^{-\frac{4s}{4s+d}} ,
\end{equation}
the rate of parts~(i) and~(ii)(a).

It remains to establish $\mathcal{E}$, for which the size of
$\|\Sigma_R\|_{\mathrm{op}}$ relative to $h^d$ is what matters. Since
$\|r\|_\infty \leq M$ and
$\widehat K(\omega)=h^d\widehat\varphi(h\omega)\le h^d\|\widehat\varphi\|_\infty$,
\eqref{eq:native} gives, for any $f \in \mathcal{H}$,
\[
  \langle \Sigma_R f, f\rangle_{\mathcal{H}}
  \;\leq\; \|f\|_{L^2(R)}^2 \;\leq\; M\|f\|_{L^2}^2
  \;\leq\; M\|\widehat\varphi\|_\infty\,h^d\,\|f\|_{\mathcal{H}}^2 ,
\]
so that $\|\Sigma_R\|_{\mathrm{op}} \lesssim h^d$ always. Together with the
assumed lower bound, this gives $\|\Sigma_R\|_{\mathrm{op}} \asymp h^d$, so the
choice $\lambda \asymp h^d$ places $\lambda$ at the order of the top of the
spectrum of $\Sigma_R$.

By Lemma~\ref{lem:B-general}, $\mathcal{E}$ holds with probability at least
$1-\delta/2$ once~\eqref{eq:lambda-range} holds, that is, once
$l^{-1}\log(l/\delta) \lesssim \lambda \leq \|\Sigma_R\|_{\mathrm{op}}$. The
upper end is met by taking the implicit constant in $\lambda \asymp h^d$
sufficiently small, by the assumed lower bound on
$\|\Sigma_R\|_{\mathrm{op}}$. The lower end reads
$l \gtrsim h^{-d}\log(h^{-d})$; in case~(i) this is met by drawing additional
covariance samples, while in case~(ii), $l \asymp N$ and the requirement
becomes $h^d \gtrsim \log N/N$, i.e.
\begin{equation}\label{eq:sob-split}
  \Delta \;\gtrsim\; \Bigl(\frac{\log N}{N}\Bigr)^{2s/d} ,
\end{equation}
the rate of part~(ii)(b). When $s > d/4$, the choice
$h \asymp N^{-2/(4s+d)}$ of part~(ii)(a) satisfies $h^d \gtrsim \log N/N$, so
that~\eqref{eq:sob-main} is the binding condition.

The stated choices of $h$ and $\lambda$ in each case are precisely those for
which $h \asymp \Delta^{1/(2s)}$ at the corresponding $\Delta$, and the
conditions of Theorem~\ref{thm:power} are thus satisfied uniformly over
$(P,Q) \in \mathcal{F}_{s,\Delta}$. \qed

\subsection{Proof of Proposition~\ref{prop:mmd}}\label{app:mmd}

With $g_\lambda \equiv \lambda^{-1}$ on $[0,2\kappa]$, \eqref{Eq:glambda} gives
$g_\lambda^{1/2}(\hat\Sigma_R) = \lambda^{-1/2}I$, so
$\langle A(x), A(y)\rangle_{\mathcal{H}} = \lambda^{-1}k(x,y)$ and
$\hat\eta = \lambda^{-1}\widehat{D}^2_{\mathrm{MMD}}$. The same identity holds
for every permutation of the pooled sample, so $\hat q_{1-\alpha}$ is
$\lambda^{-1}$ times the $(1-\alpha)$-quantile of the permutation distribution of
$\widehat{D}^2_{\mathrm{MMD}}$, and the two tests reject on the same event,
whatever $\lambda$ and $R$ are. Replacing $k$ by $ck$ with $c > 0$ multiplies
$\widehat{D}^2_{\mathrm{MMD}}$ and its permutation quantile by $c$, which leaves
the MMD test unchanged.

For the power statement, since the test depends on neither $\lambda$ nor $R$,
both may be chosen for the analysis. Take $R = (P+Q)/2$, for which
$\Sigma_P, \Sigma_Q \preccurlyeq 2\Sigma_R$ and $\|r\|_\infty \leq M$ on
$\mathcal{F}_{s,\Delta}$ (Section~\ref{sec:reference-measure}), and
$\lambda = \frac{Mc_\varphi}{4}h^d$. Since
$g_\lambda(x)(x+\lambda) = 1 + x/\lambda \geq 1$, the lower
bound~\eqref{eq:g-lower} holds with $C_\ell = 1$, and Theorem~\ref{thm:signal}
with $t = h$ gives $\eta \geq \|\psi\|_2^2/(2M)$ whenever
$\|\psi\|_2^2 \geq 4L^2h^{2s}$. The proof of Theorem~\ref{thm:separation} shows
that $\|\Sigma_R\|_{\mathrm{op}} \leq M\|\widehat\varphi\|_\infty h^d$, so
$\lambda \geq b\,\|\Sigma_R\|_{\mathrm{op}}$ with
$b = c_\varphi/(4\|\widehat\varphi\|_\infty)$, and Lemma~\ref{lem:B-flat} shows
that $\mathcal{E}(\sqrt{1+b^{-1}},1)$ holds surely, with no covariance samples.
Theorem~\ref{thm:power} therefore applies, and since
$\mathcal{N}_2(\lambda) \vee 1 \lesssim \lambda^{-1/2}$, the computation leading
to~\eqref{eq:sob-main} in the proof of Theorem~\ref{thm:separation} gives power
at least $1-2\delta$ for $\Delta \gtrsim N^{-4s/(4s+d)}$, with
$h \asymp N^{-2/(4s+d)}$ and for every $s > 0$. \qed

\subsection{Proof of Lemma~\ref{lem:scaling}}

Since $\tilde k = h^{-\gamma}k$ with $h^{-\gamma} > 0$, the two kernels have
the same set of functions as their RKHS, with norms differing by the
reciprocal factor. Indeed, on
$\mathcal{H}_0 \coloneqq \mathrm{span}\{k(\cdot,x) : x \in \mathcal{X}\}$,
the reproducing property for $\tilde k$ gives
\[
  \bigl\langle \tilde k(\cdot,x), \tilde k(\cdot,y)\bigr\rangle_{\tilde{\mathcal H}}
  = \tilde k(x,y) = h^{-\gamma}k(x,y)
  = h^{-\gamma}\bigl\langle k(\cdot,x), k(\cdot,y)\bigr\rangle_{\mathcal H},
\]
and since $\tilde k(\cdot,x) = h^{-\gamma}k(\cdot,x)$ as functions, the
left-hand side equals
$h^{-2\gamma}\langle k(\cdot,x), k(\cdot,y)\rangle_{\tilde{\mathcal H}}$,
whence $\langle f,g\rangle_{\tilde{\mathcal H}} = h^{\gamma}
\langle f,g\rangle_{\mathcal H}$ for $f, g \in \mathcal{H}_0$. As the two
norms are proportional, they induce the same Cauchy sequences and the same
completion, so $\tilde{\mathcal{H}} = \mathcal{H}$ as sets and the identity
extends to all $f, g \in \mathcal{H}$.

For the tensor product, recall $(f \otimes_{\mathcal H} g)z =
\langle z, g\rangle_{\mathcal H}\,f$ for $z \in \mathcal{H}$. Applying the
inner product identity,
\[
  (f \otimes_{\tilde{\mathcal H}} g)z
  = \langle z, g\rangle_{\tilde{\mathcal H}}\,f
  = h^\gamma\langle z, g\rangle_{\mathcal H}\,f
  = h^\gamma\,(f \otimes_{\mathcal H} g)z .
\]

For the covariance operator, $\tilde k(\cdot,x) = h^{-\gamma}k(\cdot,x)$ gives
$\tilde\mu_R = h^{-\gamma}\mu_R$, and hence
$\tilde k(\cdot,x) - \tilde\mu_R = h^{-\gamma}\bigl(k(\cdot,x)-\mu_R\bigr)$.
Substituting into the definition of $\tilde\Sigma_R$ and using bilinearity of
the tensor product,
\[
  \tilde\Sigma_R
  = h^{-2\gamma}\int_{\mathcal X}\bigl(k(\cdot,x) - \mu_R\bigr)
    \otimes_{\tilde{\mathcal H}}
    \bigl(k(\cdot,x) - \mu_R\bigr)\,dR(x),
\]
and converting the tensor product to $\mathcal{H}$ by the identity just
established contributes a further factor $h^{\gamma}$, giving
$\tilde\Sigma_R = h^{-\gamma}\Sigma_R$. The same argument applied to the
U-statistic estimator gives $\hat{\tilde\Sigma}_R = h^{-\gamma}\hat\Sigma_R$.
Finally, for any $\tilde\lambda > 0$,
\[
  \tilde\Sigma_{R,\tilde\lambda}
  = h^{-\gamma}\Sigma_R + \tilde\lambda I
  = h^{-\gamma}\bigl(\Sigma_R + h^\gamma\tilde\lambda I\bigr)
  = h^{-\gamma}\Sigma_{R, h^\gamma\tilde\lambda} . \qed
\]

\subsection{Proof of Proposition~\ref{prop:normalization}}

By Lemma~\ref{lem:scaling}, $\hat{\tilde\Sigma}_R = h^{-\gamma}\hat\Sigma_R$ as
operators on the common set $\mathcal{H} = \tilde{\mathcal H}$. The two inner
products are proportional, so the two operators have the same eigenvectors, the
eigenvalues of $\hat{\tilde\Sigma}_R$ are those of $\hat\Sigma_R$ multiplied by
$h^{-\gamma}$, and the orthogonal projections in~\eqref{Eq:glambda} are the same
in both spaces. Since
$g_{\tilde\lambda}(h^{-\gamma}x)
= \tilde\lambda^{-1}G\bigl(x/(h^\gamma\tilde\lambda)\bigr) = h^\gamma g_\lambda(x)$
with $\lambda = h^\gamma\tilde\lambda$, it follows that
$g_{\tilde\lambda}(\hat{\tilde\Sigma}_R) = h^\gamma g_\lambda(\hat\Sigma_R)$.
Hence
$\tilde A(x) = g_{\tilde\lambda}^{1/2}(\hat{\tilde\Sigma}_R)\,\tilde k(\cdot,x)
= h^{\gamma/2}g_\lambda^{1/2}(\hat\Sigma_R)\,h^{-\gamma}k(\cdot,x)
= h^{-\gamma/2}A(x)$, and by Lemma~\ref{lem:scaling},
\[
  \langle \tilde A(x), \tilde A(y)\rangle_{\tilde{\mathcal H}}
  = h^{\gamma}\,\langle \tilde A(x), \tilde A(y)\rangle_{\mathcal H}
  = \langle A(x), A(y)\rangle_{\mathcal H}
\]
for all $x, y$. The two statistics~\eqref{eq:test-stat} are therefore equal, for
the observed sample and for every permutation of it, and so are their
permutation quantiles. \qed

\subsection{Normalization for a general regularizer}\label{app:normalization}

For a general $g_\lambda$, the reparametrization $\lambda \coloneqq
h^\gamma\tilde\lambda$ still matches the two tests, in the following sense.
Define the test statistic and population discrepancy under $\tilde k$ with
regularization $\tilde\lambda$ exactly as in~\eqref{eq:test-stat}
and~\eqref{eq:eta-pop}, but with $k, \lambda$ replaced by
$\tilde k, \tilde\lambda$, and write $\hat{\tilde\eta}_{\tilde\lambda}$ and
\[
  \tilde\eta_{\tilde\lambda}
  \coloneqq \bigl\|g_{\tilde\lambda}^{1/2}(\tilde\Sigma_R)
  (\tilde\mu_P - \tilde\mu_Q)\bigr\|_{\tilde{\mathcal H}}^2
\]
for the resulting quantities.

\begin{proposition}[Reduction to the bounded kernel]
\label{prop:equivalence}
Let $\tilde\lambda > 0$ and define $\lambda \coloneqq h^\gamma\tilde\lambda$,
\[
  \tilde\B \coloneqq g_{\tilde\lambda}^{1/2}\bigl(\hat{\tilde\Sigma}_R\bigr)\,
  \tilde\Sigma_{R,\tilde\lambda}^{1/2},
  \qquad
  \tilde\B_0 \coloneqq g_{\tilde\lambda}^{1/2}\bigl(\tilde\Sigma_R\bigr)\,
  \tilde\Sigma_{R,\tilde\lambda}^{1/2},
\]
and let $\tilde A(x) \coloneqq g_{\tilde\lambda}^{1/2}
(\hat{\tilde\Sigma}_R)\,\tilde k(\cdot, x)$ be the analogue of $A(x)$ under
$\tilde k$. Then for any $x, y \in \mathcal{X}$,
\begin{equation}\label{eq:inner-equiv}
  \langle \tilde A(x), \tilde A(y)\rangle_{\tilde{\mathcal H}}
  = \bigl\langle
  \tilde\B\,\Sigma_{R,\lambda}^{-1/2}\,k(\cdot, x),\;
  \tilde\B\,\Sigma_{R,\lambda}^{-1/2}\,k(\cdot, y)\bigr\rangle_{\mathcal H}.
\end{equation}
Consequently,
\[
  \E\bigl(\hat{\tilde\eta}_{\tilde\lambda}\,\big|\, \mathbb{Z}_l\bigr)
  = \bigl\|\tilde\B\,\Sigma_{R,\lambda}^{-1/2}(\mu_P - \mu_Q)
  \bigr\|_{\mathcal H}^2,
  \qquad
  \tilde\eta_{\tilde\lambda}
  = \bigl\|\tilde\B_0\,\Sigma_{R,\lambda}^{-1/2}(\mu_P - \mu_Q)
  \bigr\|_{\mathcal H}^2,
\]
and the conditional variance
$\mathrm{Var}(\hat{\tilde\eta}_{\tilde\lambda}\,\big|\, \mathbb{Z}_l)$
admits the same bound as in Proposition~\ref{prop:var-bound}, with $\B$
replaced by $\tilde\B$.
\end{proposition}
The variance bound of Proposition~\ref{prop:var-bound} carries over verbatim
with $\tilde{\B}$ in place of $\B$, since its proof leaves the operator norm
of $\B$ in the conclusion and uses only the dominance conditions. The
remaining results transfer because the regularity event $\mathcal{E}$ does:
Lemma~\ref{lem:B-general} applies to $g_{\tilde\lambda}$ on the spectrum of
$\tilde\Sigma_R$ exactly as it does to $g_\lambda$ on that of $\Sigma_R$, and
the upper constraint on the regularization parameter there is invariant under
the rescaling, since $\tilde\lambda \leq \|\tilde\Sigma_R\|_{\mathrm{op}}$ and
$\lambda \leq \|\Sigma_R\|_{\mathrm{op}}$ are the same inequality multiplied
by $h^{-\gamma}$. Consequently the variance bound, threshold bound, power
criterion, and separation rates of Sections~\ref{sec:power}
and~\ref{sec:minimax} all apply to $\hat{\tilde\eta}_{\tilde\lambda}$, with
$\lambda$ replaced by $h^\gamma\tilde\lambda$ throughout. In particular, the
condition $\lambda \lesssim h^d$ of Theorem~\ref{thm:signal} with $t = h$ becomes
$\tilde\lambda \lesssim h^{d-\gamma}$, and the optimal Sobolev regularization
$\lambda \asymp N^{-2d/(4s+d)}$ becomes
$\tilde\lambda \asymp N^{-2(d-\gamma)/(4s+d)}$.

\begin{proof}
We first establish~\eqref{eq:inner-equiv}. Since $\tilde k(\cdot,x)
= h^{-\gamma}k(\cdot,x)$,
\[
  \tilde A(x)
  = g_{\tilde\lambda}^{1/2}(\hat{\tilde\Sigma}_R)\,\tilde k(\cdot,x)
  = h^{-\gamma}\,g_{\tilde\lambda}^{1/2}(\hat{\tilde\Sigma}_R)\,k(\cdot,x).
\]
By definition of $\tilde{\B}$, $g_{\tilde\lambda}^{1/2}(\hat{\tilde\Sigma}_R)
= \tilde{\B}\,\tilde\Sigma_{R,\tilde\lambda}^{-1/2}$, so
\[
  \tilde A(x) = h^{-\gamma}\, \tilde{\B}\,\tilde\Sigma_{R,\tilde\lambda}^{-1/2}\,
  k(\cdot,x).
\]
Applying Lemma~\ref{lem:scaling}, $\tilde\Sigma_{R,\tilde\lambda}^{-1/2}
= h^{\gamma/2}\Sigma_{R,\lambda}^{-1/2}$ (with $\lambda = h^\gamma
\tilde\lambda$), giving
\begin{equation}\label{eq:tildeA-id}
  \tilde A(x) = h^{-\gamma/2}\, \tilde{\B}\,\Sigma_{R,\lambda}^{-1/2}\,k(\cdot,x).
\end{equation}
Both $\tilde{\B}\,\Sigma_{R,\lambda}^{-1/2}\,k(\cdot,x)$ and $\tilde A(x)$ are
elements of $\mathcal{H} = \tilde{\mathcal H}$ (viewed as a set). Applying
Lemma~\ref{lem:scaling} to convert the inner product on $\tilde{\mathcal H}$
to one on $\mathcal H$, and then~\eqref{eq:tildeA-id},
\begin{align*}
  \langle \tilde A(x), \tilde A(y)\rangle_{\tilde{\mathcal H}}
  &= h^\gamma\,\langle \tilde A(x), \tilde A(y)\rangle_{\mathcal H} \\
  &= h^\gamma\cdot h^{-\gamma}\bigl\langle
     \tilde{\B}\,\Sigma_{R,\lambda}^{-1/2}k(\cdot,x),\;
     \tilde{\B}\,\Sigma_{R,\lambda}^{-1/2}k(\cdot,y)\bigr\rangle_{\mathcal H} \\
  &= \bigl\langle \tilde{\B}\,\Sigma_{R,\lambda}^{-1/2}k(\cdot,x),\;
                  \tilde{\B}\,\Sigma_{R,\lambda}^{-1/2}k(\cdot,y)
     \bigr\rangle_{\mathcal H},
\end{align*}
which is~\eqref{eq:inner-equiv}.

For the conditional mean, $\E(\hat{\tilde\eta}_{\tilde\lambda}\mid
\mathbb{Z}_l) = \|g_{\tilde\lambda}^{1/2}(\hat{\tilde\Sigma}_R)
(\tilde\mu_P - \tilde\mu_Q)\|_{\tilde{\mathcal H}}^2$ since
$\hat{\tilde\eta}_{\tilde\lambda}$ is, conditionally on $\mathbb{Z}_l$, a
two-sample U-statistic of the inner product kernel
$\langle\tilde A(\cdot), \tilde A(\cdot)\rangle_{\tilde{\mathcal H}}$. The
same calculation as above with $\mu_P - \mu_Q$ in place of $k(\cdot,x)$ gives
\[
  g_{\tilde\lambda}^{1/2}(\hat{\tilde\Sigma}_R)(\tilde\mu_P - \tilde\mu_Q)
  = h^{-\gamma/2}\,\tilde{\B}\,\Sigma_{R,\lambda}^{-1/2}(\mu_P - \mu_Q),
\]
and converting the norm,
\[
  \E\bigl(\hat{\tilde\eta}_{\tilde\lambda}\,\big|\,\mathbb{Z}_l\bigr)
  = h^\gamma\cdot h^{-\gamma}\bigl\|\tilde{\B}\,\Sigma_{R,\lambda}^{-1/2}
  (\mu_P - \mu_Q)\bigr\|_{\mathcal H}^2
  = \bigl\|\tilde{\B}\,\Sigma_{R,\lambda}^{-1/2}(\mu_P - \mu_Q)
  \bigr\|_{\mathcal H}^2 .
\]
The corresponding identity for $\tilde\eta_{\tilde\lambda}$ follows by the
same computation with $\tilde\Sigma_R$ in place of $\hat{\tilde\Sigma}_R$
throughout, giving $\tilde{\B}_0$ in place of $\tilde{\B}$.

For the variance, observe that $\hat{\tilde\eta}_{\tilde\lambda}$ written
term-by-term using~\eqref{eq:inner-equiv} is identical in form to
$\hat\eta$ with $\B$ replaced everywhere by $\tilde{\B}$. The proof of
Proposition~\ref{prop:var-bound} uses only the dominance conditions
$\Sigma_P \preccurlyeq c\,\Sigma_R$ and $\Sigma_Q \preccurlyeq c\,\Sigma_R$,
which are unchanged, and leaves the operator norm of $\B$ in the bound; the
same bound therefore holds with $\|\tilde{\B}\|_{\mathrm{op}}$ in its place.
\end{proof}

\subsection{Proof of Theorem~\ref{thm:sobolev-fixed}}

The result follows from Theorem~\ref{thm:power}, whose condition we verify
below, together with the regularity event $\mathcal{E}$ supplied by
Lemma~\ref{lem:B-general}.

Since $(P,Q) \in \mathcal{F}_{s,\Delta}$, we have $\|\psi\|_2^2 \geq \Delta$.
For $\lambda \le M c_K 2^{-\tau-2}$, let
$t := \bigl(2^{\tau+2}\lambda/(Mc_K)\bigr)^{1/(2\tau)} \le 1$.
By~\eqref{eq:Khat-lower} and $1+t^{-2}\le 2t^{-2}$,
\[
  \frac{M}{4}\,m_K(1/t) \;\geq\; \frac{M c_K}{4}\,(1+t^{-2})^{-\tau}
  \;\geq\; \frac{M c_K}{4}\,2^{-\tau}\,t^{2\tau} \;=\; \lambda ,
\]
so the first condition in~\eqref{eq:resolution} holds, and the second,
$\|\psi\|_2^2 \geq 4L^2t^{2s}$, holds whenever $\Delta \geq C_1\lambda^{s/\tau}$
with $C_1 := 4L^2\bigl(2^{\tau+2}/(Mc_K)\bigr)^{s/\tau}$. The lower
bound~\eqref{eq:g-lower} is implied by~\eqref{eq:g-two-sided}, so
Theorem~\ref{thm:signal} gives
$\eta \geq \frac{C_\ell^2}{2M}\|\psi\|_2^2 \geq \frac{C_\ell^2}{2M}\Delta$.
By~\eqref{eq:eigdecay} and \citet[Lemma~B.9]{sriperumbudur2022approx},
$\mathcal{N}_2(\lambda) \vee 1 \lesssim \lambda^{-d/(4\tau)}$ for
$\lambda \le 1$. Condition~\eqref{cond:eta} of Theorem~\ref{thm:power} is
therefore implied by
\begin{equation}\label{eq:fixed-cond}
  \Delta \;\gtrsim\; \max\Bigl\{\lambda^{s/\tau},\;
  \frac{1}{\delta N \lambda^{d/(4\tau)}}\Bigr\}.
\end{equation}
For $\lambda \asymp N^{-4\tau/(4s+d)}$ both terms are of order
$N^{-4s/(4s+d)}$, after absorbing $\delta$ into the implicit constant.

It remains to establish $\mathcal{E}$ through Lemma~\ref{lem:B-general}, which
requires $\frac{140\kappa}{l}\log\frac{128\kappa l}{\delta} \le \lambda \le
\|\Sigma_R\|_{\mathrm{op}}$. Since $\lambda \to 0$ and
$\inf\|\Sigma_R\|_{\mathrm{op}} > 0$, the upper end holds for $N$ large enough.
In case~(i) the lower end is met by drawing
$l \gtrsim \lambda^{-1}\log(\lambda^{-1}\delta^{-1})$ covariance samples. In
case~(ii), $l \asymp N$ and the lower end requires $\lambda \gtrsim \log N/N$.
For $\lambda \asymp N^{-4\tau/(4s+d)}$ this holds if and only if
$4\tau < 4s+d$, i.e.\ $s > \tau - d/4$, which is case~(a). Otherwise take
$\lambda \asymp \log N/N$: then $\lambda^{s/\tau} \asymp (\log N/N)^{s/\tau}$,
and $\lambda^{-d/(4\tau)}N^{-1} \lesssim (\log N/N)^{s/\tau}$ because
$(s+d/4)/\tau \le 1$, so~\eqref{eq:fixed-cond} holds for
$\Delta \gtrsim (\log N/N)^{s/\tau}$, which is case~(b). \qed

\subsection{Proof of Proposition~\ref{prop:asymmetric}}\label{app:ref-gof}

\emph{Dominance and covariance budget.} The law of total variance applied to
$R = \varepsilon P + (1-\varepsilon)P_0$ gives
\begin{equation}\label{eq:asym-cov}
  \Sigma_R = \varepsilon\Sigma_P + (1-\varepsilon)\Sigma_{P_0}
  + \varepsilon(1-\varepsilon)\,(\mu_P-\mu_{P_0})\otimes_{\mathcal{H}}(\mu_P-\mu_{P_0}),
\end{equation}
so $\Sigma_R \succcurlyeq \varepsilon\Sigma_P$ and
$\Sigma_R \succcurlyeq (1-\varepsilon)\Sigma_{P_0}$, and the dominance condition
holds with $c = \max\{\varepsilon^{-1}, (1-\varepsilon)^{-1}\} = \varepsilon^{-1}$
for $\varepsilon \in (0,1/2]$. Since this constant varies with $N$, we track it
through Proposition~\ref{prop:var-bound}, whose bound reads
\[
  \mathrm{Var}(\hat\eta \mid \mathbb{Z}_l) \lesssim
  \frac{c^2\|\B\|^4_{\mathrm{op}}\mathcal{N}_2^2(\lambda)}{N^2}
  + \frac{c\,\|\B\|^2_{\mathrm{op}}\zeta}{N} .
\]
Repeating the argument of the proof of Theorem~\ref{thm:power} with the
constant retained replaces~\eqref{cond:eta} by
\begin{equation}\label{eq:gof-cond1}
  \eta \gtrsim \frac{\mathcal{N}_2(\lambda) \vee 1}{\delta\,\varepsilon N} .
\end{equation}
The sample budget enters through Lemma~\ref{lem:B-general}, which requires
$l^{-1}\log(l/\delta) \lesssim \lambda \leq \|\Sigma_R\|_{\mathrm{op}}$ for the
regularity event $\mathcal{E}$. The constraint $\varepsilon l \lesssim n \asymp N$
permits $l \asymp N/\varepsilon$, so, for the choices $\varepsilon \geq N^{-1}$
made below, the lower end becomes
\begin{equation}\label{eq:gof-cond2}
  \lambda \gtrsim \frac{\varepsilon\log N}{N} .
\end{equation}
The upper end requires no further assumption: by~\eqref{eq:asym-cov},
$\Sigma_R \succcurlyeq \tfrac12\Sigma_{P_0}$ with $P_0$ fixed, so
$\|\Sigma_R\|_{\mathrm{op}} \geq \tfrac12\|\Sigma_{P_0}\|_{\mathrm{op}}$. For a
fixed kernel, as in cases~(b) and~(c), this is a fixed positive constant while
$\lambda \to 0$; for a bandwidth family, as in case~(a), it is of order $h^d$ by
Remark~\ref{rem:op-lower} applied to the fixed measure $P_0$, while
$\lambda \asymp h^d$ with a small implicit constant.

\emph{Reduction to two exponents.} In each case, the lower bound on $\eta$
(Theorem~\ref{thm:signal} or Proposition~\ref{prop:signal-rkhs}) and the bound on
the effective dimension take the form
\begin{equation}\label{eq:gof-exponents}
  \eta \gtrsim \Delta \ \text{ whenever } \ \Delta \gtrsim \lambda^{a_1},
  \qquad
  \mathcal{N}_2(\lambda) \vee 1 \lesssim \lambda^{-a_2},
\end{equation}
for exponents $a_1, a_2 > 0$ given as follows.
\begin{enumerate}[label=\textnormal{(\alph*)}]
  \item Bandwidth family: $a_1 = 2s/d$ and $a_2 = 1/2$. With $t = h$ and
  $\lambda \asymp h^d$, Theorem~\ref{thm:signal} applies as in the proof of
  Theorem~\ref{thm:separation} whenever
  $\Delta \geq 4L^2h^{2s} \asymp \lambda^{2s/d}$, and
  $\mathcal{N}_2(\lambda) \leq (\kappa/\lambda)^{1/2}$ with
  $\kappa = \varphi(0)$.
  \item Fixed kernel: $a_1 = s/\tau$ and $a_2 = d/(4\tau)$. The lower bound is
  that of the proof of Theorem~\ref{thm:sobolev-fixed}, using
  $\|r\|_\infty \leq \max\{\|p\|_\infty, \|p_0\|_\infty\} \leq M$. For the
  effective dimension, \eqref{eq:asym-cov} and Weyl's inequality give
  \[
    \lambda_{2i}(\Sigma_R)
    \leq \lambda_{2i-1}\bigl(\varepsilon\Sigma_P + (1-\varepsilon)\Sigma_{P_0}\bigr)
    \leq \varepsilon\lambda_i(\Sigma_P) + (1-\varepsilon)\lambda_i(\Sigma_{P_0})
    \leq A\,i^{-2\tau/d},
  \]
  where the first inequality holds because adding a rank-one operator shifts
  the eigenvalue indices by at most one. Hence
  $\lambda_k(\Sigma_R) \leq \max\{\kappa, 3^{2\tau/d}A\}\,k^{-2\tau/d}$ for all
  $k \geq 1$, uniformly in $\varepsilon$, and
  \citet[Lemma~B.9]{sriperumbudur2022approx} gives
  $\mathcal{N}_2(\lambda) \vee 1 \lesssim \lambda^{-d/(4\tau)}$.
  \item General-domain class: $a_1 = 2\theta$ and $a_2 = 1/2$, by
  Proposition~\ref{prop:signal-rkhs} together with the uniform bound on
  $\|\mathcal{T}^{-\theta}w\|_{L^2(R)}$, and the generic bound
  $\mathcal{N}_2(\lambda) \lesssim \lambda^{-1/2}$.
\end{enumerate}
By~\eqref{eq:gof-cond1} and~\eqref{eq:gof-exponents}, the power is at least
$1-2\delta$ as soon as
\begin{equation}\label{eq:gof-rates}
  \Delta \gtrsim \max\Bigl\{\lambda^{a_1},\;
  \frac{\lambda^{-a_2}}{\delta\,\varepsilon N}\Bigr\}
\end{equation}
for some $\lambda$ satisfying~\eqref{eq:gof-cond2}.

\emph{Choice of $\lambda$ and $\varepsilon$.} Let $a := a_1 + a_2$ and absorb
$\delta$ into the implicit constants. For fixed $\varepsilon$, the right-hand
side of~\eqref{eq:gof-rates} is minimized at
$\lambda_\varepsilon := (\varepsilon N)^{-1/a}$, where both terms equal
$(\varepsilon N)^{-a_1/a}$. If $a > 1$, take $\varepsilon = 1/2$: then
$\lambda_{1/2} \asymp N^{-1/a}$ satisfies~\eqref{eq:gof-cond2}, and
$\Delta \gtrsim N^{-a_1/a}$ suffices. If $a \leq 1$, the best admissible
$\lambda$ is $\max\{\lambda_\varepsilon, \varepsilon\log N/N\}$. When
$\lambda_\varepsilon \geq \varepsilon\log N/N$ this gives
$\Delta \asymp (\varepsilon N)^{-a_1/a}$, which is decreasing in
$\varepsilon$; otherwise the first term of~\eqref{eq:gof-rates} dominates at
$\lambda = \varepsilon\log N/N$ and gives
$\Delta \asymp (\varepsilon\log N/N)^{a_1}$, which is increasing in
$\varepsilon$. The optimal $\varepsilon$ therefore equates
$\lambda_\varepsilon = \varepsilon\log N/N$, which gives
\[
  \lambda \asymp \Bigl(\frac{\log N}{N^2}\Bigr)^{1/(1+a)},
  \quad
  \varepsilon \asymp \frac{\lambda^{-a}}{N}
  \asymp N^{(a-1)/(a+1)}(\log N)^{-a/(a+1)},
  \quad
  \Delta \asymp \lambda^{a_1}
  = \Bigl(\frac{\sqrt{\log N}}{N}\Bigr)^{2a_1/(1+a)} .
\]
Since $a > 0$, this $\varepsilon$ is at most $1/2$ and at least $N^{-1}$ for
$N$ large enough.

\emph{The three cases.} In case~(a), $a = (4s+d)/(2d)$, so $a > 1$ if and only
if $s > d/4$, with $N^{-a_1/a} = N^{-4s/(4s+d)}$ and
$\lambda \asymp N^{-2d/(4s+d)}$; otherwise $(a-1)/(a+1) = (4s-d)/(4s+3d)$ and
$2a_1/(1+a) = 8s/(4s+3d)$. In case~(b), $a = (4s+d)/(4\tau)$, so $a > 1$ if
and only if $s > \tau - d/4$, with $N^{-a_1/a} = N^{-4s/(4s+d)}$ and
$\lambda \asymp N^{-4\tau/(4s+d)}$; otherwise
$(a-1)/(a+1) = (4s+d-4\tau)/(4s+d+4\tau)$ and
$2a_1/(1+a) = 8s/(4s+d+4\tau)$. In case~(c), $a = (4\theta+1)/2$, so
$a > 1$ if and only if $\theta > 1/4$, with
$N^{-a_1/a} = N^{-4\theta/(4\theta+1)}$ and
$\lambda \asymp N^{-2/(4\theta+1)}$; otherwise
$(a-1)/(a+1) = (4\theta-1)/(4\theta+3)$ and
$2a_1/(1+a) = 8\theta/(4\theta+3)$. These are the choices and
rates stated in the proposition. \qed

\subsection{Proof of Theorem~\ref{thm:adaptation}}\label{app:adaptation}

\textbf{Level.} Fix $(\lambda,h) \in \Lambda \times \mathcal{W}$ and write
$\alpha' \coloneqq \alpha/|\Lambda \times \mathcal{W}|$. Conditionally on
$\mathbb{Z}_l$, which is independent of the pooled sample, the pooled sample
is exchangeable under $H_0$, so the permutation construction gives
$P_{H_0}\{\hat\eta_{\lambda,h} > \hat q^{\,\lambda,h}_{1-\alpha'}\}
\leq \alpha'$. A union bound over the $|\Lambda \times \mathcal{W}|$ pairs
yields
\[
  P_{H_0}\Big\{\bigcup_{(\lambda,h)\in\Lambda\times\mathcal{W}}
  \hat\eta_{\lambda,h} > \hat q^{\,\lambda,h}_{1-\alpha'}\Big\}
  \;\leq\; |\Lambda \times \mathcal{W}|\,\alpha' \;=\; \alpha .
\]

\textbf{Step 1: the power criterion at the corrected level.} The union test
rejects whenever a single member of the grid rejects, so it suffices to
exhibit one pair $(\lambda,h) \in \Lambda \times \mathcal{W}$ for which the
individual test at level $\alpha'$ has power at least $1-2\delta$. Fix such a
pair and write $\hat\eta$, $\hat q_{1-\alpha'}$, $\B$ and $\zeta$ for the
corresponding quantities. We do not apply Theorem~\ref{thm:power} as stated,
since its constant $c_2$ grows like $\log^2(1/\alpha')$
through~\eqref{eq:dom2}; using the threshold bound~\eqref{eq:I2-bound}
directly keeps the dependence on $\log(1/\alpha')$ linear.

By~\eqref{eq:qs-bound} at level $\alpha'$, \eqref{eq:I2-bound} and Markov's
inequality, with probability at least $1-\delta/2$,
\begin{equation}\label{eq:adapt-thresh}
  \hat q_{1-\alpha'} \;\leq\; C\log(1/\alpha')\,
    \frac{\|\B\|_{\mathrm{op}}^2\,\mathcal{N}_2(\lambda)
  + \|\B\|_{\mathrm{op}}\sqrt{\zeta} + \zeta}{\sqrt{\delta}\,N},
\end{equation}
and by Proposition~\ref{prop:var-bound} and Chebyshev's inequality, with
probability at least $1-\delta$,
\[
  \hat\eta \;>\; \zeta - C\,
  \frac{\|\B\|_{\mathrm{op}}^2\,\mathcal{N}_2(\lambda)
  + \|\B\|_{\mathrm{op}}\sqrt{N\zeta}}{\sqrt{\delta}\,N},
\]
for a constant $C$ depending only on $c$ and $D$. Work on the event
$\mathcal{E}$ of~\eqref{eq:event}, which has probability at least
$1-\delta/2$ under the conditions of Lemma~\ref{lem:B-general}, so that
$\|\B\|_{\mathrm{op}} \leq C_u$. For
$N \geq 4C\log(1/\alpha')/\sqrt{\delta}$, the term
$C\log(1/\alpha')\,\zeta/(\sqrt{\delta}N)$
in~\eqref{eq:adapt-thresh} is at most $\zeta/4$; for
$N \geq \log^2(1/\alpha')$, the term
$C\log(1/\alpha')\|\B\|_{\mathrm{op}}\sqrt{\zeta}/(\sqrt{\delta}N)$ there is at
most the term $C\|\B\|_{\mathrm{op}}\sqrt{N\zeta}/(\sqrt{\delta}N)$ of the
second display; and, using
$\log(1/\alpha') \geq 1$ and $\delta \leq 1$, the remaining terms total at
most $\zeta/2$ as soon as
$\zeta \geq C'\log(1/\alpha')\,(\mathcal{N}_2(\lambda)\vee1)/(\delta N)$ for a
constant $C'$ depending only on $C$ and $C_u$. On the intersection of the
three events, which has probability at least $1-2\delta$, this gives
$\hat\eta > \hat q_{1-\alpha'}$. Since $\eta \leq (C_u^2/C_\ell^2)\,\zeta$ on
$\mathcal{E}$ by~\eqref{eq:zeta-lb}, the individual test therefore has power
at least $1-2\delta$ at level $\alpha'$ provided
\begin{equation}\label{eq:adapt-cond1}
  \eta \;\gtrsim\; \frac{\log(|\Lambda\times\mathcal{W}|/\alpha)\,
  \bigl(\mathcal{N}_2(\lambda)\vee1\bigr)}{\delta N}
\end{equation}
and the conditions of Lemma~\ref{lem:B-general} hold. Under the stated grid
endpoints, both $|\Lambda| = 1 + \log_2(\lambda_U/\lambda_L)$ and
$|\mathcal{W}| = 1 + \log_2(h_U/h_L)$ are of order $\log N$, so
$|\Lambda \times \mathcal{W}| \asymp \log^2 N$ and, using
$\log(x/\alpha) \leq 2\log(1/\alpha)\log x$ for $x \geq e$ and
$\alpha \leq e^{-1}$,
\[
  \log(|\Lambda\times\mathcal{W}|/\alpha) \;\lesssim\;
  \log(1/\alpha)\,\log\log N .
\]
In particular, the requirement on $N$ above holds for $N \geq N_0$.
Substituting this bound and $\mathcal{N}_2(\lambda)\vee1 \lesssim \lambda^{-1/2}$
for $\lambda \leq \kappa$ into~\eqref{eq:adapt-cond1} gives
\begin{equation}\label{eq:adapt-cond2}
  \eta \;\gtrsim\; \frac{\log\log N}{\delta N \sqrt{\lambda}},
\end{equation}
with $\log(1/\alpha)$ absorbed into the implicit constant.

\textbf{Step 2: the target pair.} Consider first
$(P,Q) \in \mathcal{F}_{s,\Delta}$ with $s > s_l$. Repeating the proof of
Theorem~\ref{thm:separation} with~\eqref{eq:adapt-cond2} in place of
condition~\eqref{cond:eta}, and retaining $\lambda \asymp h^d$,
$h \asymp \Delta^{1/(2s)}$, condition~\eqref{eq:adapt-cond2} reduces to
$\Delta^{1+d/(4s)} \gtrsim \log\log N/N$, while the requirement
$\lambda \gtrsim l^{-1}\log(l/\delta)$ of Lemma~\ref{lem:B-general} reduces,
for $l \asymp N$, to $\Delta \gtrsim (\log N/N)^{2s/d}$. The stated Sobolev
rate is therefore sufficient, at
\begin{equation}\label{eq:hstar}
  h^* \;\asymp\; \max\Big\{
  \big(\tfrac{\log\log N}{N}\big)^{2/(4s+d)},\;
  \big(\tfrac{\log N}{N}\big)^{1/d}\Big\},
  \qquad \lambda^* = c_*\,(h^*)^d ,
\end{equation}
where $c_* > 0$ is small enough that every $\lambda \leq c_*2^d h^d$ satisfies
the first condition in~\eqref{eq:resolution} with $t = h$ and the upper
constraint $\lambda \leq \|\Sigma_R\|_{\mathrm{op}}$ of
Lemma~\ref{lem:B-general}; this is possible by~\eqref{eq:mK-bandwidth} and
the assumed lower bound $\|\Sigma_R\|_{\mathrm{op}} \gtrsim h^d$. The
constants in the second terms of~\eqref{eq:hstar} do not depend on $s$, so
$h^* \geq h_L$ and $\lambda^* \geq \lambda_L$ for suitable $r_1, r_3$. For
$(P,Q) \in \mathcal{P}_{\theta,\Delta}$ with $\theta > \theta_l$, the same
argument applied to Theorem~\ref{thm:separation-rkhs} gives the stated
general-domain rate at
\[
  \lambda^* \;\asymp\; \max\Big\{
  \big(\tfrac{\log\log N}{N}\big)^{2/(4\theta+1)},\;
  \tfrac{\log N}{N}\Big\},
\]
with no constraint on $h$, and again $\lambda^* \geq \lambda_L$.

\textbf{Step 3: the grid attains the target.} Set
$v^*_\lambda \coloneqq \max\{x \in \Lambda : x \leq \lambda^*\}$ and
$v^*_h \coloneqq \max\{x \in \mathcal{W} : x \leq h^*\}$, which are well
defined by Step~2. Since both grids are geometric with ratio $2$,
$\min\{\lambda^*/2, \lambda_U\} \leq v^*_\lambda \leq \lambda^*$ and
$\min\{h^*/2, h_U\} \leq v^*_h \leq h^*$. Choose $r_2 \leq c_* r_4^d$ and
$r_2 \leq \inf\|\Sigma_R\|_{\mathrm{op}}$, the infimum being over the
general-domain class. We check that $(v^*_\lambda, v^*_h)$ satisfies the
conditions of Step~1 under the stated rates.

For the lower bound on $\eta$: since $v^*_h \leq h^*$ and
  $v^*_\lambda \leq \lambda^*$, the second condition in~\eqref{eq:resolution}
  and the condition of Proposition~\ref{prop:signal-rkhs} only become weaker. For
  the first condition in~\eqref{eq:resolution}: if $h^* \leq h_U$, then
  $v^*_\lambda \leq c_*(h^*)^d \leq c_*2^d(v^*_h)^d$; otherwise
  $v^*_h = h_U$ and $v^*_\lambda \leq \lambda_U \leq c_*h_U^d$. For the variance: if $\lambda^* \leq \lambda_U$, then
  $v^*_\lambda \geq \lambda^*/2$, and~\eqref{eq:adapt-cond2} holds at
  $v^*_\lambda$ under the stated rates, up to constants. Otherwise
  $v^*_\lambda = \lambda_U$, and the right-hand side
  of~\eqref{eq:adapt-cond2} is of order $\log\log N/N$, which is at most the
  first term of each stated rate. For Lemma~\ref{lem:B-general}, the lower constraint in~\eqref{eq:lambda-range} holds since
  $v^*_\lambda \geq \min\{\lambda^*/2, \lambda_U\} \gtrsim \log N/N$. The
  upper constraint holds for the Sobolev class since
  $v^*_\lambda \leq c_*2^d(v^*_h)^d$, by the choice of $c_*$, and for the
  general-domain class since $v^*_\lambda \leq \lambda_U \leq
  \inf\|\Sigma_R\|_{\mathrm{op}}$.

Hence the individual test at $(v^*_\lambda, v^*_h)$ has power at least
$1 - 2\delta$ at level $\alpha'$ under the stated rate, and since the union
test rejects whenever this test does,
\[
  P_{H_1}\Big\{\bigcup_{(\lambda,h)\in\Lambda\times\mathcal{W}}
  \hat\eta_{\lambda,h} > \hat q^{\,\lambda,h}_{1-\alpha'}\Big\}
  \;\geq\; P_{H_1}\big\{\hat\eta_{v^*_\lambda, v^*_h} >
  \hat q^{\,v^*_\lambda,v^*_h}_{1-\alpha'}\big\} \;\geq\; 1 - 2\delta .
\]
For $(P,Q) \in \mathcal{P}_{\theta,\Delta}$ any $h \in \mathcal{W}$ is
admissible, since the bandwidth does not enter
Proposition~\ref{prop:signal-rkhs}. Taking the infimum over the union of the two
classes gives the result; none of $r_1, \dots, r_4$ depends on $N$, $s$ or
$\theta$. \qed

\subsection{Proof of Corollary~\ref{cor:adaptation-fixed}}

The level follows exactly as in the proof of Theorem~\ref{thm:adaptation},
with $|\Lambda|$ in place of $|\Lambda\times\mathcal{W}|$, and Step~1 of that
proof gives~\eqref{eq:adapt-cond1} with the same replacement, so that power
at least $1-2\delta$ at the corrected level holds provided
$\eta \gtrsim \log\log N\,(\mathcal{N}_2(\lambda)\vee1)/(\delta N)$ and the
conditions of Lemma~\ref{lem:B-general} hold.

For $(P,Q) \in \mathcal{F}_{s,\Delta}$ with $s > s_l$, repeat the proof of
Theorem~\ref{thm:sobolev-fixed} with this condition. Since
$\mathcal{N}_2(\lambda)\vee1 \lesssim \lambda^{-d/(4\tau)}$ for
$\lambda \leq 1$ by~\eqref{eq:eigdecay}, it suffices that
\[
  \Delta \gtrsim \max\Big\{\lambda^{s/\tau},\;
  \frac{\log\log N}{N\,\lambda^{d/(4\tau)}}\Big\},
\]
and Lemma~\ref{lem:B-general} with $l \asymp N$ requires
$\lambda \gtrsim \log N/N$. Balancing the two terms subject to this
constraint gives the stated rate, at
\[
  \lambda^* \asymp \max\Big\{
  \big(\tfrac{\log\log N}{N}\big)^{4\tau/(4s+d)},\;
  \tfrac{\log N}{N}\Big\}.
\]
The constant in the second term does not depend on $s$, so
$\lambda^* \geq \lambda_L$ for suitable $r_1$. Set
$v^*_\lambda \coloneqq \max\{x \in \Lambda : x \leq \lambda^*\}$ and choose
$r_2 \leq 1 \wedge \inf\|\Sigma_R\|_{\mathrm{op}}$, the infimum being over both
classes. If $\lambda^* \leq \lambda_U$, then
$v^*_\lambda \in [\lambda^*/2, \lambda^*]$ and the conditions above hold at
$v^*_\lambda$ up to constants. Otherwise $v^*_\lambda = \lambda_U$: the first
term decreases, and the second becomes of order $\log\log N/N$, which is at
most the first term of the stated rate. In both cases
$\log N/N \lesssim v^*_\lambda \leq \|\Sigma_R\|_{\mathrm{op}}$, so
Lemma~\ref{lem:B-general} applies, and the union test rejects whenever the
test at $v^*_\lambda$ does. For $(P,Q) \in \mathcal{P}_{\theta,\Delta}$ the
argument is the general-domain part of the proof of
Theorem~\ref{thm:adaptation}, which uses $\Lambda$ alone. \qed

\bibliographystyle{plainnat}
\bibliography{ref}

\end{document}